\documentclass[11pt, oneside]{amsart}
\usepackage{makecell}
\usepackage{amsmath,ifthen, amsfonts, amssymb, 
srcltx, amsopn,tikz, tkz-euclide, url}
\usetikzlibrary{decorations.markings,arrows, cd}
\usepackage{tikz-cd}
 \usepackage{dutchcal}
\usepackage{xfrac} 
\usepackage{graphicx, ifsym}
\usepackage{mathtools}
\usepackage{bm,amsmath,amssymb,color,subfigure}
\usepackage{fancyhdr, ifpdf}
\usepackage[normalem]{ulem}

\newcommand{\showcomments}{yes}
\renewcommand{\showcomments}{no}
\usepackage{hyperref}
\usepackage[capitalise, noabbrev]{cleveref}

\newsavebox{\commentbox}
{\ifthenelse{\equal{\showcomments}{yes}}%
{\footnotemark
        \begin{lrbox}{\commentbox}
        \begin{minipage}[t]{1.25in}\raggedright\sffamily\tiny
        \footnotemark[\arabic{footnote}]}
{\begin{lrbox}{\commentbox}}}%
{\ifthenelse{\equal{\showcomments}{yes}}%
{\end{minipage}\end{lrbox}\marginpar{\usebox{\commentbox}}}
{\end{lrbox}}}

\usepackage{amsmath,amsfonts, amssymb, graphicx, tikz, bbding,ifthen,srcltx,amsopn, url}
\usepackage{bm}
\usetikzlibrary{decorations.markings,arrows, intersections}
\usetikzlibrary{decorations.pathreplacing}

\newcommand{\N}{\mathbb N}
\newcommand{\R}{\mathbb R}

\newcommand{\inv}[1]{#1^{-1}}

\newcommand{\embeds}{\hookrightarrow}

\newcommand{\mb}{\partial_*}

\newcommand{\Morsetree}{Morse hierarchy}
\newcommand{\incompleteMorsetree}{weak Morse hierarchy}

\newcommand{\stab}{\text{stab}}
\newcommand{\acts}{\curvearrowright}
\newcommand{\emb}{\hookrightarrow}
\newcommand{\mc}{\mathcal}
\newcommand{\mf}{\mathfrak}

\newcommand{\blowup}{\overline{\mc{H}}}

\theoremstyle{definition}
\newtheorem{thm}{Theorem}[section]
\newtheorem{lem}[thm]{Lemma}
\newtheorem{prop}[thm]{Proposition}

\newtheorem{question}{Question}

\newtheorem{remark}[thm]{Remark}
\newtheorem{defn}[thm]{Definition}
\newtheorem{example}[thm]{Example}
\newtheorem{cor}[thm]{Corollary}

\newtheorem{con}{Convention}

\newtheorem{introthm}{Theorem} 
\newtheorem{introcor}[introthm]{Corollary}

\newcommand{\g}{\gamma}
\newcommand{\cu}{\subseteq} 
\renewcommand{\H}{\mc{H}} 
\newcommand{\ra}{\rightarrow} 
\usepackage{dirtytalk} 
\usepackage[version=4]{mhchem}
\newcommand{\E}{\mc{E}} 

\makeatletter
\let\@wraptoccontribs\wraptoccontribs
\makeatother

\author{Matthew Cordes}
\author{Annette Karrer}
\author{Kim Ruane}
\address{Maxwell Institute for Mathematical Sciences and Department of Mathematics, Heriot-Watt University, Edinburgh, UK}
\email{m.cordes@hw.ac.uk}
\address{Institute of Mathematics, University of Wrocław, Wrocław, Poland}
\email{annette.karrer@uwr.edu.pl}
\address{Department of Mathematics, Tufts University, Medford, MA  02155, USA}
\email{kim.ruane@tufts.edu}
\contrib[with an appendix by]{Vivian He \and Annette Karrer}

\title{Morse accessibility with Applications to Coxeter groups}

\begin{document}

\begin{abstract}
In this paper we begin developing a theory of Morse hierarchies for finitely generated groups. Very roughly, a Morse hierarchy is a iterative splitting of a group over subgroups that do not contribute to the Morse boundary of the ambient group. 

We apply this theory to Coxeter groups, establishing new combinatorial criteria that distinguish certain  Coxeter groups from other Coxeter groups as well as right-angled Artin groups up to quasi-isometry. We highlight the consequences for the right-angled case. To build the reader's intuition, the last section provides many examples of Morse hierarchies and their subtleties. 

Finally, there is an appendix by the second author and He that extends these ideas to the metrizable Cashen--Mackay topology on the Morse boundary.
\end{abstract}

\maketitle

\section{Introduction}

\subsection{Hierarchies and Morse accessibility} 
Suppose a finitely generated group $\Gamma$ splits as a graph of groups whose vertex groups split again as a graph of groups and perhaps the vertex groups of \emph{that} splitting also split again, and so on. 
This type of iterated graphs of groups splitting leads to the notion of a \emph{hierarchy} of splittings for the group $\Gamma$. There are generally restrictions on the types of splittings allowed in this process. The goal is to obtain a finite hierarchy so that the vertex groups in the final splitting, do not split anymore in the allowable way. 

In the development of various types of hierarchies, \emph{boundaries} have proven to be a powerful tool. They frequently reveal splittings and hierarchies inherent in the group's structure. For instance, topological properties of different types of boundaries can reveal splittings and even full JSJ decompositions for different classes of groups, including hyperbolic groups \cite{Bowditch_canonical_splittings}, relatively hyperbolic groups \cite{Bowditch_peripheral_splitting} \cite{Haulmark_cutpoints}  \cite{Haulmark_Hruska_canonical_splittings}, and $\mathrm{CAT}(0)$ groups \cite{Papasoglu_Swenson}.
The most important example for this paper, however, is the Dunwoody hierarchy for hyperbolic groups: Every hyperbolic group admits a finite hierarchy over finite edge groups  where the final vertex groups are hyperbolic and either one-ended or finite \cite{Dunwoody}. In terms of the Gromov boundary, this hierarchy can be described as one in which the edge groups have empty Gromov boundary and where every vertex group has connected (possible empty) Gromov boundary. The Gromov boundaries of the one-ended vertex groups embed in the Gromov boundary of the ambient group and yield its infinite connected components. This is the origin of the idea in  \cite{karrer2024subgroups}  to investigate hierarchies of this kind not only for the Gromov boundary of a hyperbolic group but also for the \emph{Morse boundary} for a general finitely generated group.  

The Morse boundary was introduced by Charney--Sultan~\cite{Charney-Sultan}, and Cordes \cite{Cordes} with the intention of providing a tool analogous to the Gromov boundary of a hyperbolic group for groups that are not Gromov hyperbolic.  It is a generalization in the sense that the two boundaries coincide for hyperbolic groups and both are a quasi-isometry invariant for the group. 

When transferring the above from Gromov boundaries to Morse boundaries, the natural subgroups to split over are subgroups $H$ whose \emph{relative} Morse boundary $(\mb H, G)$ is empty, i.e., subgroups $H$ that do not have hyperbolic-like behavior when viewed in the ambient group. A tree that keeps track of such iterative splittings is called a \emph{weak Morse hierarchy}. One might hope that such a splitting process terminates at some point with final vertex groups that have connected, possibly empty (relative) Morse boundary. If it does, we call the tree a \emph{Morse hierarchy}. In \cite{karrer2024subgroups}, the authors began to investigate which finitely generated groups admit such a hierarchy, and called such groups \textit {Morse accessible}. In this  paper, we formalize the notion of Morse accessibility in detail, using and extending the results in \cite{fioravanti2022connected}. This allows us to verify that Morse hierarchies indeed generalize the Dunwoody hierarchy for hyperbolic groups in the following sense:

\begin{introthm}\label{thmMorseaccessiblegroups_Intro}
    Let $\Gamma$ be a Morse accessible group and $\mathcal H$ a \Morsetree~ of $\Gamma$. Then the following two sets are identical: 
    \begin{itemize}
        \item $\mathfrak{L}$, the set of leaf groups $L$ of $\mc H$ with $(\mb L, \Gamma)\neq \emptyset$
        \item $\mathfrak{S}$, the set of stabilizers of infinite connected components of $\mb \Gamma$. 
        \end{itemize}
        Moreover, every infinite connected component $C$ in $\mb \Gamma$ satisfies \[C = (\mb \stab (C), \Gamma).\]
  \end{introthm}

\begin{figure}[h]
    \centering
    \includegraphics[width=0.7\textwidth]{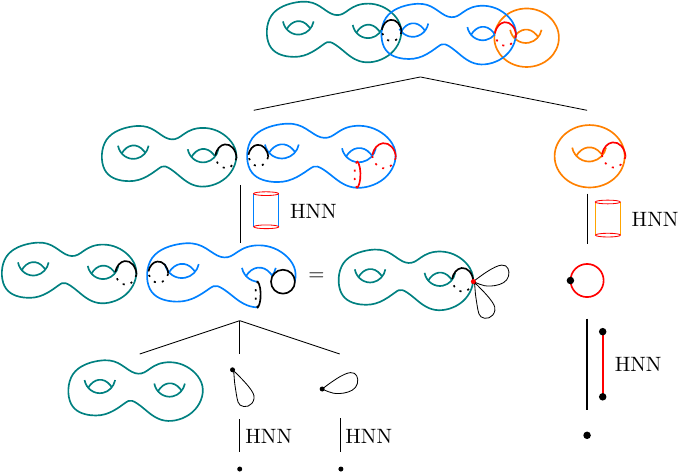}
    \caption{The pictured tree yields a Morse hierarchy for the fundamental group of the amalgam at the top. }
    \label{fig:surfaceamalgam1_intro}
\end{figure}

  For example, suppose $\Gamma$ is the fundamental group of the surface amalgam at the top of \cref{fig:surfaceamalgam1_intro}. One Morse hierarchy of $\Gamma$ over $\mathbb{Z}$-subgroups (with empty relative Morse boundary) can be obtained by decomposing the space successively 
  as illustrated in the finite tree of spaces below. 

There are four groups associated to the leaves of this tree, the fundamental groups of the genus two surface, the torus, and two one-point spaces. Importantly, the set $\mf L$ of leaf groups with non-empty relative Morse boundary consists only of the set of all conjugates of $\pi_1(S)$ where $S$ denotes the surface at the left bottom of \cref{fig:surfaceamalgam1_intro}. The theorem says that the set $\mf L$ coincides with the set of stabilizers of infinite Morse boundary components of $\Gamma$.  In particular, the set $\mf L$ is independent of the choice of the Morse hierarchy. A detailed explanation of this and other examples can be found in \cref{sec:examples}. 

  A direct consequence of the theorem above is that we can characterize when a Morse-accessible group has totally disconnected Morse boundary. This extends and uses results from \cite{fioravanti2022connected}.

\begin{introcor}\label{corMorseaccessiblegroups_Intro}
         Let $\Gamma$ be a Morse accessible group with a \Morsetree~ $\mathcal H$. The Morse boundary of $\Gamma$ is totally disconnected if and only if for every leaf group $L$ of $\mc H$ we have $(\mb L, \Gamma)=\emptyset$. 
    \end{introcor}

Just as with group splittings, it is natural to consider how a Morse hierarchy behaves with respect to quasi-isometry. For instance, the quasi-isometry invariance of accessibility follows from the characterization of accessibility in \cite{THOMASSEN1993248}. 
 We apply the results from \cite{karrer2024subgroups} to show the quasi-isometric rigidity statements in \cref{thm:point-convergence_intro} and \cref{corSurfaces_intro} below.
In order to obtain these results, we show that the set of infinite Morse boundary components is {point-convergent}.  \emph{Point-convergence} is a notion that was introduced in \cite{karrer2024subgroups}, and generalizes the notion of a \emph{null-family of subsets}. Null-families have been important in the study of boundaries. In particular, they were used to characterize quasi-convex subgroups of hyperbolic groups \cite{Gitik–Mitra–Rips–Sageev, Mj08}. 

\begin{introthm} \label{thmintro:point-convergende}
    The set of all infinite connected components of the Morse boundary of a Morse accessible group is point-convergent.
\end{introthm}

By Corollary 1.2, Corollary 1.3, and Theorem 1.5 in \cite{karrer2024subgroups}, \cref{thmintro:point-convergende} directly implies the following two statements. 

\begin{introcor}\label{thm:point-convergence_intro}
 Let $\Gamma$ be a Morse accessible group, and $\Gamma'$  a finitely generated group quasi-isometric to $\Gamma$. Then 
    
    \begin{itemize}
        \item if $C$ is an infinite Morse boundary component of $\mb \Gamma'$ with $H$ its stabilizer,  then $C$ coincides with both the limit set of $H$ in $\mb \Gamma'$ and the relative Morse boundary $(\mb H, \Gamma')$;
        \item an infinite component $C$ of $\mb \Gamma'$ is compact if and only if its stabilizer is stable;
        \item if  $C$ is a compact infinite connected component of $\mb \Gamma'=\mb \Gamma$, then the stabilizer of $C$ in $\Gamma$ and $\Gamma'$ are quasi-isometric to each other. 
    \end{itemize}
\end{introcor}

Applying the previous theorem to Morse accessible group whose Morse boundary contains a connected component homeomorphic to a circle yields:

\begin{introcor}\label{corSurfaces_intro}
 Let $\Gamma$ be a Morse accessible group whose Morse boundary $\mb \Gamma$ contains a connected component $C$ homeomorphic to a circle. Then any finitely generated group quasi-isometric to $\Gamma$ has a stable subgroup that contains as a subgroup of finite index the fundamental group of a closed hyperbolic surface. 
\end{introcor}

     Note, there are two natural topologies one can put on the Morse boundary. In this paper, we use the direct limit topology introduced in \cite{Charney-Sultan} and \cite{Cordes}. A metrizable topology was introduced in \cite{cashen2019metrizable} and we refer to this as the Cashen-Mackay topology. 
     Almost all Morse hierarchies discussed in this paper remain valid when considered with respect to the Cashen–Mackay topology. This is discussed in the Appendix \cref{appendix}.

\subsection{(Weak) Morse hierarchies of Coxeter groups and applications} 
\label{seubsec:introAppCox}

We show that every Coxeter groups has a weak Morse hierarchy. If the Coxeter group does not contains a \emph{spherical nexus} (a condition restricting the behavior of special spherical subgroups) this weak Morse hierarchy is a \emph{genuine} Morse hierarchy. Many Coxeter groups do not contain a spherical nexus, as for instance all right-angled Coxeter groups (RACGs), and are therefore Morse accessible.  We derive consequences for their Morse boundaries and quasi-isometry classification in \cref{subsec:introRACGRAAG}. For definitions related to Coxeter groups, see \cref{sec:MorseAccessAFCoxeter}.

To obtain a finite Morse hierarchies for Coxeter groups, we use ideas from \cite{karrer} to develop a graph decomposition algorithm over subgraphs that correspond to special subgroups having empty relative Morse boundary in the ambient group. See \cref{fig:proof_morse_accessibility} for how the algorithm works.  Then we apply the results from Cordes--Levcovitz~\cite{CorLev25} to show that the vertex groups we are left with one-ended and have connected relative Morse boundary (this is detected by the Coxeter graph being \emph{relatively wide-avoidant}). These vertex groups are special subgroups that are either finite or the stabilizer of an infinite connected component of the Morse boundary. Conversely, we can identify the stabilizer of each infinite connected component of the Morse boundary as a conjugate of one of these non-finite special subgroups that remain after applying the algorithm. This is a remarkable relation between the topology of the Morse boundary and the combinatorics of the graph defining the Coxeter group.

 \begin{introthm}   
\label{thmRACG_intro} (\cref{thmRACG_}) Every Coxeter group has a weak Morse hierarchy
   in which each leaf group $L$ is conjugate to a special subgroup $W_\Lambda$ that satisfies one of the following: 

   \begin{enumerate}
       \item  $\Lambda$ is a relatively wide-avoidant subgraph of $\Delta$  and $W_\Lambda$ is one-ended. In particular, $(\mb W_\Lambda, W_\Delta) \neq \emptyset$; 
       \label{property1}
       \item  $\Lambda$ is a clique that is either spherical or contained in a wide subgraph of $\Delta$. In particular, $(\mb W_\Lambda, W_\Delta) = \emptyset$. \label{property2}
   \end{enumerate}

Moreover, if $W_\Delta$ does not contain any spherical nexus then for the groups in (\ref{property1}), we additionally have that $(\mb W_\Lambda, W_\Delta)$ is connected and locally connected. This gives a Morse hierarchy for Coxeter groups that do not contain
any spherical nexus.
\end{introthm}

\begin{figure}[h]
    \centering
    \includegraphics[width=0.3\textwidth]{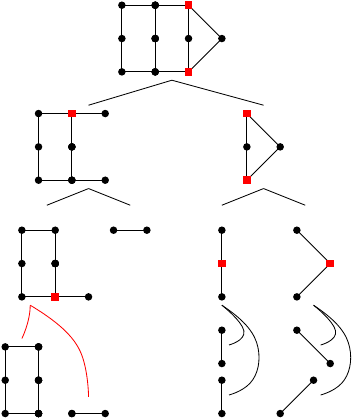}
    \caption{The pictured tree yields a Morse hierarchy of the RACG with the defining graph at the top. At each level we split over the subgroup generated by the red vertices.} 
    \label{fig:Morse_hierarchie_RACGs}
\end{figure}

We illustrate the process for obtaining a finite Morse hierarchy in the RACG case in \cref{fig:Morse_hierarchie_RACGs}. In this case, we iteratively split the defining graph along cliques and subgraphs contained in joins which represent a direct products of two infinite subgroups. For more details on this example, see \cref{subsec:exampes_RACG}.

A direct consequence of \cref{thmMorseaccessiblegroups_Intro} and \cref{thmRACG_intro} is the following.

\begin{introcor} \label{stabCor_intro} 
Every infinite connected component $C$ of the Morse boundary of a Coxeter group $W_\Delta$  that does not contain a spherical nexus satisfies:
    \begin{itemize}
        \item $C$ is locally connected, 
        \item $C$ equals the relative Morse boundary of its stabilizer
        \item has a stabilizer isomorphic to a one-ended special subgroup whose Coxeter graph is relatively wide-spherical-avoidant in $\Delta$.
    \end{itemize}
\end{introcor} 

Again, we illustrate this corollary with the RACG example in \cref{fig:Morse_hierarchie_RACGs}. The leaf corresponding to the $6$-cycle (on the bottom left of the figure) is the defining graph of a special subgroup $H$ that is the stabilizer of an infinite Morse boundary component $C$ homeomorphic to its Gromov boundary, $\partial H \cong \mathbb{S}^1$. It is possible, though, that a stabilizer $H$ is hyperbolic, but its corresponding connected component is not homeomorphic to $\partial H$, see \cref{ex:topologies different}.

As a direct consequence of \cref{stabCor_intro} (and implicitly \cref{corMorseaccessiblegroups_Intro}), we generalize the results in \cite{karrer} for RACGs as follows: 
\begin{introcor}\label{thm:AFCoxeterTotDisconn_intro}
    If every leaf group of the weak Morse hierarchy of $W_\Delta$ constructed in the proof of \cref{thmRACG_intro}  (\cref{thmRACG_}) is spherical or contained in a special subgroup of $W_\Delta$ that is wide, then $\mb W_\Delta$ is totally disconnected. Moreover if $W_\Delta$ does not contain any spherical nexus, then the reverse implication is also true. 
\end{introcor}

In the \cref{subsec:introRACGRAAG}, we will apply \cref{thm:AFCoxeterTotDisconn_intro} and the proof of \cref{thmRACG_intro} to show that there is an algorithm to determine when certain Coxeter groups are not quasi-isometric to any right-angled Artin group (RAAG). This enables us to differentiate many new RACGs groups from each other up to quasi-isometry, even in the non triangle-free case. 

\begin{figure}[h]
    \centering
    \includegraphics[]{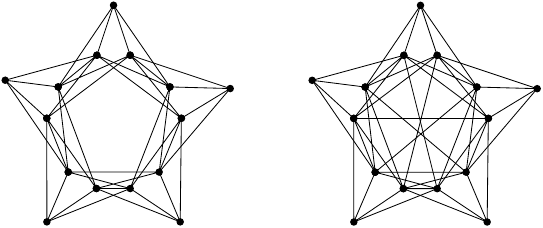}
    \caption{The two RACGs corresponding to the pictured graphs are not quasi-isometric to each other. The graph on the left was studied in \cite[Fig 10]{MihalikCamp}. }
    \label{fig:Main_exampleintro}
\end{figure}

 For instance, the graph on the left in \cref{fig:Main_exampleintro} corresponds to a RACG with  totally disconnected Morse boundary  \cite{karrer}. We show in \cref{subsec:introRACGRAAG} that the RACG corresponding to the graph on the right of \cref{fig:Main_exampleintro} has an infinite connected component in its Morse boundary. So, the pictured graphs correspond to two RACGs whose Morse boundaries are not homeomorphic. Hence, the corresponding RACGs are not quasi-isometric.

For more examples, see \cref{fig:Main_examples} in \cref{subsec:introRACGRAAG}.

\subsection{Outline}
The paper is organized as follows.  \cref{sec:prelims} contains the necessary definitions and properties of Morse boundaries, relative Morse boundaries, and graphs of groups. Additionally, the theorems of \cite{fioravanti2022connected} and \cite{karrer2024subgroups} are recalled, and we provide slight variants of Theorem 2.1 from the former and Lemma 5.8 from the latter, which are needed in this paper. (\cref{thmAvariant}). 

In \cref{sec:Morse-accessibility} we give the formal definition of Morse accessibility which is a modified version of the notion introduced in \cite{karrer2024subgroups}. 
 In order to formulate our definition, we introduce \textit{Morse hierarchies} and methodically transfer all the results proven for group splittings along groups with empty relative Morse boundaries from \cite{fioravanti2022connected} to Morse hierarchies using induction. This will give \cref{thmMorseaccessiblegroups_Intro} and \cref{corMorseaccessiblegroups_Intro}, which provide structural consequences of a group being Morse accessible. 
 Additionally, \cref{lem:Cardrelmorse} shows that the relative Morse boundary of a finitely generated undistorted subgroup of a finitely generated group cannot consist of just one point. Finally, in \cref{sec:rigidity}, we prove \cref{thmintro:point-convergende} which is necessary for obtaining the quasi-isometric rigidity statements in 
 \cref{thm:point-convergence_intro} and \cref{corSurfaces_intro}.
 
 In \cref{sec:MorseAccessAFCoxeter} we recall and adjust the relevant results from Cordes-Levcovitz \cite{CorLev25} and use them to prove \cref{thmRACG_intro}, \cref{stabCor_intro}, and \cref{thm:AFCoxeterTotDisconn_intro}.

In \cref{subsec:introRACGRAAG} we illustrate how to use Morse hierarchies to distinguish groups up to quasi-isometry, with a focus on right-angled Artin groups (RAAGs) and right-angled Coxeter groups (RACGs).

In \cref{sec:examples}, we discuss how to obtain Morse hierarchies for several classes of groups including RAAGs, 3-manifold groups, surface amalgams.  We give several examples of Morse hierarchies for RACGs that illustrate properties and consequences of having a Morse hierarchy. 

Finally, in \cref{sec:AlternativeMorse} we consider alternative definitions of Morse accessibility. 

\subsection*{AI Declaration}
Artificial intelligence tools were used solely for grammatical editing and, in one instance, to assist with nomenclature. No AI tools, including LLMs, were used to develop, verify, or write any of the mathematical arguments presented in this article.

\subsection*{Acknowledgments}
The authors would like to thank Emily Stark and Motiejus Valiunas for fruitful discussions. 
\section{Preliminary definitions and results}\label{sec:prelims}

In this section, we will give the necessary definitions and background on Morse boundaries and relative boundaries for metric spaces and groups. We will also set up our notation and terminology for graphs of groups and recall the work in \cite{fioravanti2022connected} needed here.  We provide variants of the results in that paper that will allow us to work with hierarchies rather than just graphs of groups. 

\subsection{Morse Boundaries and  Relative Morse Boundaries of Metric Spaces}\label{subsec:metric}

In this subsection, we give the definition of the Morse boundary.
There are two natural topologies one can consider on the Morse boundary. The first is the direct-limit topology, which we use throughout this paper. An alternative topology is the metrizable topology introduced by Cashen and Mackay \cite{cashen2019metrizable}. 
The main advantages of the direct limit topology is its 
ease of computation. However, it has some notable drawbacks: It is not metrizable, and does not interact well with subgroups. We address which results carry over to the Cashen--Mackay topology in \cref{appendix}.

\begin{defn} 
Let $N \colon [1,\infty) \times [0,\infty) \to [0,\infty)$ be a function which is non-decreasing and continuous in the second factor. The quasi-geodesic $\gamma \colon I \to X$ 
is an  \emph{$N$-Morse quasi-geodesic} if for any sub-interval $[s,t] \subseteq I$ of the interval $I$, any $(k,c)$-quasi-geodesic with endpoints $\gamma(s)$ and $\gamma(t)$ is contained in the $N(k,c)$-neighborhood of the image of $\gamma\vert_{[s,t]}$.
The function $N$ is called a \emph{Morse gauge}. 
\end{defn}

As a set, the Morse boundary of a proper geodesic metric space $X$ is the set of all Morse geodesic rays in $X$ based at a point 
$o\in X$ and modulo asymptotic equivalence; given a geodesic ray $\gamma \colon [0,\infty)\to X$, 
we denote the equivalence class as $[\alpha]$. 
This collection is called the \emph{Morse boundary} of $X$ and is denoted by $\partial_{*} X$ or $\partial_{*} X_o$ 
if we wish to emphasize the basepoint $o$. 

To topologize the boundary, first fix a Morse gauge $N$ and consider the subset of the Morse boundary that consists of all rays in $X$ 
with Morse gauge at most $N$:  
\begin{equation*} \partial_{*}^N X_o= \{[\alpha] \mid \exists \beta \in [\alpha] 
    \text{ that is an $N$-Morse geodesic ray with base point }o\}. \end{equation*} 
One can topologize this set with the compact-open topology. 
This topology is equivalent to the one defined by a neighborhood basis: 
at each point $[\alpha] \in \partial_{*}^N X_o$, and every $n \in \mathbb{N}$, let $V_n^N( \alpha)$ 
denote the set of $N$-Morse geodesic rays $\gamma \colon [0,\infty) \to X$ with $\gamma(0)=o$ and $d(\alpha(t), \gamma(t))< \delta_N$ 
for all $t<n$, where $\delta_N$ is a constant that depends only on $N$. 
A neighborhood basis of $[\alpha]$ in $\partial_{*}^N X_o$
is the set $$\{V_n^N(\alpha) \mid n \in \mathbb{N} \}.$$
The constant $\delta_N$ is computed in \cite{Cordes, corrigendum} and is chosen so that any ray in the 
equivalence class of $\alpha$ is included in $V_n^N(\alpha)$, i.e., so that the topology is well defined.

Let $\mathfrak M$ be the set of all Morse gauges. 
We put a partial ordering on $\mathfrak M$ so that for two Morse gauges $N, N' \in \mathfrak M$, we have that $N \leq N'$ if and only if 
$N(k,c) \leq N'(k,c)$ for all $k \geq 1$ and $c \geq 0$. We define the Morse boundary of $X$ to be
 \begin{equation*} \partial_{*} X_o=\varinjlim_\mathfrak{M} \partial^N_* X_o \end{equation*} with the induced direct limit topology, i.e., a set $U$ is open in $\partial_{*} X_o$ if and only if $U \cap \partial^N_* X_o$ is open for all $M \in \mathfrak{M}$.  For more details on the Morse boundary see \cite{Cordes, corrigendum}. 

The first salient feature of the Morse boundary is that it is invariant under quasi-isometry regardless of the geometry of the space. 

\begin{thm}[\cite{Cordes}]\label{thm:Morse_boundar_QI_invariant}
Let $X$ and $Y$ be proper geodesic metric spaces. Every quasi-isometry $f\colon X \to Y$ induces a homeomorphism $\partial f \colon \partial_{*} X_o \to \partial_{*} Y_{f(o)}$.
\end{thm}

As a corollary of this theorem, a change in basepoint induces a homeomorphism on the boundary. This justifies the omission of the basepoint from the notation. More generally, we will usually assume the basepoint is fixed and omit it from the notation $\partial^N_* X_o$ as well. The reader is warned, however, that unlike $\partial_* X$, the substrata, $\partial^N_* X$, \emph{do} depend on a choice of basepoint.

\begin{defn}[relative Morse boundary, see \cite{karrer2024subgroups}]\label{definition:relative_morse_boundary}
Let $Y$ be a subset of a proper metric space $X$. 
The \textit{relative Morse boundary} of $Y$ (with respect to $X$), denoted by $(\mb Y, X)$, is the set 
\[(\mb Y, X)=\{[\gamma]\in \mb X \mid  \exists ~ k \text{ so that } \gamma \subseteq  N_k(Y)\}\]   where $ N_k(Y)$ denotes the $k$-neighborhood of $Y$. We endow $(\mb Y, X)\subseteq \mb X$ with the subspace topology. 
\end{defn}

To understand how a subset $Y$ of $X$ determines a subset of $\mb X$, one might also consider the relationship between $\mb Y$ and $\mb X$.  In general, $\mb Y$ does not embed in $\mb X$. Indeed, there could be Morse rays in $Y$ that are no longer Morse when viewed inside $X$.  
This motivates the following definition of an important subspace of $\partial_* Y$.

\begin{defn}[surviving Morse boundary] \label{def:surviving_morse_boundary}
Let $\iota$ denote the inclusion map of $Y$ into $X$. We define the \emph{surviving Morse boundary of $Y$ with respect to the embedding in $X$} to be the subspace of $\mb Y$ containing those Morse rays from $Y$ that remain Morse when included in $X$. We denote this as follows:  \[\mathcal M(Y,X) \coloneqq \{[\gamma] \in \mb Y \mid \iota(\gamma) \text{ is Morse in }X\} \subseteq \mb Y. \]
\end{defn}

\begin{remark} \label{rem:differentConventions} 
A warning to the reader, in \cite{fioravanti2022connected}, the notation $(\mb Y, X)$ refers to the surviving Morse boundary, $\mathcal M(Y,X)$, of $Y$ with respect to its embedding in $X$. 
\end{remark}
\begin{remark}[Choice of topology on subspace]\label{rem_sub_top}
Note that these two definitions of topologies on the boundary of a subset $Y$ of $X$ are not the same. An example of this difference is demonstrated in \cite{fioravanti2022connected}[Example on page 341]. Resolving this discrepancy requires care and is where much of the technical difficulties lie in the proofs in this paper. 
\end{remark}

\begin{defn}\label{defn:undistorted}
Let $X$ be a metric space.  A subset $Y$ of $X$ is \emph{undistorted} in $X$ if the inclusion of $Y$ into $X$ is a quasi-isometric embedding. 
\end{defn}
If $[\gamma]\in \mathcal M(Y,X)$ then $\iota([\gamma])$ is at finite Hausdorff distance from a Morse geodesic ray and hence defines an injection $\iota_*: \mathcal M(Y,X) \to \mb X$.
The following lemma says that $(\mb \iota(Y), X)$ is exactly the image of the set $\mathcal M(Y,X)$ under $\iota_*$.

\begin{lem}[see \cite{karrer2024subgroups}]\label{lem:RelMorse}
    Let $(X, d_X)$ and $(Y, d_Y)$ be proper geodesic metric spaces and let $\iota \colon Y\to X$ be a quasi-isometric embedding. Let $\gamma \colon [0, \infty) \to X$ be a Morse geodesic ray. The following are equivalent: 
    \begin{enumerate}
        \item $[\gamma]\in (\mb \iota(Y), X)$
        \item there exists a geodesic ray $\lambda\colon [0, \infty)\to Y$ which is Morse in $Y$ and whose image $\iota\circ\lambda$ is at bounded Hausdorff distance from $\gamma$
    \end{enumerate}
\end{lem}
This combined with the observations in \cite{karrer2024subgroups} yields: 

\begin{lem}[Lemma 2.3 in \cite{fioravanti2022connected}]\label{lem_jota}
     Let $(X, d_X)$ and $(Y, d_Y)$ be proper geodesic metric spaces and let $\iota \colon Y\to X$ be a quasi-isometric embedding. Then $\iota$ induces a bijection $\iota_* \colon \mathcal M(Y,X) \to  (\mb \iota(Y), X)$ which is closed and open. 
\end{lem}

We now discuss the case where $X$ is a finitely generated group $\Gamma$ endowed with a word metric and $Y$ is an undistorted subgroup $H$ of $\Gamma$. In this setting, $\iota_*$ is $H$-equivariant and if $H$ is stable or Morse (see below), then $\iota_*$ is not only a closed and open bijection, it is a homeomorphism onto its image in $\mb\Gamma$ \cite[Thm 1.4]{Tran}. 

\begin{defn}[stable subgroup] \label{defstable}
Let $H<\Gamma$ are finitely generated groups, and further assume $H$ is an undistorted, quasi-convex subgroup, we say $H$ is a \emph{stable subgroup} of $G$ if there exists a Morse gauge $N$ such that every pair of points in $H$ can be connected by an $N$-Morse quasi-geodesic in $G$. 
\end{defn}
Tran shows that an infinite subgroup of a finitely generated group $H$ is stable in $\Gamma$ if and only if $H$ is hyperbolic and Morse in $\Gamma$ \cite[Thm 4.8]{Tran}. 

\begin{defn}[Morse subgroup \cite{Tran}, \cite{Genevois}] \label{def:MorseSubgroup}
Let $\Gamma$ be a finitely generated group and $H$ a subgroup of $\Gamma$. We say that $H$ is \textit{Morse} in $\Gamma$ if for every $K\ge 1$, $C\ge 0$, there is some $M = M(K,C)$ such that every $(K,C)$-quasi-geodesic in $\Gamma$ with endpoints on $H$ is contained in the $M$-neighborhood of $H$. 
\end{defn}

We end this subsection with three lemmas about how the surviving and relative Morse boundaries of subgroups are related to the Morse boundary of the ambient group. These will be used in order to extend the results on graphs of groups discussed in the next subsection to hierarchies in \cref{sec:Morse-accessibility}. 
 
\begin{lem}\label{lemma:iota_compatible} Suppose $\Gamma$ is a finitely generated group with finitely generated subgroups $\Gamma''\le \Gamma'$.  Let $\bar \iota: \Gamma'' \embeds \Gamma'$, $\iota:\Gamma' \embeds \Gamma$,  and $\iota\circ \bar \iota :\Gamma'' \embeds \Gamma$be quasi-isometric embeddings inducing the maps 
\begin{align*}
\bar \iota_*\colon& \mathcal{M}(\Gamma'', \Gamma') \to (\mb \Gamma'', \Gamma') \subseteq \mb \Gamma'\\
\iota_*\colon&\mathcal{M}(\Gamma',\Gamma) \to (\mb \Gamma', \Gamma) \subseteq \mb \Gamma \\
(\iota \circ \bar \iota)_*\colon&\mathcal{M}(\Gamma'',\Gamma) \to (\mb \Gamma'', \Gamma) \subseteq \mb \Gamma.
\end{align*}
Then $\iota_*((\mb \Gamma'', \Gamma') \cap \mc M(\Gamma', \Gamma))=(\mb \Gamma'', \Gamma)$. 
 \end{lem} 

 \begin{proof}
     Let $\xi\in  \iota_*((\mb \Gamma'', \Gamma') \cap \mc M(\Gamma', \Gamma))$. Then there exists $[\alpha] \in (\mb \Gamma'', \Gamma') \cap \mc M(\Gamma', \Gamma)$ such that $\iota_*([\alpha])=\xi$. Since $[\alpha] \in \mc M(\Gamma', \Gamma)$, $[\alpha]\in \mb \Gamma'$ such that $\iota(\alpha)$ is Morse in $\Gamma$. Since  $[\alpha] \in \mc (\mb \Gamma'', \Gamma')$, there exists $[\beta]\in \mb \Gamma''$ such that $\bar \iota_*(\beta)$ is Morse in $\Gamma'$ and at bounded Hausdorff distance to $\alpha$. In particular, $\iota(\bar \iota(\beta))$ is a Morse quasi-geodesic in $\Gamma$ at bounded Hausdorff distance to $\iota(\alpha)$.  Hence $\xi=\iota_*([\alpha]) \in (\mb \Gamma'', \Gamma)$.\\
     On the other hand, let $[\gamma]\in (\mb \Gamma'', \Gamma)$. We aim to find a Morse geodesic ray $\gamma'$ of $\mb \Gamma'$ with $[\gamma']\in (\mb \Gamma'', \Gamma') \cap \mc M(\Gamma', \Gamma)$ such that $\iota_*([\gamma'])= [\gamma]$. \\
     Sine $[\gamma] \in (\mb \Gamma'', \Gamma)$, there exists a Morse ray $\alpha$ in $\Gamma''$ such that $\iota (\bar \iota(\alpha))$ is Morse in $\Gamma$ and at bounded Hausdorff distance to  $\gamma$. Since $\iota$ and $\bar\iota$ are quasi-isometries, this implies that $\bar \iota(\alpha)$ is a Morse quasi-geodesic in $\Gamma'$. Let $\gamma'$ be a Morse geodesic ray starting at the basepoint of $\Gamma'$ that is at bounded Hausdorff distance to $\bar \iota(\alpha)$. Then $[\gamma'] = \bar\iota_*([\alpha]) \in (\mb \Gamma'', \Gamma')$.   As $\iota$ is a quasi-isometry,  $\iota(\gamma')$ is a Morse quasi-geodesic in $\Gamma$, Hence, $[\gamma'] \in \mc M(\Gamma', \Gamma)$.  Moreover, $\iota(\gamma')$ is bounded Hausdorff distance to $\iota (\bar \iota(\alpha))$. Hence, $\iota_*([\gamma']) = [\gamma]$. 
 \end{proof}
\begin{lem}\label{lemma:transitive} Suppose $\Gamma$ is a finitely generated group with finitely generated subgroups $\Gamma''\le \Gamma'$.  If $\Gamma''$ is undistorted in $\Gamma'$ and $\Gamma'$ is undistorted in $\Gamma$, then $\Gamma''$ is undistorted in $\Gamma$ and $(\mb \Gamma'', \Gamma)\subseteq (\mb \Gamma', \Gamma)$ in $\mb \Gamma$.
 \end{lem}

 \begin{proof}
     The composition of quasi-isometric embeddings is a quasi-isometric embedding so $\Gamma''$ is undistorted in $\Gamma$.  It remains to show that $(\mb \Gamma'', \Gamma)\subseteq (\mb \Gamma', \Gamma)$. 

     To show this, we assume $\gamma\colon[0,\infty)\to \Gamma$ is a Morse geodesic ray in $\Gamma$ such that $[\gamma]$ is in $(\mb \Gamma'',\Gamma)$.  Since $\Gamma''$ is undistorted in $\Gamma$, \cref{lem:RelMorse} implies there exists a geodesic ray $\lambda\colon[0,\infty)\to \Gamma''$ that is Morse in $\Gamma''$ whose image in $\Gamma$ is bounded Hausdorff distance (say $D$) from $\gamma$.  

    We must show $[\gamma]$ is also in $(\mb \Gamma', \Gamma)$.  As $\gamma$ is Morse in $\Gamma$, to show this, we only need to show that there exists a $k$ such that $\gamma\subseteq N_k(\Gamma')$ according to the usual definition of $(\mb \Gamma',\Gamma)$. 

    Each point of $\gamma$ is within $D$ of a point of $\Gamma''$. But since $\Gamma''\subseteq \Gamma'$, each point of $\gamma$ is also within $D$ of $\Gamma'$ as needed.
 \end{proof}

Since $\Gamma$ acts on its Cayley graph by isometries, it follows that the action of $\Gamma$ on the relative Morse boundary of any subgroup $H$, commutes with the action of $\Gamma$ on H by left multiplication and conjugation:

\begin{lem}\label{lem:gC}
 If $H$ is a subgroup of a finitely generated group $\Gamma$, then every $ g \in \Gamma$ satisfies: 
\[  g(\mb H, \Gamma)=(\mb gH, \Gamma)=(\mb gH\inv{g}, \Gamma). \]
\end{lem}

\subsection{Preliminaries on graph of groups}\label{sec:graphsofgroups}

\begin{con}\label{conv}
      We use the terminology and notation for graph of groups as in \cite{fioravanti2022connected}. In particular, we use the following characterization of group splittings:  A finitely generated group $\Gamma$ splits as a graph of groups if $\Gamma$ admits a non-elliptic action on a simplicial tree $T$ without inversions. In this case, there is a unique smallest $G$-invariant subtree. Restricting to it, we assume that the action of $G$ on $T$ is minimal, i.e. that no proper subtree is $G$-invariant.   We call $T$ a \textit{Bass--Serre tree} for the splitting. As in  \cite{fioravanti2022connected}, we refer to the vertex stabilizers  for the action $G \acts T$ as \textit{vertex groups}.   Analogously, we refer to the edge stabilizers for the action $G \acts T$ as \textit{edge groups}. 

    The quotient of $T$ is a finite simplicial graph. Labeling each vertex and edge of this graph with the stabilizer of one of its preimages as in \cite[Sec 5.4]{Serre} yields a finite graph of groups which we denote by $\mc G (\Gamma \backslash T)$. If the action on the tree $T$ is clear, we just write $\mc G (\Gamma)$. In this situation, we say that $\Gamma$ splits as graph of groups $\mc G(\Gamma)$. 

    Alternatively, given a finite graph of groups $\mc G (\Gamma)$, we can build a simplicial tree $T$ on which $\Gamma$ admits a non-elliptic action without inversions. For more information on this construction see \cite[Sec 5.3]{Serre}. We will use this viewpoint in Section \ref{sec:Morse-accessibility} and assume the action of $\Gamma$ onto this tree to be minimal. 
\end{con}

\subsection{Modified versions of \cite{fioravanti2022connected} and \cite{karrer2024subgroups}}
\label{sec:recallthms}

In this subsection, we recall the theorems of \cite{fioravanti2022connected} and \cite{karrer2024subgroups} and provide modified versions of some the results that are necessary for the proofs in \cref{sec:Def_Morse_hierarchy}.

\begin{thm}[Theorem A in \cite{fioravanti2022connected}]\label{thmA}
    Suppose the action of a finitely generated group $\Gamma$ on a Bass--Serre tree $T$ (as in \cref{conv}) results in a graph of groups splitting with the following properties:
    \begin{enumerate}
        \item[(a)] Each edge group $E$ is finitely generated and undistorted in $\Gamma$. 
        \item[(b)] $(\mb E, \Gamma)=\emptyset$ for every edge group $E$. 
    \end{enumerate}
     Let $x$ and $y$ be distinct vertices of $T$ and $V_x$ and $V_y$ be their stabilizers. 
     Then $(\mb V_x, \Gamma)\cap (\mb V_y, \Gamma)=\emptyset$. In particular, if $V_x$ and $V_y$ are the same subgroup of $\Gamma$, then that subgroup must have empty relative Morse boundary in $\Gamma$.
     
     Furthermore, if a connected component $C \subseteq \mb \Gamma$ contains at least two points, then there exists a vertex group $V$ so that  $C$ is contained in $(\mb V, \Gamma)$.  
\end{thm}

\begin{cor} \label{cor:gC}
Let $\Gamma$ be as in \cref{thmA}. 
    If a connected component $C \subseteq \mb \Gamma$ with $|C|>2$ is contained in $(\mb V, \Gamma)$ for a vertex group $V \le \Gamma$, then $\stab(C) \subseteq V$ where $\stab (C)$ denotes the stabilizer of $C$ in the action $\Gamma \acts \mb \Gamma$.
\end{cor}
\begin{proof}
    Let $g \in \Gamma$ so that $gC=C$. 
    By \cref{thmA}, there exists a vertex group $V$ so that $C \subseteq (\mb V, \Gamma)$. 
    Accordingly, $gC \subseteq g(\mb V, \Gamma)$. By \cref{lem:gC}, $g(\mb V, \Gamma)=(\mb gV\inv{g}, \Gamma)$. 
    Hence, $gC \subseteq (\mb gV\inv{g}, \Gamma)$. As $C = gC$, we conclude that $C \subseteq (\mb V, \Gamma) \cap  (\mb gV\inv{g}, \Gamma)$. However, $gV\inv{g}$ is a vertex group and hence, \cref{thmA} implies that $V = gV\inv{g}$, i.e. $g \in V$. 
    \end{proof}

Motivated by the \cref{thmA} and \cref{cor:gC} above, we give the following definition that will be generalized to hierarchies in \cref{sec:Morse-accessibility}.

\begin{defn}\label{def:weak-Morse-splittings}
       Let $\Gamma$ be a finitely generated group that splits as graph of groups $\mc G (\Gamma)$ with the following properties:
    \begin{enumerate}
        \item[(a)] Each edge group $E$ is finitely generated and undistorted in $\Gamma$. 
        \item[(b)] $(\mb E, \Gamma)=\emptyset$ for every edge group $E$. 
    \end{enumerate}
    Then we call the graph of groups $\mc G (\Gamma)$ (and splitting) \emph{weak Morse}. If in addition, we have that $(\mb V, \Gamma)$ is connected for every vertex group $V$, then we call the graph of groups (and splitting) \emph{Morse}.
\end{defn}

When generalizing to hierarchies in \cref{sec:Morse-accessibility}, we need a variant of \cref{thmA}. In \cref{thmA}, no edge group contains a Morse ray. In the variant presented here, we skip this condition and instead study subspaces of $\mb \Gamma$ that do not contain Morse rays in edge groups. The intuition behind this approach is that removing all equivalence classes of geodesic rays from edge groups effectively has the same consequence as starting with a group that does not contain Morse rays in the first place. Indeed, the proof of the variant discussed here can be obtained from the proof in \cite{fioravanti2022connected} with minor adjustments. 

\begin{thm}[Variant of Theorem A in \cite{fioravanti2022connected}]\label{thmAvariant}
    Suppose the action of a finitely generated group $\Gamma$ on a Bass--Serre tree $T$ (as in \cref{conv}) results in graph of group splitting with all edge groups finitely generated and undistorted in $\Gamma$.
    Let $ \mathcal A$ be the subspace of 
    \[\mb \Gamma - \left (\bigcup_{E \text{ edge group}} (\mb E, \Gamma)\right )\]
    Let $x$ and $y$ be distinct vertices of $T$ and $V_x$ and $V_y$ be their stabilizers. 
    Then $((\mb V_x, \Gamma)\cap \mathcal A)\cap ((\mb V_y, \Gamma) \cap \mathcal A)=\emptyset $. Furthermore, if a connected component $C \subseteq  \mathcal A$ contains at least two points, then there exists a vertex group $V$ so that  $C$ is contained in $(\mb V, \Gamma)\cap \mathcal A$;
\end{thm}
\begin{proof}
Adjust the proof in Section 3 in \cite{fioravanti2022connected} as follows: 
\begin{enumerate}
    \item In Assumption 3.1, remove the condition that the relative Morse boundary of every edge group is empty; 
    \item Note that the proofs of Lemmas 3.2 and 3.3 do not need the assumption that the edge groups have empty relative Morse boundary;
    \item After the proof of Lemma 3.3, whenever the Morse boundary of the considered group is mentioned, replace it with $\mathcal A$. In particular, whenever $\partial_M G$ occurs in the text, replace it by  $\mathcal A$; 
    \item By \cref{lem:RelMorse}, a Morse geodesic ray is of finite distance to an edge group $E$ if and only if its equivalence class is contained in the relative Morse boundary (in the sense of \cite{fioravanti2022connected}) of $E$;
    \item Whenever the assumption that the relative Morse boundaries (in the sense of \cite{fioravanti2022connected}) of the edge groups are empty, we instead use that $\mathcal A$ does not contain any equivalence class of a geodesic ray of finite distance of an edge group. This is in particular important for the proof of Corollary 3.4;
    \item Whenever a statement is made about the relative Morse boundary of a subgroup, delete from this set all equivalence classes of Morse geodesic rays that are at finite distance of an edge group. 
   This is in particular important for Lemma 3.9 and Corollary 3.10. 
\end{enumerate}
After these adjustments, the proof in Section 3 in \cite{fioravanti2022connected} can be copied word by word.
The notation for the relative Morse boundary in \cite{fioravanti2022connected} is different from the notation in this paper (see \cref{rem:differentConventions}).
Thus, \cref{thmAvariant} follows after transferring the notation accordingly.  
\end{proof}

We will also need the following variant of Lemma 6.7 in \cite{karrer2024subgroups}.

If $\gamma: \R \to X$ denotes a Morse bi-infinite line in a proper metric space $X$, $\gamma[0,\infty)$ and $\gamma(-\infty, 0]$ are the images of two Morse geodesic rays $\gamma^+$ and $\gamma^-$ of $X$. We denote the corresponding boundary points in $\mb X$ by $[\gamma^-]$ and $[\gamma^+]$.

\begin{lem} \label{lem:point-convergence_prem}
 Suppose a finitely generated group $\Gamma$ splits as graph of groups with all edge groups finitely generated and undistorted in $\Gamma$. 

    Let $ \mathcal A$ be the subspace of $\mb\Gamma$ with the relative Morse boundaries of the edge spaces removed:  
    \[\mb \Gamma - \left (\bigcup_{E \text{ edge group}} (\mb E, \Gamma)\right )\]
     Let $N$ be a Morse gauge and let $v$ be a vertex in a Cayley graph $X$ of $\Gamma$. Let $(\gamma_n)$ be a sequence of $N$-Morse geodesic lines in $X$  converging to a line $\gamma$ which contain $v$ and whose endpoints  $[\gamma^+]$, $[\gamma^-]$ lie in $\mc A$. Suppose that $[\gamma_n^-]$ and $[\gamma_n^+]$ lie in a connected component $C_n\subseteq \mc A$ for all $n$. Then there exists $i_0$ and a vertex group $V$ of $\Gamma$ so that for all $i \ge i_0$, $[\gamma^+]$, $[\gamma^-]$, $[\gamma_i^+]$, $[\gamma_i^-] \in (\mb V, \Gamma) \cap \mc A$. 
\end{lem}
\begin{proof}
   Follow the proof by contradiction of Lemma 6.7 in \cite{karrer2024subgroups}, which relies on Lemmas 3.7 and 3.2 from \cite{fioravanti2022connected}. Instead of these original lemmas, use the corresponding variants that are obtained by changing the setup in \cite{fioravanti2022connected} to the setup here.  This will lead to the conclusion that $[\gamma^-]$ and $[\gamma^+]$ must lie in the relative Morse boundary of an edge group. However, this is a contradiction, as by assumption, $[\gamma^+]$, $[\gamma^-]$ are contained in $\mc A$.
\end{proof}

\section{Morse accessibility and its consequences} 
\label{sec:Morse-accessibility}

A proposed definition of Morse accessibility was introduced in \cite{karrer2024subgroups}. In this section, we provide a formal definition of \emph{Morse accessibility} which is a slight variation of the definition given in that paper (see \cref{rem_varMorse-acc} for further explanation of the difference). In order to formulate our definition, we introduce \emph{Morse hierarchies} and then methodically transfer all the results proven for group splittings along finitely generated undistorted subgroups with empty relative Morse boundaries from \cite{fioravanti2022connected} to Morse hierarchies. This will give \cref{thmMorseaccessiblegroups} (\cref{thmMorseaccessiblegroups_Intro}) and \cref{corMorseaccessiblegroupsdiscon} (\cref{corMorseaccessiblegroups_Intro}) which provide structural consequences of a group being Morse accessible.

Additionally, \cref{lem:Cardrelmorse} shows that the relative Morse boundary of a finitely generated undistorted subgroup of a finitely generated group cannot consist of just one point---i.e., it will contain 0, 2, or infinitely many points. 

\subsection{Definition of Morse hierarchies}\label{sec:Def_Morse_hierarchy}
Recall we will use \cref{conv} throughout this section when discussing splittings; in particular, actions on Bass--Serre trees are assumed to be non-elliptic and minimal.
\begin{defn}[Hierarchies]
A \emph{hierarchy} $\mc H$ of a group $\Gamma$ is a labeled rooted tree where each vertex is labeled by a subgroup of $\Gamma$ and the root $r$ is labeled by $\Gamma$. For each non-leaf vertex $v \in \mathcal V(\mc H)$, the associated group $\Gamma_v$ splits as finite graph of groups $\mathcal{G}(\Gamma_v)$ and the descendants of $v$ in $\mathcal V(\mc H)$ are labeled by the groups associated to the vertices in $\mathcal{G}(\Gamma_v)$. A hierarchy is \emph{finite} if $\mathcal{H}$ a is a finite tree.  In this case, we say the hierarchy has height $h$ where \[h=\max\{d(r,v)\mid v \in \mathcal V(\mathcal H)\},\] with $r$ denoting the root of $\mc H$. 
A subgroup $V$ of $\Gamma$ that is conjugate to $\Gamma_v$ for some vertex $v$ is called a \textit{vertex group} of $\mc H$ similar to our convention used for a graph of groups decompositions in \cref{sec:graphsofgroups}. If, in addition, $v$ is a leaf of $\mc H$, we call any subgroup $L$ conjugate to $\Gamma_v$ a \textit{leaf group}. 
\end{defn}

\begin{remark}
       Let $H$ be a  subgroup of a finitely generated group $\Gamma$ with connected relative boundary $(\mb H, \Gamma)$. Since the Morse boundary of a finitely generated group is Hausdorff, the connected subset $(\mb H, \Gamma)$ is empty, a single point, or consists of infinitely many elements.
\end{remark}

\begin{defn}[Morse hierarchies] \label{def_Morsehierarchy}
    A hierarchy $\mc H$ of a finitely generated group $\Gamma$ is called  \emph{weak Morse} if it is finite and for each non-leaf vertex $v\in \mathcal V(\mathcal H)$, the (finite) graphs of groups $\mathcal{G}(\Gamma_v)$ corresponds to a splitting of $\Gamma_v$ over finitely generated undistorted subgroups of $\Gamma_v$ that have empty relative Morse boundary  in $\Gamma$. A weak Morse hierarchy is called a \emph{Morse hierarchy} if for each leaf $v$ of $\mathcal H$, $(\mb\Gamma_v,\Gamma)$ is a connected subspace of $\mb\Gamma$.
\end{defn}

\begin{defn}
    A finitely generated group $\Gamma$ is \textit{Morse accessible} if there exists a \Morsetree~ of $\Gamma$.
\end{defn}

\begin{remark}\label{rem:hierarchy_defn_properties}We make the following remarks about the definition of Morse hierarchies. 
\begin{itemize}
    \item If every splitting in a hierarchy $\mc H$ is weak Morse in the sense of Definition \ref{def:weak-Morse-splittings}, then we have a weak Morse hierarchy, but it is important to note, that this is a stronger condition than is asked for a weak Morse hierarchy in this sense: in a weak Morse hierarchy, we ask that every edge group has empty relative Morse boundary in the ambient group $\Gamma$, which does not imply it has empty relative Morse boundary in the vertex group being split. For example, see Remark \ref{rem_relMor}. 
    This implies that in a (weak) Morse hierarchy of height at least $2$ every splitting need not be weak Morse. 
    \item A (weak) Morse hierarchy of height $1$, is exactly a (weak) Morse splitting.
    \item It follows from \cref{lemma:transitive} and \cref{lem:undistorted} in the next subsection that all vertex and edge groups of a (weak) Morse hierarchy are also finitely generated and undistorted in $\Gamma$. 
    \item By \cref{lem:gC}, conjugate subgroups have homeomorphic relative Morse boundaries. Thus, a leaf group $L$ of a Morse hierarchy has connected relative Morse boundary if and only if the group $\Gamma_v$ conjugate to $L$ also has that property.  
    \end{itemize}
\end{remark}

We end this subsection by ruling out the possibility of a singleton connected component for a leaf group in a Morse hierarchy before moving on to the main results of the section.

\begin{lem} \label{lem:noMorseelliptic}
Suppose $H$ is an undistorted subgroup of a finitely generated group $\Gamma$. If $H$ has an undistorted infinite cyclic subgroup, then 
 $|(\mb H, \Gamma)| \neq 1$. In particular, the relative Morse boundary of that cyclic subgroup, consists of two points or is empty. 
\end{lem}

\begin{proof}
    Assume $H$ is a subgroup of $\Gamma$ as in the statement of the lemma and that $|(\mb H, \Gamma)| = [\xi]$ consists of a single point. Let $\langle g \rangle$ be an undistorted infinite, cyclic subgroup in $H$. Note that by the definition of the relative Morse boundary and the fact that $|(\mb H, \Gamma)|$ consists of only one point $[\xi]$, we must have $g^i[\xi] = [\xi]$ for all $i \in \N$. Let $\xi$ be an $N$-Morse representative of $[\xi]$ based at $e$. 
    
   We first claim that $\{g^ie\}_{\mathbb{Z}_{\geq 0}}$ is a $N'$-Morse quasi-geodesic and $[\{g^ie\}] = [\xi]$. 
   For each $i \in \N$ consider the ideal geodesic triangle $\xi \cup g^i\xi \cup [e, g^i]$. By Lemma 2.4 of \cite{CharneyCordesMurray}, since both $\xi$ and $g^i\xi$ are both $N$-Morse, we know that $[e, g^i]$ is $N'$-Morse for some $N$ depending only on $N$. 
    Since $\langle g \rangle$ is undistorted we know that  $d(e, g_i) \to \infty$ as $i \to \infty$. Furthermore, by we know that each $[e, g_i]$ is bounded Hausdorff distance (depending on $N'$ and the quasi-geodesic constants of the orbit of $\langle g \rangle$) from $\{e, g, g^2, \ldots, g^i\}$. Thus by Arzela-Ascoli and Lemma A.5 of \cite{corrigendum}, we know that $[e, g_i]$ sub-sequentially converges to a $N'$-Morse geodesic $\alpha$ bounded Hausdorff distance from the orbit of the non-negative powers of $g$. Thus $[\alpha] \in (\mb H, \Gamma)$ and $[\alpha] = [\xi]$. 

    We can repeat this argument with the non-positive powers of $g$ to get another $N'$-Morse geodesic $\beta$ bounded Hausdorff distance from the orbit of the non-positive powers of $g$ and since $\langle g \rangle$ is undistorted, it cannot be that $[\beta]=[\xi]$.
\end{proof}

\begin{cor}\label{lem:Cardrelmorse}
 If $H$ is a finitely generated subgroup of a finitely generated group $\Gamma$ and $H$ is undistorted in $\Gamma$ then $(\mb H, \Gamma)$ consists of $0$, $2$ or infinitely many points.  
\end{cor}

\begin{proof}
The subgroup $H$ is finitely generated and undistorted in  $\Gamma$ i.e. the inclusion $\iota: H \emb \Gamma$ is a quasi-isometric embedding. Recall that  $\mathcal M(H,\Gamma)$ is the subspace of $\mb H$ consisting of Morse rays from $H$ that remain Morse when viewed in $\Gamma$. 
    By \cref{lem_jota},
    the map $\iota$ induces a natural bijection
    \[\iota_*: \mathcal M(H,\Gamma) \to (\mb H, \Gamma) \subseteq \mb \Gamma \]

    By \cite[Cor 5.8]{Qing}, $\mb H$ has either $0$, $2$ or infinitely many elements. 
\begin{itemize}
    \item If $\mb H = \emptyset$ then $\mathcal M(H,\Gamma) = \emptyset$ and $(\mb H, \Gamma)= \iota_*(\mathcal M(H,\Gamma))= \emptyset$.   
    \item If $|\mb H| = 2$, then $H$ contains an infinite cyclic subgroup that is undistorted in $\Gamma$. By \cref{lem:noMorseelliptic},  $(\mb H, \Gamma)= \iota_*(\mathcal M(H,\Gamma))$ consists either of two points or $0$ points. 
    \item Finally, assume that $|\mb H| = \infty$ and that $(\mb H, \Gamma) \neq \emptyset$. Then $H$ contains a Morse geodesic ray $\gamma$ so that $\iota_*([\gamma]) \in (\mb H, \Gamma)$. By \cite[Cor 6.2]{Qing}, the action of $H$ on $\mb H$ is minimal so the orbit $H [\gamma]$ is infinite. 
    Since the map $\iota_*$ is an $H$-equivariant bijection, the image of the  orbit $H [\gamma]$ in $(\mb H,\Gamma)$ is infinite. 
\end{itemize}
Thus the corollary is proven.
 
\end{proof}

\begin{cor} \label{cor:noelliptic}
      If $\Gamma$ is a Morse accessible group, then the relative Morse boundary of any leaf group is either empty or contains infinitely many points. 
\end{cor}
\begin{proof}
    By \cref{lem:undistorted} in the next subsection, every leaf group is finitely generated and undistorted. Thus, \cref{lem:Cardrelmorse} applies.
\end{proof}

\subsection{From graphs of groups to hierarchies: the proofs of \cref{thmMorseaccessiblegroups_Intro} and \cref{corMorseaccessiblegroups_Intro}}\label{sec:proofThmA}

In this subsection, we will prove \cref{thmMorseaccessiblegroups} and \cref{corMorseaccessiblegroupsdiscon} which are \cref{thmMorseaccessiblegroups_Intro} and \cref{corMorseaccessiblegroups_Intro} from the introduction. 
The strategy is to carefully generalize the results on splittings of groups along subgroups with empty relative Morse boundaries from \cite{fioravanti2022connected} to hierarchies using induction. \\

Let $\Gamma$ be a finitely generated group with finite hierarchy $\mc H$. Recall that we call a subgroup of $\Gamma$ a vertex group of $\mc H$ if it is conjugate to a vertex label of $\mc H$. If $\mc H$ has height one, $\Gamma$ splits as graph of groups and the vertex groups of $\mc H$ other than $\Gamma$ are the stabilizers of an associated Bass--Serre tree. In order to generalize the results in \cite{fioravanti2022connected} via induction, we introduce the rooted \emph{blow-up tree} $ \overline{\mathcal H}$ of $  \mathcal H$ with an action $\Gamma \acts \blowup$ for which the vertex groups of $\mc H$  occur as vertex stabilizers of $\overline{\mathcal H}$. 
The vertices of $\overline{\mathcal H}$ will be labeled by the left cosets of $\Gamma_v$ for each vertex $v \in \mathcal H$ and the edges are induced by the edges of the hierarchy $\mathcal H$ and inclusion of cosets.

\begin{defn}[blow-up tree  of $\mathcal H$]\label{def:blowuptree}
Let $\Gamma$ be a finitely generated group with hierarchy $\mathcal H$ and root $r$.  
The blow-up tree $\overline{\mathcal H}$ of $\mathcal H$ is the labeled rooted tree with
\begin{itemize}
    \item $ \mathcal V(\overline{\mathcal H}) =\{w_{g\Gamma_{v}}\mid g \in \Gamma \text{ and } v \in \mathcal V(\mathcal H)\}$
  \item $\mathcal E(\overline{\mathcal H}) =\left\{\{w_{g\Gamma_{u}},w_{h\Gamma_{v}}\}\ \middle| \ \begin{aligned} &\{u,v\}\in \mathcal E(\mc H)  \text{ and } \\
  &g\Gamma_{u} \subseteq h\Gamma_{v} \text{ or }g\Gamma_{v} \subseteq h\Gamma_{u}\end{aligned}\right\}$
  \item \text{ The root of $\overline{\mathcal H}$ is $w_{\Gamma_{r}}$}.
\end{itemize}
\end{defn}

\begin{remark}\label{rem:blowuptree} 
   Let $w_{g\Gamma_v}$ be a vertex of $\overline {\mc H}$ with $g \in \Gamma$ and $\Gamma_v$ the group associated to the vertex $v \in \mc H$. Let $u_1,\dots,u_n$ be the descendants of $v$ in $\mc H$, associated to $\Gamma_{u_1},\dots, \Gamma_{u_n}$ respectively. 
    By definition, the set of descendants $D$ of $w_{g\Gamma_v}$ are $$D=\{w_{gh\Gamma_{u_i}} \mid h \in \Gamma_{v}, i \in \{1,\dots,n\} \}.$$   
Let $T_v$ be a Bass--Serre tree corresponding to $\mc G(\Gamma_v)$; By construction, the vertices of $T_v$ are labeled by the the left cosets of the vertex labels of $\mc G(\Gamma_v)$. Then the vertex labels of $D$ coincide with the vertex labels of the Bass--Serre tree $gT_v$, which are obtained from $T_v$ by left multiplication by $g$ on the vertex labels of $T_v$.
The conjugate $g\Gamma_vg^{-1}$ acts on $gT_v$ by multiplying the vertex labels by $g$ on the left. 
\end{remark}

Recall that the  \emph{$i$th level} of a rooted tree denotes the set of all vertices at distance $i$ from the root. \cref{rem:blowuptree} implies:
\begin{lem}
Suppose $\blowup$ is the blow-up tree of a hierarchy $\mc H$. Then the $i$th level of $\blowup$ consists of all the left cosets $g\Gamma_w$, $g \in \Gamma$ of subgroups $\Gamma_w$ that occur as a label of a vertex in the $i$th level of $\blowup$. In particular, the height of $\mc H$ and $\blowup$ coincide. 
\end{lem}

We observe:

\begin{lem}[The action of $\Gamma$ on the blow-up tree $\blowup$]
    Let $\Gamma$ be a finitely generated group with hierarchy $\mathcal H$ and $\overline {\mc H}$ the blow-up tree of $\mathcal H$. 
    Then $\Gamma$ acts on $\overline {\mc H}$ by left multiplication on the vertex labels. The action fixes the root of $\overline {\mc H}$ and is level-preserving. If $w_{g\Gamma_u}$ is a vertex of $\overline {\mc H}$  then $g\Gamma_u g^{-1}$ is its stabilizer.
\end{lem}

\begin{cor}\label{leafgroups}
     Let $\Gamma$ be a finitely generated group with finite hierarchy $\mathcal H$. 
     Then the set of vertex groups $$\{g\Gamma_{v}g^{-1} \mid g \in \Gamma, v \in \mathcal V(\mathcal H)\}$$ of $\mc H$ coincides with the set of vertex stabilizers $$\{\stab(w) \mid w \in \mathcal V(\overline{\mc H})\}.$$ 
And, the set of leaf groups $$\{g\Gamma_{v}g^{-1} \mid g \in \Gamma, v \text{ is a leaf of $\blowup$ }\}$$ coincides with the set of leaf stabilizers
 $$\{\stab(w) \mid w \text{ is a leaf of $\blowup$ }\}.$$

     If the height of $\mc H$ is $1$, then $\Gamma$ splits as graph of groups that is weak Morse and the set of vertex stabilizers of an associated Bass--Serre tree coincides with both, the leaf groups and the non-root vertex groups of $\mc H$. 
\end{cor}

The above enables us to generalize the results in  \cite{fioravanti2022connected} to weak Morse hierarchies via induction. Here, we use the following crucial property of weak Morse hierarchies: 
\begin{remark} [Induction argument for weak Morse hierarchies]
    If a group $\Gamma$ has a weak Morse hierarchy $\mc H$ of height $n$, then each sub-tree induced by all the vertices of distance $m\le n$ to the root of $\mc H$ yields a \incompleteMorsetree~ of height $m$ for $\Gamma$.  We will denote this sub-tree by $\mc H_m$. Moreover, the sub-tree of the blow-up tree $\blowup$ induced by all the vertices of distance $m\le n$ to the root of $\blowup $ yields the blow-up tree $\blowup_m$ of $\mc H_m$ and the action of $\Gamma$ on $\blowup_n$ restricts to the action of $\Gamma$ on $\blowup_m$  for all $m \le n$. That is, if one restricts the action of $\Gamma$ on $\blowup_n$ to the subtree $\blowup_m$ then this yields the action of $\Gamma$ on the blow-up tree $\blowup_m$ of the weak hierarchy $\mc H_m$ for all $m \le n$. 
\end{remark}

As a first step,   the described approach can be used to extend Lemma 2.15 from \cite{fioravanti2022connected} which shows that in a graph of groups decomposition of $\Gamma$, if all edge groups are finitely generated and undistorted, then so are all vertex groups. Similarly, an induction argument shows how the vertex groups and their relative Morse boundaries sit in the ambient group $\Gamma$ when it has a weak Morse hierarchy. In summary, we obtain:  

\begin{lem}\label{lem:undistorted}
    Let $\mc H$ be hierarchy  of a finitely generated group $\Gamma$. 
   \begin{itemize}
        \item If every edge group of $\mc H$ is finitely generated, then every vertex group of $\mc H$  is finitely generated and undistorted in $\Gamma$. 
        \item Suppose $a$ is a vertex of $\overline {\mc{H }}$ that lies on the unique path of the root of $\overline {\mc{H }}$ to another vertex $b$ of $\overline {\mc{H }}$. Then the vertex stabilizer $B$ of $b$ is properly contained in the vertex stabilizer $A$ of $a$. In particular, if $\mc H$ is weak Morse, then the relative Morse boundaries satisfy $(\mb A, \Gamma)\subseteq(\mb B, \Gamma) $.
    \end{itemize}
\end{lem}

Theorem A in \cite{fioravanti2022connected} (see \cref{thmA}) shows that if a finitely generated group $\Gamma$ has a weak Morse splitting, then the relative Morse boundaries of two distinct vertex groups have empty intersection. The following lemma extends this result to weak Morse hierarchies of arbitrary height. 

If $v$ and $w$ are two vertices in a tree, we denote the unique path between $v$ and $w$ by $[v,w]$ and its vertex set by $\mathcal V([v,w])$. We abuse notation by writing $x\in [v,w]$ to refer to a vertex $x$ along this path. 

\begin{lem}\label{lem:emptyintersection}
Let $\mc H$ be a weak Morse hierarchy of a finitely generated group $\Gamma$. 
Let $\overline{\mc H}$ be its blow-up tree with root $r$. 
If $A$ and $B$ are two vertex groups of $\mc H $ such that  $(\mb A, \Gamma) \cap (\mb B, \Gamma) \neq\emptyset$ then $A$ is the stabilizer of exactly one vertex $a$ of $\blowup$ and $B$ is the stabilizer of exactly one vertex $b$ from  $\blowup$ and $[r, a ] \subseteq [ r, b]$ or $[r, b]\subseteq [r,  a]$. 
\end{lem}
\begin{proof}
   We use induction on the height of a weak Morse hierarchy.\\
   
\textbf{Induction base:} Let $\mathcal H_1$ be a \incompleteMorsetree \ of a finitely generated group $\Gamma$ of height 1. This is simply a graph of groups decomposition of $\Gamma$ and the vertex groups of $\mc H_1$ other than $\Gamma$ are the vertex stabilizers of a corresponding Bass--Serre tree  by \cref{leafgroups}. In this case, Theorem A  of \cite{fioravanti2022connected} (see \cref{thmA}) directly implies the claim.\\

\textbf{Induction step:} Let $\mathcal H_{n+1}$ be a \incompleteMorsetree \ of a finitely generated group $\Gamma$ of height $n+1$. Let $\blowup_{n+1}$ be its blow-up tree and $A$ and $B$ be two vertex groups of $\mc H_{n+1}$ with $(\mb A,\Gamma) \cap (\mb B, \Gamma) \neq \emptyset$.  By \cref{leafgroups},  $A$ and $B$ are each the stabilizers of at least one vertex of $\blowup_{n+1}$. 

Let $a, b$ be vertices of $\blowup_{n+1}$ such that $A = \stab(a)$ and $B = \stab(b)$. The following two claims imply that then $[r, a ] \subseteq [ r, b]$ or $[r, b]\subseteq [r,  a]$. \\

\textbf{Claim 1: } 
    Either $\large([r,a]\subseteq [r,b]$ or  $[r,b]\subseteq [r,a]\large)$ OR  $a$ and $b$ are leaves of $\blowup_{n+1}$ with a common ancestor $z \in \mathcal V(\blowup_{n+1})$.\\

\textbf{Claim 2:}  
If $a$ and $b$ are leaves in $\blowup_{n+1}$ with a common ancestor then $a=b$. \\

We now prove the uniqueness of $a$ and $b$ assuming Claims 1 and 2. \\

\textbf{Uniqueness of $a$ and $b$:}  We first prove uniqueness assuming Claim 1 and Claim 2. By symmetry it suffices to show that the vertex $a$ is unique. Suppose $a$ and $a'$ are two vertices of $\blowup_{n+1}$ with $\stab(a)=\stab(a')=A$. By assumption, $(\mb A, \Gamma) \cap (\mb B, \Gamma) \neq \emptyset$. Hence, $(\mb A, \Gamma) \neq \emptyset$. Applying Claim 1 and Claim 2 for the case where $A=B$ yields that $a$ and $a'$ lie on a common path that connects the root $r$ of $\blowup_{n+1}$ with a leaf of $\blowup_{n+1}$. However, the stabilizers of two consecutive vertices of such a path are properly contained in each other by \cref{lem:undistorted}, and hence, $a=a'$. \qed \\

We will now prove Claim 1 and Claim 2.\\

\textbf{Proof of Claim 1:}
Let $\blowup_n$ be the subtree of $\blowup_{n+1}$ induced by all the vertices of distance at most $n$ from the root $r$ of $\blowup_{n+1}$. Note that $\blowup_{n}$ is the blow-up tree of the \incompleteMorsetree~ $\mc H_n$ of $\Gamma$ of height $n$ for which the induction hypothesis applies. Thus, we may assume that $a$ or $b$, say $a$, is a leaf of $\blowup_{n+1}$. Let $\bar a$ be the ancestor of $a$. Let $\bar b$ be the vertex closest to the root $r$ with distance at most $1$ to $b$, i.e. if $b = r$ then $\bar b= r$ and otherwise $\bar b$ denotes the ancestor of $b$. Note that $\bar a$ and $\bar b$ are contained in $\blowup_n$.
Let $\bar A $ and $\bar B$ be the stabilizers of $\bar a$ and $\bar b$ respectively. 
By \cref{lem:undistorted}, $(\mb B, \Gamma)\subseteq (\mb \bar B, \Gamma) $ and $(\mb A, \Gamma)\subseteq (\mb \bar A, \Gamma) $.  Thus $(\mb A, \Gamma) \cap (\mb B, \Gamma) \neq\emptyset$ implies $(\mb \bar A, \Gamma)\cap(\mb \bar B, \Gamma) \neq \emptyset $. As $\bar a$ is a leaf of $\blowup_n$, the induction hypothesis implies that $[r,\bar b] \subseteq[r,\bar a]$.  If $\bar a= \bar b$ then $a$ and $b$ are both leaves with a common ancestor as asked for in claim 1. Otherwise, $b \in [r, \bar a] \subseteq \blowup_n$ and the induction hypothesis implies that $[r,b] \subseteq  [r, \bar a] \subseteq  [r, a]$. \qed\\

\textbf{Proof of Claim 2: }
Let $a$ and $b$ be two leaves of $\blowup_{n+1}$ with common ancestor $z$. Let $A=\stab(a)$, $B=\stab(b)$ and $Z=\stab(z)$. By definition of $\blowup_{n+1}$, there exists $g \in \Gamma$ and $h_a, h_b \in \Gamma_z$ and two (not necessarily distinct) vertex labels $\Gamma_a$ and $\Gamma_b$ of $\mathcal G(\Gamma_z)$ such that $z = w_{g\Gamma_z}$, $a= w_{gh_a\Gamma_a}$ and $b = w_{gh_b\Gamma_b}$. We need to show that  $gh_a\Gamma_a= gh_b\Gamma_b$. \\ 

  Let $\mathcal G(\Gamma_z)$ be the graph of groups associated to $\Gamma_z$ and $T_z$ be a Bass--Serre tree for $\mc G(\Gamma_z)$  whose vertices are labeled with the left cosets of the vertex labels of  $\mc G(\Gamma_z)$. Let  $gT_z$ be the tree that is obtained from $T_z$ by multiplying all vertex labels with $g$ from the left. Note that $Z = g\Gamma_zg^{-1}$  splits as graph of groups $g\mc G(\Gamma) g^{-1}$ and acts on $gT_z$ by multiplying the vertex labels of $T_z$ from the left.   In particular, the tree $gT_z$ contains a vertex $x$ labeled with $gh_a\Gamma_a$ and a vertex $y$ labeled by $gh_b\Gamma_b$. We will show that $x=y$, as then  $gh_a\Gamma_a= gh_b\Gamma_b$. \\

To simplify notation, we assume that $g = \mathrm{id}$, i.e.  $A$ and $B$ are both contained in $\Gamma_z$ that splits as graph of groups $\mc G(\Gamma_z)$. The case $g\neq \mathrm{id} $ follows analogously. \\

The stabilizers $A=\Gamma_a^{h_a}$ and $B=\Gamma_b^{h_b}$ are conjugates of $\Gamma_a$ and $\Gamma_b$ in $\Gamma_z$ and as such undistorted in $\Gamma_z$. 
 By \cref{lem:undistorted}, $\Gamma_z$ is undistorted in $\Gamma$, i.e. the inclusion $\iota: \Gamma_z \emb \Gamma$ is a quasi-isometric embedding. Recall that  $\mathcal M(\Gamma_z,\Gamma)$ is the subspace of $\mb \Gamma_z$ consisting of Morse rays from $\Gamma_z$ that remain Morse when viewed in $\Gamma$. 
    By \cref{lem_jota},
    the map $\iota$ induces a natural bijection
    \[\iota_*: \mathcal M(\Gamma_z,\Gamma) \to (\mb \Gamma_z, \Gamma) \subseteq \mb \Gamma. \]

As $A$ and $B$ are undistorted in $\Gamma_z$, \[(\mb A,\Gamma) \cap (\mb B, \Gamma)=\iota_*(\mathcal M(\Gamma_z,\Gamma)\cap (\mb A, \Gamma_z))\cap \iota_*(\mathcal M(\Gamma_z,\Gamma)\cap(\mb B, \Gamma_z) ).\]
Injectivity of $\iota_*$ implies that
\begin{align*}
    (\mb A,\Gamma) \cap (\mb B, \Gamma)&= \iota_*\left((\mathcal M(\Gamma_z,\Gamma)\cap (\mb A, \Gamma_z))\cap (\mathcal M(\Gamma_z,\Gamma)\cap(\mb B, \Gamma_z) )\right)\\
    &= \iota_*(\mathcal M(\Gamma_z,\Gamma)\cap (\mb A, \Gamma_z)\cap(\mb B, \Gamma_z) ).
\end{align*}

Accordingly, 
 \[(\mb A,\Gamma) \cap (\mb B, \Gamma) \neq \emptyset \
 \Rightarrow \ \mathcal M(\Gamma_z,\Gamma)\cap (\mb A, \Gamma_z)\cap(\mb B, \Gamma_z)\neq \emptyset.\]

Since every edge group of $\mathcal G_z$ has empty relative Morse boundary in $\mb \Gamma$, 
 
\[ \mathcal M(\Gamma_z,\Gamma) \subseteq   \mathcal \mb \Gamma_z-\left(\bigcup_{E \text{ edge group of }\mathcal G_z} (\mb E, \Gamma_z)\right).\]
\cref{thmAvariant} implies that $A$ and $B$  are both the stabilizer of the same vertex of $T_z$. Hence $x=y$. \end{proof}

Theorem A in \cite{fioravanti2022connected} (see \cref{thmA}) and \cref{cor:gC} shows that the infinite Morse boundary components and their stabilizers of a weak Morse splitting arise from the vertex groups.
The following proposition generalizes this result to weak Morse hierarchies of arbitrary height. Recall that if $C$ is an infinite connected component of $\mb \Gamma$ then $\stab(C)$ denotes the stabilizer of $C$ in $\Gamma \acts \mb \Gamma$. 

\begin{prop} \label{propweakMorseaccessiblegroups} Let $\overline{\mathcal H}$ be a blow-up tree of a \incompleteMorsetree \ $\mathcal H$ of a finitely generated group $\Gamma$ with root $r$. 
     For every infinite connected component $C$ of $\mb \Gamma$ there exists exactly one leaf $v$ of $\overline{\mathcal H}$  such every vertex stabilizer $W$ of a vertex $w$ on $[r,v]$  satisfies
    \begin{align*}
        &C \subseteq  (\mb W,\Gamma) \qquad (\dagger)\\
        &\stab(C) \subseteq W  \qquad (\ddagger)
    \end{align*}

\end{prop}
\begin{proof}
We show the claim by induction on the height of a \incompleteMorsetree. \\

\textbf{Induction base:} Let $ \mc H_1$  be \incompleteMorsetree \ of height 1. This is simply a graph of groups decomposition of $\Gamma$ as in \cref{thmA}  and the leaf groups of  $\mc H_1$  are the vertex stabilizers of a corresponding Bass--Serre tree  by \cref{leafgroups}.
Thus $(\dagger)$ immediately follows from \cref{thmA} and $(\ddagger)$ from \cref{cor:gC}.\\ 

\textbf{Induction step:} Let $\mathcal H_{n+1}$ be a \incompleteMorsetree~ of height $n+1$ of a finitely generated group $\Gamma$ and $\blowup_{n+1}$ be its blow-up tree with root $r$. Let  $C$ be an infinite connected component of $\mb \Gamma$. 
The subtree $\mc H_{n}$ defines a \incompleteMorsetree \ of height $n$ of $\Gamma$ for which the induction hypothesis applies. Hence, $\mathcal \blowup_n$ has exactly one leaf $z$ so that each $w \in V([r,z])$ satisfies $(\dagger)$ and $(\ddagger)$. We will extend $[r,z]$ uniquely to a path connecting $r$ to a leaf $v$ of $\blowup_{n+1}$ so that $[r,v]$ satisfies ($\dagger$) and show afterwards that this path satisfies ($\ddagger$).\\

\textbf{Existence of a leaf $v$ so that $[r,v]$ satisfies $(\dagger)$:} 

For the existence, it suffices to show that there exists one descendant $v$ of $z$ whose stabilizer $V$ satisfies $C \subseteq (\mb V, \Gamma)$. 
By definition of $\blowup_{n+1}$, $z=w_{g_z\Gamma_z}$ with $g_z \in \Gamma$ and $\Gamma_z$ a vertex label of $\mc H$. Let $Z=g_z\Gamma_zg_z^{-1}$ be the stabilizer of $z$. 
To simplify notation, we assume in the following that $g_z = \mathrm{id}$, i.e. we assume that $C \subseteq (\mb  \Gamma_z , \Gamma)$. The case $g_z\neq \mathrm{id}$ follows analogously. \\
As in the proof of \cref{lem:emptyintersection}, we  study the natural bijection
    \[\iota_*: \mathcal M(\Gamma_z,\Gamma) \to (\mb \Gamma_z, \Gamma) \]
    By Lemma 2.3 in \cite{fioravanti2022connected} (see \cref{lem_jota}), the map $\iota_*$ is open. As $\iota_*$ is injective, this implies that there exists a connected set $S \subseteq \mathcal M(\Gamma_v,\Gamma)$   so that $\iota_*(S)=C$. Let $C_z$ be the component of $\mathcal M(\Gamma_v,\Gamma)$ containing $S$. By definition of weak Morse hierarchies, $\Gamma_z$ splits as graph of groups $\mc G(\Gamma_z)$ so that 
\[ \mathcal M(\Gamma_z,\Gamma) \subseteq   \mathcal \mb \Gamma_z-\left(\bigcup_{E \text{ edge group of }\mathcal G_z} (\mb E, \Gamma_z)\right)\]
Thus, by \cref{thmAvariant}, there exists a vertex of $\mc G_z$ labeled by $\Gamma_v$ and a group element $g_v \in \Gamma$  so that the connected component $C_z $ of $\mathcal M(\Gamma_z,\Gamma)$ satisfies: 
 \[C_z\subseteq(\mb g_v\Gamma_v \inv{g_v}, \Gamma_z)\cap  \mathcal M(\Gamma_z,\Gamma).\] 

   Using that $\Gamma_v$ is undistorted in $\Gamma_z$ and that $\Gamma_z$ is undistorted in $\Gamma$ (by \cref{lem:undistorted}), this implies: 
    \[C=\iota_*(S)\subseteq\iota_*(C_z) \subseteq \iota_*((\mb g_v\Gamma_v \inv{g_v}, \Gamma_z)\cap \mathcal M(\Gamma_z,\Gamma)) =(\mb g_v\Gamma_v \inv{g_v}, \Gamma)\]

 As the vertex $w_{g_v\Gamma_v}$ is a descendent of $z$ with stabilizer $V=g_v\Gamma_v \inv{g_v}$, this finishes the proof for the existence of a leaf $v=w_{g_v\Gamma_v} $ of $\blowup_{n+1}$ whose stabilizer $V= g_v \Gamma_v \inv{g_v}$ satisfies  $C \subseteq (\mb V, \Gamma)$. 
Thus, $[v,r]$ satisfies $(\dagger)$. \\

    \textbf{Uniqueness of [r,v]:} 
    It remains to show that the choice of the leaf $v$ and the group $V$ is unique.
   However, if $L$ is the stabilizer of a leaf  of $\blowup_{n+1}$ other than $v$, then  then  $(\mb L, \Gamma) \cap (\mb V, \Gamma)=\emptyset $ by \cref{lem:emptyintersection}. Since $C \subseteq (\mb V, \Gamma) $ this implies that $C \nsubseteq (\mb L, \Gamma)$.\\

    \textbf{The path [r,v] satisfies $(\ddagger)$:}
   In order to prove that $[r,v]$ satisfies $(\ddagger)$, it remains to show 
    that $\stab(C)\subseteq V$.  Let $g \in \Gamma$ so that $gC=C$. Recall that each $w \in [r,z]$ satisfies $(\ddagger)$. In particular, $\stab(C)\subseteq \stab(z)=\Gamma_z$. Hence, $g \in \Gamma_z$. \\
    
    Recall from the proof of $(\dagger)$ that 
 \[C_z\subseteq(\mb V, \Gamma_z)\cap  \mathcal M(\Gamma_z,\Gamma)  .\] 

 Hence,  \[g C_z\subseteq g(\mb V, \Gamma_z)\cap g\mathcal M(\Gamma_z,\Gamma).\] 

As $g \in \Gamma_z$, \cref{lem:gC} implies that $g(\mb V, \Gamma_z)= (\mb gV\inv{g}, \Gamma_z)$. Moreover, as $g \in \Gamma_z$, $g\mathcal M(\Gamma_z,\Gamma)=\mathcal M(\Gamma_z,\Gamma)$ also follows from the same reasoning. 
Since $\iota_*$ is $\Gamma_z$-equivariant, this implies:
\[ gC =g\iota_*(S)= \iota_*(gS)\subseteq \iota_*(gC_z)) \subseteq \iota_*\large(g(\mb V, \Gamma_z)\cap g\mathcal M(\Gamma_z,\Gamma)\large)\]   
\[=\iota_*((\mb gV\inv{g}, \Gamma_z)\cap \mathcal M(\Gamma_z,\Gamma)) = (\mb gV\inv{g}, \Gamma).\]

   As $[r,v]$ satisfies (*),  $C \subseteq (\mb V, \Gamma)$. As $C = gC$, we conclude that $C   \subseteq (\mb V, \Gamma) \cap  (\mb g V\inv{g}, \Gamma)$. By \cref{lem:emptyintersection}, this implies that $gV\inv{g} = V$, i.e. $g \in V$.  
 \end{proof}

\cref{thmMorseaccessiblegroups_Intro} will follow from the following two propositions.

\begin{prop}\label{Lem:sigmainL}Let $\Gamma$ be a Morse accessible group with \Morsetree~ $\mc H$. 
  For every infinite component $C$ in $\mb \Gamma$ there exists a leaf group $L$ so that $C = (\mb L, \Gamma)$ and $\stab(C)=L$. 
\end{prop}
\begin{proof}
  Let $C$ be an infinite connected component of $\mb \Gamma$.
By \cref{propweakMorseaccessiblegroups}, there exists exactly one leaf group $L$ of $\mc H$ so that
      $C \subseteq  (\mb L,\Gamma)$ and
         $\stab(C)\subseteq L$.
      By the definition of Morse hierarchies, $(\mb L,\Gamma)$ is connected. Hence $(\mb L,\Gamma) \subseteq C$.  We conclude that that  $C = (\mb L,\Gamma)$.    
      
      It remains to show that $L\subseteq \stab(C)$: Let $g' \in L $ and $\xi \in C$. We must show that $g'\xi \in C$. Let $\gamma:\mathbb R^+ \to\Gamma$ be a geodesic ray with $\xi=[\gamma]$.       
      As  $C = (\mb L,\Gamma)$, a $D$-neighborhood of  $L$ contains $\gamma$. In other words, there exists $D>0$ satisfying the following: For any $t\ge 0$ there exists $h \in L$ with $d(\gamma(t),h)\le D$. 
      Accordingly, $d(g'\cdot \gamma(t),g'\cdot h)\le D$. As $g' \cdot h \in L$, this implies that $g'\cdot \gamma(t)$ is at bounded distance from $L$ for all $t\ge 0$. This shows that $g' \xi = g'[\gamma]=[g'\gamma] \in (\mb L,\Gamma)=C$.
\end{proof}

\begin{prop}\label{Lem:LinSigma}
    Let $\Gamma$ be a Morse accessible group with \Morsetree~$\mc H$. If $L$ is a leaf group with nonempty relative Morse boundary $(\mb L, \Gamma)$, then $(\mb L, \Gamma)$ is an infinite connected component of $\mb \Gamma$ and $L$ is its stabilizer in $\Gamma \acts \mb \Gamma$. 
\end{prop}

\begin{proof}
    
Let $L$ be a leaf group such that $(\mb L, \Gamma)$ is connected and contains at least two points. Then  there exists exactly one infinite connected component $C$ in $\mb \Gamma$  that contains $(\mb L, \Gamma)$. By \cref{Lem:sigmainL}  there exists a leaf group $L'$ so that $C = (\mb L', \Gamma)$. Note that $(\mb L, \Gamma) \cap (\mb L', \Gamma) \neq \emptyset$. By \cref{lem:emptyintersection}, $L=L'$.
\end{proof}

\begin{thm}[\cref{thmMorseaccessiblegroups_Intro}]\label{thmMorseaccessiblegroups}
    Let $\Gamma$ be a Morse accessible group and $\mathcal H$ a \Morsetree~ of $\Gamma$. Then the following two sets are identical: 
    \begin{itemize}
        \item $\mathfrak{L}$, the set of leaf groups $L$ of $\mc H$ with $(\mb L, \Gamma)\neq \emptyset$
        \item $\mathfrak{S}$, the set of stabilizers of infinite connected components of $\mb \Gamma$. 
        \end{itemize}
        Moreover, every infinite connected component $C$ in $\mb \Gamma$ satisfies \[C = (\mb \stab (C), \Gamma).\]
  \end{thm} 

\begin{proof}
    Follows directly from \cref{Lem:sigmainL} and \cref{Lem:LinSigma}.
\end{proof}

In \cref{lem:noMorseelliptic} we prove that no leaf group of a Morse accessible group can have a singleton as relative Morse boundary, therefore we can conclude 
\cref{corMorseaccessiblegroups_Intro}:

\begin{cor}[\cref{corMorseaccessiblegroups_Intro}]\label{corMorseaccessiblegroupsdiscon} 
         Let $\Gamma$ be a Morse accessible group with a \Morsetree~ $\mathcal H$. The Morse boundary of $\Gamma$ is totally disconnected if and only if for every leaf group $L$ of $\mc H$ we have $(\mb L, \Gamma)=\emptyset$.
    \end{cor}

\subsection{Point-convergence of the set of infinite Morse boundary components} \label{sec:rigidity}

In this section, we prove that the set of infinite Morse boundary components of a Morse accessible group is point-convergent. By Corollary 1.2, Corollary 1.3, and Theorem 1.5 in \cite{karrer2024subgroups}, this directly implies the quasi-isometry rigidity statements \cref{thm:point-convergence_intro} and \cref{corSurfaces_intro} in the introduction. 

\textit{Point-convergence}, a notion introduced in \cite{karrer2024subgroups} which generalizes the notion of a null-family.  

\begin{defn}
    Let $Y$ be a topological space and let $(A_n)_{n\in \N}$ be a sequence of subsets of $Y$. The limit 
    \begin{align*}
        L = \lim_{n\to \infty} A_n
    \end{align*}
    is the set of points $x\in Y$ which can be written as $x =\displaystyle{\lim_{n\to \infty} x_n}$ for points $x_n\in A_n$.
\end{defn}

\begin{defn}\label{definition:point-convergence}
    Let $Y$ be a topological space. A family of subsets $\{A_i\}_{i\in I}$ of $Y$ is \emph{point-convergent} if for all sequences $(A_{i_n})_{n\in \N}$ of distinct elements of $\{A_i\}_{i\in I}$ we have that
    \begin{align*}
        |{\lim_{n\to \infty} A_{i_n}}|\leq 1.
    \end{align*}
\end{defn}

We show the following stronger variant of \cref{thmintro:point-convergende}: 
\begin{thm} \label{prop_point-converg}
    Let $V$ be a vertex group of a Morse hierarchy of a finitely generated group $\Gamma$. Then the set of infinite connected components of the subspace $(\mb V, \Gamma)$ is point-convergent. In particular, if $\Gamma$ is Morse accessible then the set of its infinite Morse boundary components is point-convergent.
\end{thm}

In order to show \cref{prop_point-converg} we need the following compatibility lemma. 
In this lemma one might consider the special case where $\Gamma'$ is a vertex group of a Morse hierarchy. In that case, the lemma shows that if $\Gamma'$ has a Morse hierarchy, then the leaf groups with non-empty relative Morse boundary in $\mb\Gamma'$ are compatible with the leaf groups of the Morse hierarchy of the ambient group $\Gamma$.

\begin{lem}\label{lem:thmVar_subgroup}
Let $\Gamma$ be a finitely generated group and $\Gamma'$ a finitely generated subgroup of $\Gamma$ which is undistorted in $\Gamma$. Suppose that $\Gamma'$ splits non-trivially as a graph of groups so that every edge group has empty relative Morse boundary in $\Gamma$ and is finitely generated and undistorted in $\Gamma'$. 
 For every infinite connected component $C$ of $(\mb \Gamma', \Gamma)$ there exists exactly one vertex group $V$  in the splitting of $\Gamma'$ so that  $C$ is contained in $(\mb V, \Gamma)$;
\end{lem}

\begin{proof} 
Let $C$ be an infinite component of $(\mb \Gamma', \Gamma)$ and $\iota \colon \Gamma' \to \Gamma$ be a the inclusion map. By assumption, $\iota$ is a quasi-isometric embedding. We consider 
\[\mathcal M(\Gamma',\Gamma) \coloneqq \{[\gamma] \in \mb \Gamma' \mid \iota(\gamma) \text{ is Morse in }\Gamma\} \subseteq \mb \Gamma'. \]
Let $\iota_*$ be the map $\mathcal M(\Gamma',\Gamma)\to (\mb \Gamma', \Gamma) \subseteq \mb \Gamma$ induced by $\iota$, i.e $\iota_*[\gamma]=[\iota_*(\gamma)]$.
By the discussion in \cref{sec:prelims}, $\iota_*$ is a bijection. Moreover, $\iota_*$ is an open map by \cref{lem_jota}. Thus, $\iota^{-1}(C)$ is a connected subset of $\mathcal M(\Gamma',\Gamma)$.

 By assumption, $(\mb E,\Gamma)=\emptyset$ for every edge group $E$ of $\Gamma'$. Thus, $C \cap (\mb E,\Gamma) =\emptyset$ for every edge group $E$ of $\Gamma'$ and
  \[\mathcal M(\Gamma',\Gamma)   \subseteq \mb \Gamma' - \left (\bigcup_{E \text{ edge group of } \Gamma'} (\mb E, \Gamma')\right ).\]

By \cref{thmAvariant}, there exists exactly one vertex group $V$ of $\Gamma'$ so that $ \iota^{-1}(C)$ is
contained in $(\mb V, \Gamma') \cap \mathcal M(\Gamma',\Gamma) $. 
Applying $\iota_*$ gives that $C$ is contained in $\iota_*((\mb V, \Gamma')\cap \mathcal M(\Gamma',\Gamma) )$. It remains to show that $\iota_*((\mb V, \Gamma')\cap \mathcal M(\Gamma',\Gamma) )\subseteq (\mb V, \Gamma)$.
By assumption, every edge group of $\Gamma'$ is finitely generated and undistorted in $\Gamma'$. Thus, the vertex group $V$ is undistorted in $\Gamma'$ by Lemma 2.15 in \cite{fioravanti2022connected}. As $\Gamma'$ is undistorted in $\Gamma$ by assumption, the claim follows from \cref{lemma:iota_compatible}.
\end{proof}

The following lemma can be viewed as a variant of  Lemma 6.7 in \cite{karrer2024subgroups}.

\begin{lem}\label{lem:point-convergence}
Let $\Gamma$ be a finitely generated group and $\Gamma'$ a finitely generated subgroup of $\Gamma$ which is undistorted in $\Gamma$. Suppose that $\Gamma'$ splits non-trivially as a graph of groups so that every edge group has empty relative Morse boundary in $\Gamma$ and is finitely generated and undistorted in $\Gamma'$. 

Let $(\xi_i)_{i \in \mathbb N}$ and $(\zeta_i)_{i \in \mathbb N}$ two convergent sequences in $(\mb \Gamma',\Gamma) $  whose limits $\xi=\lim_{i \to \infty}{\xi_i}$ and $\zeta=\lim_{i \to \infty}{\zeta_i}$ are distinct and contained in $(\mb \Gamma',\Gamma) $.  Suppose that $\xi_i$ and $\zeta_i$ are contained in a common connected component of $(\mb \Gamma',\Gamma) $ for every $i$. After passing through a subsequence, there exists $i_0$ and a vertex group $V$ of $\Gamma'$ such that for all $i\geq i_0$, $\xi_i \in (\mb V, \Gamma)$ and  $\zeta_i \in (\mb V, \Gamma)$ and $\zeta, \xi \in (\mb V, \Gamma)$. 
\end{lem}

\begin{proof} By assumption, $(\mb E,\Gamma)=\emptyset$ for every edge group $E$ of $\Gamma'$. Thus,
  \[\mathcal M(\Gamma',\Gamma) \subseteq \mb \Gamma' - \left (\bigcup_{E \text{ edge group of }\Gamma'} (\mb E, \Gamma')\right ).\]
Let $\bar \xi_i \coloneq \iota_*^{-1}(  \xi_i)$, $\bar \zeta_i \coloneq \iota_*^{-1}(\zeta_i)$, $\bar \xi \coloneq \iota_*^{-1}(\xi) $ and $\bar \zeta \coloneq \iota_*^{-1}(\zeta) $. 
As the map $\iota_*$ is an open bijection, $(\bar \xi_i)_i$ converges to $\bar \xi \in M(\Gamma',\Gamma) $ and $(\bar \zeta_i)_i$ converges to $\bar \zeta \in M(\Gamma',\Gamma) $. 
Moreover, $\bar \xi_i$ and $\bar \zeta_i$ lie in a common connected component of $\mathcal M(\Gamma',\Gamma)$  as $\xi_i$ and $\zeta_i$ lie in a common connected component of $(\mb \Gamma',\Gamma) $ and $\iota_*$ is an open bijection.
As $\iota_*$ is a bijection, $\bar \xi$ and $\bar \zeta$ are distinct. Up to possibly passing to a subsequence we might assume that $\bar \xi_i$ and $\bar \zeta_i$ are distinct for all $i$. 

 For each $i$ let $\gamma_i$ be a Morse geodesic line in a Cayley graph $X'$ of $\Gamma'$ with endpoints $\gamma_i^- = \bar \xi_i$ and $\gamma_i^+ = \bar \zeta_i$ (these lines exist as the Morse boundary is a visibility space by \cite{Cordes}). Following the proof of Lemma 3.14 in \cite{karrer2024subgroups} one can show that there exists a Morse gauge $N$ so that all the lines $\gamma_i$ are $N$-Morse in $X'$ and that all the lines $\gamma_i$ pass through a vertex $x \in X'$ after possibly passing to a subsequence.    
 Up to passing to a subsequence, we can assume that $(\gamma_i)_i$ is a sequence converging to an $N$-Morse geodesic line $\gamma$ in $\Gamma'$ with endpoints $\gamma^- = \lim_{i\to \infty} \gamma_i^-= \bar \xi$ and $\gamma^+ = \lim_{i\to \infty} \gamma_i^+= \bar \zeta$ by locally finiteness. 
 By \cref{lem:point-convergence_prem}, there exists $i_0$ and a vertex group $V$ of $\Gamma'$ such that for all $i\geq i_0$, $[\bar \xi_i]=[\gamma_i^+] \in (\mb V, \Gamma') \cap \mathcal M(\Gamma',\Gamma) $, $[\bar \zeta_i]=[\gamma_i^-] \in (\mb V, \Gamma') \cap \mathcal M(\Gamma',\Gamma) $ and $\bar \xi, \bar \zeta \in (\mb V, \Gamma') \cap \mathcal M(\Gamma',\Gamma) $.
Now we apply the map $\iota_*$ and conclude that  for all $i\geq i_0$, $[ \xi_i] \in \iota_*((\mb V, \Gamma)\cap \mathcal M(\Gamma',\Gamma) )$, $[\zeta_i] \in \iota_*((\mb V, \Gamma)\cap \mathcal M(\Gamma',\Gamma) )$ and $ \xi, \zeta \in \iota_*((\mb V, \Gamma)\cap \mathcal M(\Gamma',\Gamma) )$.

  It remains to show that $\iota_*((\mb V, \Gamma')\cap \mathcal M(\Gamma',\Gamma) )\subseteq (\mb V, \Gamma)$.
By assumption, every edge group of $\Gamma'$ is finitely generated and undistorted in $\Gamma'$. Thus, the vertex group $V$ is undistorted in $\Gamma'$ by Lemma 2.15 in \cite{fioravanti2022connected}. As $\Gamma'$ is undistorted in $\Gamma$ by assumption, the claim follows from \cref{lemma:iota_compatible}.
\end{proof}

\begin{proof}[Proof of \cref{prop_point-converg}]
  Let $V$ be a vertex group of a Morse hierarchy $\mc H$ of a finitely generated group $\Gamma$, i.e. $V=g\Gamma_vg^{-1}$ with $g \in \Gamma$ and $\Gamma_v$ a label of a vertex $v$ of $\mc H$. 
   By \cref{lem:gC}, $(\mb V, \Gamma) = (\mb g\Gamma_v g^{-1}, \Gamma)=g(\mb \Gamma_v,\Gamma)$. As $\Gamma$ acts by homeomorphisms on its Morse boundary, this implies that $(\mb V, \Gamma)$ and $(\mb \Gamma_v, \Gamma)$ are homeomorphic. Thus, it suffices to show that for each group $\Gamma_v$ labeling a vertex $v$ of $\mc H$, the set of infinite components of $(\mb \Gamma_v, \Gamma)$ is point-convergent.

    Let $\mc H$ be a Morse hierarchy and $v$ a vertex of $\mc H$.  We argue by induction on the distance of $v$ to the leaf of $\mc H$ closest to $v$. If $v$ is itself a leaf, then this distance is $0$ and $(\mb \Gamma_v, \Gamma)$ is empty or consists of at most one infinite connected component. In particular, the set of infinite connected components of  $(\mb \Gamma_v, \Gamma)$ is point-convergent.
    
Let $v$ be a vertex of $\mc H$ of distance $n+1$ to the leaf of $\mc H$ closest to $v$. Let $\Gamma_v$ be the label of $v$. To proceed a contradiction suppose that $(\mb \Gamma_v,\Gamma)$ is not point-convergent.   Then  there exists a sequence $(C_i)_{i \in \mathbb N}$ of pairwise distinct connected components $C_i$ of $(\mb \Gamma_v,\Gamma) $ and two convergent sequences $(\xi_i)_{i \in \mathbb N}$ and $(\zeta_i)_{i \in \mathbb N}$ in $(\mb \Gamma_v,\Gamma) $ with $\xi_i,\zeta_i \in C_i$ whose limits $\xi=\lim_{i \to \infty}{\xi_i}$ and $\zeta=\lim_{i \to \infty}{\zeta_i}$ are contained in $(\mb \Gamma_v,\Gamma) $ and distinct. 
 By \cref{lem:point-convergence}, there exists $i_0$ and a vertex group $V$ of $\Gamma_v$ such that for all $i\geq i_0$, $\xi_i \in (\mb V, \Gamma)$, $\zeta_i \in (\mb V, \Gamma)$ and $\xi, \zeta \in (\mb V, \Gamma)$. Then \cref{lem:thmVar_subgroup} implies that the connected components $C_i$ are entirely  contained in  $(\mb V, \Gamma)$ for all $i\ge n_0$. In particular, every such $C_i$ is a component of  $(\mb V, \Gamma)$ and $(C_i)_i$ is a sequence of connected components of  $(\mb V, \Gamma)$ that have two distinct points $\xi$ and $\zeta$ in their limit set. Accordingly, the set of infinite connected components of $(\mb V, \Gamma)$ is not point-convergent. By the definition of Morse hierarchies, the vertex group $V$ is conjugate to a group $\Gamma_u$ that occurs as a label of a descendant $u$ of $v$. Repeating the argument of the beginning of this proof yields that  $(\mb V, \Gamma)$ is homeomorphism to $(\mb \Gamma_u, \Gamma)$. Thus, the set of infinite connected components of  $(\mb \Gamma_u, \Gamma)$ is not point-convergent. However, this is a contradiction to the induction hypothesis which applies to the vertex $u$ because the distance of $u$ to the leaf closest to it is $n<n+1$.
\end{proof}

\begin{remark}
   The above shows that the set of infinite connected components of $(\mb V, \Gamma)$ is point-convergent for every vertex group $V$ of a Morse hierarchy. One might think that the injectivity and openness of $\iota_*$ is sufficient for concluding the point-convergence of $\mc M(V, \Gamma)$ in $\mb V$ via the pull-back of $\iota_*$. However, this might not be the case as $\mc M(V, \Gamma)$ might contain connected sets whose images under $\iota_*$ are not connected in $(\mb V, \Gamma)$. 
\end{remark}

\section{(Weak) Morse hierarchies of Coxeter groups}\label{sec:MorseAccessAFCoxeter}  

In this section we prove all Coxeter groups have a weak Morse hierarchy, which we prove in the spherical nexus free (\cref{defn:spherical-nexus}) case to be a \emph{genuine} Morse hierarchy (\cref{thmRACG_intro}). This shows that every Coxeter group without a spherical nexus is Morse accessible. We will begin this section with preliminaries on Coxeter groups, followed by results that ``relativize" the necessary results from \cite{CorLev25}, and then prove the theorem in \cref{subsec:RACGMorseaccessible}. For a discussion of consequences of this theorem and examples see Sections \ref{subsec:introRACGRAAG} and \ref{sec:examples}.

\subsection{Preliminaries on Coxeter groups} \label{subsec:coxeter_definitions}

A \emph{Coxeter graph} is a simplicial graph $\Delta$ whose edges are labeled by integers in $\{2,3, 4, \dots\}$. Given two vertices $s, t \in \Delta$, we let $m_{st}$ denote the label between them. We let $V(\Delta)$ denote the vertex set of $\Delta$.

The \emph{Coxeter group $W_\Delta$} associated to the Coxeter graph $\Delta$ is the group with presentation
\[\left \langle V(\Delta) ~|~ s^2 = 1 ~\text{ for all } s \in V(\Delta),~(st)^{m_{st}} = 1 \text{ for each } (s,t) \in E(\Delta) \right \rangle\] 

A Coxeter group is called a \emph{right-angled Coxeter group} (RACG) if all edges of its Coxeter graph are labeled with $2$. 
In this case, the usual labeling becomes unnecessary for extracting the group's presentation. Consequently, if $\Delta$
denotes the Coxeter graph of a right‑angled Coxeter group, we may omit all edge labels and refer to the resulting simple graph as the \emph{defining graph} of the RACG.
\begin{remark}
    Note that there is another common convention for Coxeter graphs where the are unlabeled edges and edges labeled by $\{4,5,6, \ldots \} \cup\{\infty\}$, if two vertices $s$ and $t$ are connected by an unlabeled edge, $m_{st} = 3$, and if two vertices $s$ and $t$ are not joined by an edge then $m_{st} = 2$ (see \cite[D\'efinition 4]{bourbaki} for instance). In this paper we do \emph{not} use this convention. 
\end{remark}

Given a subgraph $\Lambda \subseteq \Delta$, we say that $\Lambda$ is \emph{induced} if every edge in $\Delta$ between vertices of $\Lambda$ is also an edge of $\Lambda$. 
The subgraph \emph{induced by $\Lambda$} is the the induced subgraph of $\Delta$ with vertex set $V(\Lambda)$.
Since a subgroup generated by a subset of the generators is isomorphic to the Coxeter group on the induced subgraph (see \cite{Davis} for instance), we let $W_\Lambda$ denote the subgroup of $W_\Delta$ generated by the vertices of $\Lambda$. Subgroups arising this way are called \emph{special subgroups} and conjugates of special subgroups are called \emph{parabolic subgroups}.

Let $W_\Delta$ be a Coxeter group. 
We say that $W_\Delta$ is \emph{irreducible} if it does not split as a product of groups, and $W_\Delta$ is \emph{affine} if it is irreducible and a Euclidean reflection group. 
The classification of affine Coxeter groups is well-known 
(see \cite{Davis}).
The \emph{rank} of a Coxeter group $W_\Delta$ is the number of vertices in $\Delta$.
We say that $W_{\Delta}$ is \emph{affine-free} if for any subgraph $\Lambda \subseteq \Delta$, $W_\Lambda$ is not an affine Coxeter group of rank $3$ or larger. 
Note, this is equivalent to saying that $W_\Delta$ has no irreducible affine reflection subgroups of rank $3$ or larger \cite{CapraceWalls}.

Given a Coxeter graph $\Delta$, an induced subgraph $\Lambda$ is called a \emph{spherical subgraph} if $W_\Lambda$ is a finite group. In a RACG, every spherical subgraph is an induced clique.  If $\Delta -\Lambda$ is disconnected, we say that $\Lambda$ is a \emph{separating spherical subgraph}. 
In the proof of \cref{thmRACG_intro} and the two corollaries,  we use the following convention: 
\begin{con} \label{conv:clique}
   If $\Lambda$ does not contain any vertex, we say that $\Lambda$ is \emph{empty}. If $\Lambda$ is empty, then $W_\Lambda$ is the trivial group. Moreover, the empty set is an induced subgraph of any graph which is induced by the empty set, and the empty set is a clique on $0$ vertices. 
\end{con}

\begin{remark}\label{rem:one-endedCoxeter}  
    Given a Coxeter graph $\Delta$, recall that $W_\Delta$ is one-ended if and only if $\Delta$ is not spherical and \emph{does not have} a separating spherical subgraph \cite[Theorem 8.7.2]{Davis}.  In view of \cref{conv:clique}, if $\Delta$ is disconnected, then $\emptyset$ is a separating spherical subgraph and thus $W_\Delta$ is not one-ended. We will use this characterization of one-endedness throughout this section.
\end{remark}

The main ingredient of the proof are the results of Cordes--Levcovitz discussed in \cref{sec:results}. Using these results, we will prove Morse-accessibility of affine-free Coxeter groups in \cref{subsec:RACGMorseaccessible}. Afterwards, we will discuss consequences and examples in \cref{sec:examples}.

\subsection{Relativizing results of Cordes--Levcovitz}\label{sec:results}

In this subsection, we review and ``relativize'' the results on Morse boundaries of affine free Coxeter groups of Cordes--Levcovitz which we later use to give a Morse-hierarchy for certain affine-free Coxeter groups, which includes all RACGs. 

\begin{defn} \label{def:wide_subgroups} 
Let $W_\Delta$ be a Coxeter group.
We call an induced subgraph $\Lambda \subseteq \Delta$ \emph{wide} if there is a decomposition $V(\Lambda) = P \sqcup Q$ (with $Q$ possibly empty) such that there is an edge labeled by $2$ between  any vertex in $P$ and any vertex in $Q$, and additionally either
\begin{enumerate}
    \item  both $W_{P}$ and $W_{Q}$ are infinite, or
    \item $W_P$ is irreducible, affine of rank at least $3$. 
\end{enumerate}
 We call $(P,Q)$ a \emph{wide decomposition} of $\Lambda$. Note that, when $Q$ is not empty, $W_\Lambda \cong W_P \times W_Q$. 
 We also say that $W_\Lambda$ is wide if $\Lambda$ is wide. 
\end{defn}

The following lemma describes which special subgroups have empty relative Morse boundaries in Coxeter groups. 
\begin{lem}\label{emptygraphs}
    Let $\Delta$ be a Coxeter graph, and let $\Lambda$ be an induced subgraph of a graph $\Delta$. The relative Morse boundary of $W_\Lambda$ with respect to $W_\Delta$,
     $(\mb W_\Lambda, W_\Delta)=\emptyset$ if and only if $\Lambda$ is contained in a wide or spherical subgraph. 
\end{lem}
\begin{proof}
As wide subgroups have no Morse rays \cite[Lemma 3.2]{CorLev25}, if $\Lambda$ is contained in a wide subgraph $\overline{\Lambda}$, $W_{\overline{\Lambda}}$, then $(\mb W_\Lambda, W_\Delta)$ will not have any Morse rays. If $\Lambda$ is a spherical subgraph, $W_\Lambda$ is finite and thus does not contain any Morse rays. So $(\mb W_\Lambda, W_\Delta)= \emptyset$.
For the other direction, suppose $\Lambda$ is not contained in a wide or spherical subgraph. Let $v_1,\dots,v_n$ be the vertices of $\Lambda$. Then the group element $v_1v_2,\dots v_n$ has an axis that is not bounded by a flat, so it is a rank-one isometry of $W_\Delta$, implying that  $(\mb W_\Lambda, W_\Delta)\neq\emptyset$ \cite[Proposition 4.5]{Caprace-Fujiwara}. 
\end{proof}

\begin{defn}
    We say that the Coxeter graph $\Delta$ and the Coxeter group $W_\Delta$ are \emph{wide-avoidant} if given any wide subgraph $\Lambda$ and any pair of distinct vertices $s, t \in \Delta$, then there exists a path $\alpha \subseteq \Delta$ from $s$ to $t$ such that $\alpha \cap \Lambda \subseteq \{s, t\}$, i.e., $\alpha \cap E(\Lambda) = \emptyset$.
\end{defn}

\begin{defn}[Special Join]
    Given a Coxeter graph $\Delta$, a \emph{special join} is a triplet $(P, Q, K)$ where $P$, $Q$ and $K$ are disjoint vertex sets of $\Delta$ such that: 
    \begin{enumerate}
        \item $(P,Q)$ is a wide decomposition of the induced subgraph spanned by $P \cup Q$.
        \item $K$ is either empty or $W_K$ is finite
        \item For all $k \in K$ and $s \in P$, there is an edge (with any label) in $\Delta$ from $k$ to $s$.
    \end{enumerate}
\end{defn}

\begin{remark} \label{rmk:special_joins_in_racgs}
    If $W_\Delta$ is a right-angled Coxeter group, then special joins correspond to wide subgroups and vice versa.
    Indeed, let $(P, Q, K)$ be a special join in $\Delta$. Then $W_{P \cup Q \cup K} = W_P \times W_{Q \cup K}$
    is a wide subgroup. Conversely, suppose that $W_{\Lambda}$ is a wide subgroup of a right-angled Coxeter group.
    Then $\Lambda$ has a wide decomposition $(P, Q)$, and $(P, Q, \emptyset)$ is a special join. 
\end{remark}

\begin{defn}[wide-spherical-avoidant] \label{def:wide-spherical-avoidant}
	We say that the Coxeter graph $\Delta$ and the Coxeter group $W_\Delta$ are \emph{wide-spherical-avoidant} if given any special join $(P, Q, K)$ and any pair of vertices $s, t \in \Delta \setminus K$, then there exists a path $\alpha$ in $\Delta$ from $s$ to $t$ such that $\alpha \cap (K \cup P \cup Q) \subseteq \{s, t\}$, i.e., $\alpha \cap E(K \cup P \cup Q) = \emptyset$.
\end{defn}
\begin{thm}[Theorems B and C \cite{CorLev25}]\label{thmCorLev}
    If $W_\Delta$ is an affine-free, one-ended, wide-spherical-avoidant Coxeter group, then $W_\Delta$ has connected (and locally-connected) Morse boundary. 
    
    In the case when $W_\Delta$ is a one-ended right-angled Coxeter group with defining graph $\Delta$, we get the following strengthening. 
    The following are equivalent: 
    \begin{enumerate}
        \item The Morse boundary $\partial_* W_\Delta$ is locally connected and connected.
        \item $\Delta$ is wide-avoidant.
    \end{enumerate}
\end{thm}

Recently Cashen, Chu and Tao proved a full characterization of Morse geodesics \cite{CashenChuTao}; importantly in Corollary 1.5 they give a criterion in the affine case. This corollary allows us to extend Theorem \ref{thmCorLev} to the affine free case:
\begin{cor}\label{thmCorLev-CashenChuTao}
    If $W_\Delta$ is a one-ended, wide-spherical-avoidant Coxeter group, then $W_\Delta$ has connected (and locally-connected) Morse boundary. 
\end{cor} 
\begin{proof}
    In Corollary 1.5 of \cite{CashenChuTao}, they prove that a geodesic is Morse if any only if there is an upper bound on the lengths of subsegments $\gamma'$ of $\gamma$ for which there exists $\Lambda \subset \Delta$ such that $W_\Lambda$ is a wide special subgroup and all edges of $\gamma'$ have label in $\Lambda$. In the construction in \cite{CorLev25}, this is exactly the characterization they use to ensure the ``filter construction" embeds stably in the Cayley graph (\cite[Section 6.2]{CorLev25}). Thus, we can drop the affine free condition from the Cordes--Levcovitz result. 
\end{proof}

This theorem can be generalized to a statement about \textit{relative Morse boundaries} of special subgroups of $W_\Delta$. 

\begin{defn}\label{defn:rel-wide-avoidant}
    Let $\Delta$ be a graph and $\Lambda$ an induced subgraph of $\Delta$. The subgraph $\Lambda$ is \textit{relatively wide-avoidant in $\Delta$} if given any induced wide subgraph $J \subseteq \Delta$, and any vertices $s$, $t\in \Lambda$, there exists a path $\alpha$ in $\Lambda$ from $s$ to $t$ so that $\alpha \cap J \subseteq\{s,t\}$. 

           We also say the subgraph $\Lambda$ is \textit{relatively wide-spherical-avoidant in $\Delta$} if given any special join $(P, Q, K)$ in $\Delta$ and any pair of vertices $s, t \in \Lambda \setminus K$, there exists a path $\alpha$ in $\Lambda$ from $s$ to $t$ such that $\alpha \cap (K \cup P \cup Q) \subseteq \{s, t\}$. 
\end{defn}

\begin{remark} \label{wide-avoidant}
    If $\Lambda$, an induced subgraph of $\Delta$, is relatively wide-spherical-avoidant in $\Delta$, this means that $\Lambda$ as a stand-alone graph is wide-spherical-avoidant. If, in addition, $W_{\Delta}$ is one-ended, then then $\mb W_\Lambda$ is connected and locally connected by \cref{thmCorLev-CashenChuTao}.  
    But the definition of relatively wide-spherical-avoidant requires that $\Lambda$ satisfy an even stronger property: that there must exist a path $\alpha$, as in the definition, that avoids \emph{every} induced wide subgraph in $\Delta$, not just the wide subgraphs in $\Lambda$. This property is the key observation in the proof of the following corollary of \cref{thmCorLev-CashenChuTao}: 
\end{remark}

\begin{cor}[Relative version of \cref{thmCorLev-CashenChuTao}]
\label{CorCorLev}
   Let $W_\Delta$ be a Coxeter group with Coxeter graph $\Delta$.
    Let $\Lambda$ be an induced subgraph of $\Delta$. 
    If $W_\Lambda$ is one-ended and $\Lambda$ is relatively wide-spherical-avoidant in $\Delta$, then the relative Morse boundary $(\mb W_\Lambda, W_\Delta)$ is connected and locally connected. 
\end{cor}

\begin{proof}
    In this proof we will be using the notation and constructions (e.g., fans/filters) from \cite{CorLev25}. Let $\alpha, \beta \in (\mb W_\Lambda, W_\Delta)$. Since $W_\Lambda$ is isometrically embedded in $W_\Delta$, by \cref{lem:RelMorse}, we can assume that $\alpha$ and $\beta$ are geodesics in $W_\Delta$ written in the letters of $\Lambda$. The construction of the fans and thus the filters for wide-spherical-avoidant graphs can proceed as described for induced subgraphs $\Lambda$ that do not contain a separating spherical subgraph and is relatively wide-spherical-avoidant. Therefore the filter $\mathcal{F}$ will be made of Morse geodesics in $W_\Delta$ written in the letters of $\Lambda$. Thus the image of $\partial_* \mathcal{F}$ will be a locally connected, connected subset of $(\mb W_\Lambda, W_\Delta)$. Hence, following the proofs of Theorem 7.4 and Theorem 7.11 in  \cite{CorLev25} yields the claim.
\end{proof}

The following lemma characterizes when an induced subgraph of a Coxeter group is one-ended and relatively wide-avoidant or contained in a spherical or wide subgroup. If such a relatively wide-avoidant subgraph of $\Delta$ is also relatively wide-spherical avoidant, then both of these conditions mean that the relative Morse boundary is connected, see \cref{CorCorLev} for the first case (non-empty relative Morse boundary) and \cref{emptygraphs} for the second (empty relative Morse boundary).
\begin{lem}
   \label{rel-wide-avoidant} Let $\Lambda$ be an induced subgraph of a Coxeter graph $\Delta$. 
 The following are equivalent:
    \begin{enumerate}
        \item $\Lambda$ does not contain any (possibly empty) {proper} induced subgraph $S$ so that
    \begin{itemize}
        \item $S$ is contained in a wide or (possibly empty) spherical subgraph of $\Delta$ and 
        \item $\Lambda - S$ is disconnected
    \end{itemize}
    \item There do not exist two {proper} induced subgraphs $\Delta_1$ and $\Delta_2$ of $\Lambda$ so that 
    \begin{itemize}
        \item $\Lambda = \Delta_1 \cup \Delta_2$ and
        \item $\Delta_1\cap \Delta_2$ is contained in a wide or (possibly empty) spherical subgraph of $\Delta$. 
    \end{itemize}
    \item $\Lambda$ satisfies exactly one of the following: 
    \begin{itemize}
        \item $\Lambda$ is a clique that is contained in a wide subgraph of $\Delta$;
        \item $\Lambda$ is a connected graph without separating spherical subgraphs and is relatively wide-avoidant in $\Delta$.
    \end{itemize}
    \end{enumerate}
\end{lem}

   \begin{proof}
    $(1) \Rightarrow (3)$: Suppose $\Lambda$ satisfies $(1)$. Since it has no separating spherical subgraphs, then $\Lambda$ is connected (in view of \cref{conv:clique}). Let $J$ be a  wide subgraph of $\Delta$, and $s,t \in V(\Lambda)$ be two distinct vertices in $\Lambda$.
    If $\Lambda \cap J \neq \Lambda$, then by assumption $\Lambda - J$ is connected, and therefore there exists a path $\alpha$ connecting $s$ and $t$ in $\Lambda - J$ so that $\alpha \cap J \subseteq \{s,t\}$. If $\Lambda \cap J = \Lambda$, then $\Lambda$ is a subset of a wide subgraph and is thus not relatively wide-avoidant. Moreover, since $\Lambda$ is a subset of a wide subgraph $J$ of $\Delta$, any proper induced subgraph is also a subgraph of $J$, so by assumption cannot disconnect $\Lambda$. Thus $\Lambda$ can't have any separating subgraph, i.e., it is a clique.
    
     $(3) \Rightarrow (1)$: Suppose $\Lambda$ satisfies (3). To arrive at a contradiction, suppose  $\Lambda$ contains a proper induced subgraph $S$ that is forbidden by $(1)$. Then $\Lambda - S$ is disconnected, and hence $\Lambda$ cannot be a clique. By $(3)$ this means that $\Lambda$ is connected and does not contain separating spherical subgraphs. Thus, $S$ can't be a spherical subgraph or empty. Accordingly, the graph $S$ is contained in a wide subgraph $J$ of $\Delta$ so that $\Lambda- S$ is disconnected. Pick two vertices $v$ and $w$ of $\Lambda- S$ that lie in different connected components of $\Lambda- S$. Then every path connecting $v$ and $w$ in $\Lambda$ contains at least one vertex in $J$. However, this means that $\Lambda$ is not relatively wide-avoidant in $\Delta$---a contradiction.   
     
     $(2) \Rightarrow (1)$: We show the contrapositive, i.e., $\neg (1) \Rightarrow \neg (2)$: 
     Let $S$ be a (possibly empty) proper induced subgraph of $\Lambda$ that is contained in a (possible empty) spherical subgraph or a wide subgraph of $\Delta$ so that $\Lambda - S$ is disconnected. Let $\Delta_1$ be the union of $S$ with one connected component of $\Lambda - S$. Let $\Delta_2$ be the union of $S$ with $\Lambda - \Delta_1$. As $\Lambda - S$ is disconnected, $\Delta_1  \neq \Lambda$ and $\Delta_1 \neq \Lambda$.  Moreover, $\Lambda = \Delta_1 \cup \Delta_2$ and $\Delta_1 \cap \Delta_2= S$.  Since $S$ is contained in  a (possible empty) spherical subgraph or wide subgraph of $\Delta$, this shows $\neg (2)$. 

     $(1) \Rightarrow (2)$:  Again, we show the contraposition, i.e., $\neg (2) \Rightarrow \neg (1)$: Assume we have two induced proper subgraphs of $\Lambda$, $\Delta_1$ and $\Delta_2$ whose union is all of $\Lambda$ and that $\Delta_1 \cap \Delta_2$ is contained in a (possibly empty) spherical subgraph or wide subgraph of $\Delta$.  Since $\Delta_1$ and $\Delta_2$ are induced subgraphs whose union is all of $\Lambda$, we know that either $\Delta_1 \cap \Delta_2 = \emptyset $ or $\Delta_1 \cap \Delta_2$ disconnects $\Lambda$. Thus setting $S= \Delta_1 \cap \Delta_2$, we get a proper induced subgraph contained in a (possibly empty) spherical subgraph or a wide subgraph of $\Delta$, this shows $\neg (1)$.
    \end{proof}

In fact, the lemma above gives us a way to check if a induced subgraph corresponding to a one-ended subgroup is relatively wide-avoidant. Indeed, if $\Lambda$ is such an induced subgraph of a Coxeter graph $\Delta$ and $\Lambda$ is not contained in a spherical or wide subgraph of $\Delta$, then any of the three equivalent conditions will verify if a subgraph is relatively wide avoidant. 

The most general result characterizing connectivity of the (relative) Morse boundary, \cref{thmCorLev-CashenChuTao}, requires the stronger wide-spherical-avoidant condition on graphs. There are Coxeter graphs, though, that are wide-avoidant, not wide-spherical-avoidant, but surprisingly(!) have connected Morse boundary. Since we do not know how to detect such graphs, and we are interested in the case when these groups have connected (relative) Morse boundary, we give a condition on a Coxeter graph (spherical-nexus free) that will allow us to not worry about the case when a graph is wide-avoidant but not wide-spherical-avoidant.

\begin{defn}[Spherical nexus] \label{defn:spherical-nexus}
    We call an induced subgraph $K$ of a Coxeter graph $\Delta$ and the corresponding special subgroup $W_K$ a \emph{spherical nexus} if $\Delta$ has a special join $(P,Q,K)$ and a relatively wide-avoidant subgraph $\Lambda$ in $\Delta$ with two vertices $s,t \in \Lambda - K$ such that any path $\alpha \subseteq \Lambda$ from $s$ to $t$ passes through $K$.  
\end{defn}

The existence of a spherical nexus means that a graph cannot be wide-spherical-avoidant. The following lemma follows from the definitions:

\begin{lem} \label{lem_nexus}
    Suppose a Coxeter graph $\Delta$ does not contain any spherical nexus. Then an induced subgraph $\Lambda$ of $\Delta$ is relatively wide-avoidant if and only if it is relatively wide-spherical-avoidant. 
\end{lem}

As pointed out in \cref{rmk:special_joins_in_racgs}, if $\Delta$ is the Coxeter graph of a RACG, then it does not contain a spherical nexus.

\subsection{Proof of weak Morse hierarchies for Coxeter groups} \label{subsec:RACGMorseaccessible}
In this subsection, we give the proof of \cref{thmRACG_intro}, which we restate here: 
\begin{thm}[\cref{thmRACG_intro}]   
\label{thmRACG_} Every Coxeter group has a weak Morse hierarchy in which each leaf group $L$ is conjugate to a special subgroup $W_\Lambda$ that satisfies one of the following: 

   \begin{enumerate}
       \item  $\Lambda$ is a relatively wide-avoidant subgraph of $\Delta$  and $W_\Lambda$ is one-ended. In particular, $(\mb W_\Lambda, W_\Delta) \neq \emptyset$; 
       \label{property1}
       \item  $\Lambda$ is a clique that is either spherical or contained in a wide subgraph of $\Delta$. In particular, $(\mb W_\Lambda, W_\Delta) = \emptyset$. \label{property2}
   \end{enumerate}

Moreover, if $W_\Delta$ does not contain any spherical nexus then for the groups in (\ref{property1}), we additionally have that $(\mb W_\Lambda, W_\Delta)$ is connected and locally connected. 
This gives a Morse hierarchy for Coxeter groups that do not contain any spherical nexus.
\end{thm}
\begin{remark}
    Note by \cref{wide-avoidant}, for a leaf group $W_\Lambda$ with no spherical nexus, the Morse boundary of $W_\Lambda$ as a stand-alone group, $\mb W_\Lambda$, is connected and locally connected. 
\end{remark}

The diagram in \cref{fig:proof_morse_accessibility} illustrates an algorithm on a Coxeter graph $\Delta$  that will yield a weak Morse hierarchy of $W_\Delta$.

\begin{proof}[Proof of \cref{thmRACG_intro} (\cref{thmRACG_})]
\begin{figure}[h]
\label{fig:proof}
    \centering
    \includegraphics[width=0.9\textwidth]{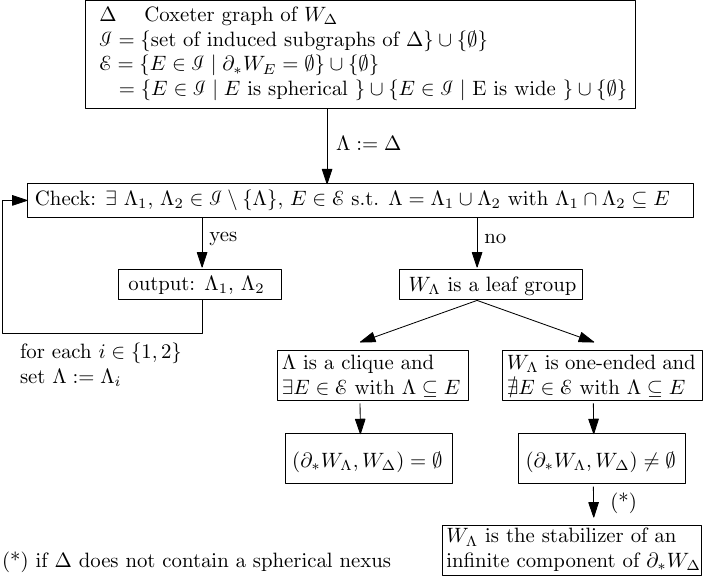}
    \caption{Construction of a weak Morse hierarchy for a Coxeter group $W_\Delta$ of height at most $2^{|V|}$. If $\Delta$ has no spherical nexus, this give a Morse hierarchy. }
    \label{fig:proof_morse_accessibility}
\end{figure}
Let $\Delta$ be a Coxeter graph and let $W_\Delta$ be the associated Coxeter group. Let
\begin{align*}
\mathcal{I} &\coloneqq \{E \subseteq \Delta \mid  E \text{ is an induced subgraph of } \Delta \}\cup \{\emptyset\}\\
    \mathcal E &\coloneq \{E \in \mathcal{I}\mid  \mb W_E = \emptyset\}\cup \{\emptyset\}\\&=~~ \{E \in \mathcal{I}\mid E \text{ is a spherical subgraph or a wide subgraph}\}\cup\{\emptyset\}.
\end{align*}

We build the weak \Morsetree~ $\mathcal H$ we are looking for following the algorithm illustrated in \cref{fig:proof_morse_accessibility}. We will start with a hierarchy $\mathcal H_0$ and will construct hierarchies $\mathcal H_i$ from $\mathcal H_{i-1}$ by possibly adding finitely many edges to the leaves of $\mathcal{H_{i-1}}$. The process will end after finitely many steps, yielding the required hierarchy $\mathcal H$. \\

\textbf{Base case:}
Let $\mathcal H_0$ be the labeled rooted tree that consists of a single vertex (the root itself) labeled it with $W_\Delta$. This is a hierarchy of height $0$. \\

 \textbf{The $i+1^{\text{st}}$ step:} Suppose we constructed the hierarchy $ \mc H_i$ of height $i$. 
 For any leaf $v$ of $\mathcal H_i$ we proceed as follows: By construction, $v$ is labeled by a special subgroup $W_{\Lambda}$ of $W_\Delta$.   We check:

 \begin{align}\label{check}\exists~\Lambda_1,~\Lambda_2 \in \mc I \setminus \{\Lambda\} \text{ and } E \in \mc E \text{ with }
    \Lambda= \Lambda_1 \cup \Lambda_2 \text{ and } 
        \Lambda_1 \cap \Lambda_2 \subseteq E.
\end{align}

 \textbf{Yes: } If there exist such $\Lambda_1$
 and $\Lambda_2$ then $W_\Lambda = W_{\Delta_1}*_{W_{\Lambda_1\cap\Lambda_2}} W_{\Lambda_2}$. In this case, we add two new edges $\{v, v_1\}$ and $\{v, v_2\}$ to $\mathcal H_i$ and label $v_1$ with the group $W_{\Lambda_1}$ and $v_2$ with the group $W_{\Lambda_2}$.\\

\textbf{No: }If no such $\Lambda_1$ and $\Lambda_2$ exist, then $\Lambda$ is a clique that is contained in a wide subgraph of $\Delta$,  or $\Lambda$ does not have any spherical separating subsets and is relatively wide-avoidant in $\Delta$ by \cref{rel-wide-avoidant} . In this case, we do not attach anything to the vertex $v$ in $\mathcal T_i$, i.e. $v$ will be a leaf of every hierarchy $\mathcal H_{i+m}$, $m \ge 1$.  \\

As $\Delta$ is finite, this procedure ends after finitely many steps. Let $\mathcal H$ be the tree obtained after the last step and let $n$ be its height.
\cref{check} together with \cref{emptygraphs} ensures that every $\mc H_i$ with $i \le n$ is a weak Morse hierarchy. 
Let $\Lambda$ be a graph corresponding to a leaf of $\mc H$. By construction, $\Lambda$ is either a clique that is contained in a wide subgraph of $\Delta$ or $\Lambda$ does not contain any spherical separating subgraphs and is relatively wide-avoidant by \cref{rel-wide-avoidant}. In the latter case $\Lambda$ is either spherical or one-ended by \cite[Thm 8.7.2]{Davis}. If $\Lambda$ is spherical or contained in a wide subgraph then $(\mb W_\Lambda) = \emptyset$  by \cref{emptygraphs}. If $W_\Lambda$ is one-ended and $\Lambda$ is relatively wide-avoidant in $\Delta$, then $\Lambda$ can't be contained in a wide or spherical subgraph and $(\mb W_\Lambda, W_\Delta) \neq \emptyset$ by \cref{emptygraphs}.

 If in addition, $\Delta$ does not contain any spherical nexus, then $\Lambda$ being relatively wide-avoidant implies $\Lambda$ being relatively wide-spherical-avoidant by  \cref{lem_nexus}. Then,  \cref{thmCorLev-CashenChuTao} implies that all leaf groups have connected and locally connected (possibly empty) relative Morse boundary.  
 \end{proof}

In view of \cref{rel-wide-avoidant}, the algorithm in the proof decides whether a Coxeter graph $\Delta$ contains a relatively wide-avoidant induced subgraph whose corresponding special subgroup $W_\Delta$ is one-ended: 

\begin{lem}\label{cor:Combinatorics_tot}
 Let $\Delta$ be a Coxeter graph. 
 The following are equivalent: 
 \begin{itemize}
     \item The algorithm in \cref{fig:proof_morse_accessibility} terminates in at least one graph that is not spherical and is not contained in a wide subgraph of $\Delta$;
     \item $\Delta$ contains a relatively wide-avoidant induced subgraph whose corresponding special subgroup  $W_\Delta$ is one-ended.
 \end{itemize}
If, in addition, $\Delta$ does not contain any spherical nexus, then the above is equivalent to:
\begin{itemize}
     \item $\mb W_\Delta$ contains an infinite connected component.
\end{itemize}

\end{lem}
\begin{proof}
    The claim follows from  \cref{rel-wide-avoidant} and \cref{CorCorLev}.
\end{proof}
\begin{cor}[\cref{thm:AFCoxeterTotDisconn_intro}]\label{totdiscMorseboundariesAllCoxeter}
   If every leaf group of the weak Morse hierarchy of $W_\Delta$ constructed in the proof of \cref{thmRACG_intro}  (\cref{thmRACG_}) is spherical or contained in a special subgroup of $W_\Delta$ that is wide, then $\mb W_\Delta$ is totally disconnected. Moreover if $W_\Delta$ does not contain any spherical nexus, then the reverse implication is also true. 
\end{cor}

\begin{remark} \label{rem:coxeter_leafgroups_spherical} More specifically, our algorithm recursively defines a graph class that contains Coxeter graphs whose associated groups have totally disconnected Morse boundary. This class contains all Coxeter graphs for which the algorithm in \cref{fig:proof} ends with cliques that are either spherical or contained in wide subgraphs of $\Delta$.
In the right-angled case, this graph class was introduced in \cite{karrer}. 
\end{remark}

\section{Application: Distinguishing Coxeter groups and RAAGs up to quasi-isometry} 
\label{subsec:introRACGRAAG}

In this section we illustrate how to use Morse hierarchies to distinguish groups up to quasi-isometry, with a focus on right-angled Artin groups (RAAGs) and right-angled Coxeter groups (RACGs). Note that RACGs are examples of Coxeter groups that do not contain any spherical nexus.  

In \cite{karrer}, a procedure for constructing graphs that define RACGs with totally disconnected Morse boundary was given. \cref{thm:AFCoxeterTotDisconn_intro} shows that all RACGs with totally disconnected boundary arise in this way (see \cref{rem:coxeter_leafgroups_spherical}). 

Charney--Cordes--Sisto~\cite{ChaCorSis23} proved that the Morse boundary of every right-angled Artin group (RAAG) is totally disconnected. Accordingly, if the Morse boundary of a given group $G$ is \emph{not} totally disconnected, then it cannot be quasi-isometric to RAAG. Thus, \cref{cor:Combinatorics_tot} yields:

\begin{cor}\label{cor:Combinatorics}
Let $\Delta$ be the Coxeter graph without a spherical nexus (e.g., the defining graph of a RACG). 
  If the algorithm in \cref{fig:proof_morse_accessibility} decides that $\Delta$ contains a relatively wide-avoidant induced subgraph whose corresponding special subgroup is one-ended, then $\mb W_\Delta$ contains an infinite connected component,  and $W_\Delta$ is not quasi-isometric to any RAAG. 
\end{cor}
 Another method for finding infinite connected components in Morse boundaries of RACGs is to look for stably embedded virtually closed surface RACG-subgroups \cite{Behrstockstablecycle, RAAGEDYRACGS}. Such RACG-subgroups arise, for example, from \textit{stable cycles} in the Coxeter graph, i.e., cycles of length at least 5 in which no non-adjacent vertices are contained in an induced 4-cycle. Our approach via Morse hierarchies, however, offers a more robust alternative, because the presence of a topologically embedded circle in the Morse boundary does not imply the existence of a stable cycle in the Coxeter graph \cite{GraKarLaz21}.  

Besides Morse boundaries, there are existing tools that are known to distinguish RACGs from RAAGs up to quasi-isometry such as divergence, maximal product region graphs, and JSJ-decompositions. In \cite{RAAGEDYRACGS}, definitions and an updated outline on what is known about RAAGedy RACGS is given.

We now apply our tools for Morse hierarchies we have developed to distinguish groups using connected components of the Morse boundary.  If $W_\Delta$ and $W_{\Delta'}$ are quasi-isometric Coxeter groups, their Morse boundaries are homeomorphic. Thus, the occurrence of an infinite Morse boundary component $C$ in $\mb W_\Delta$ implies the existence of an infinite Morse boundary component $C'$ in $\mb W_{\Delta'}$ that is homeomorphic to $C$. In particular, we can distinguish two Coxeter groups by identifying the existence of, or the compactness$\slash$non-compactness of any infinite connected components of their Morse boundaries. For example, if one contains a compact infinite connected component and the other does not, then the two groups cannot be quasi-isometric.

If $W_\Delta$ has no spherical nexus, then it is Morse accessible and \cref{thmRACG_intro}    (\cref{thmRACG_}) yields a corresponding Morse hierarchy. In particular, we can determine all infinite connected components of $\mb W_\Delta$ using the algorithm on the defining graph $\Delta$ illustrated in \cref{fig:proof_morse_accessibility}. The algorithm has as input the defining graph $\Delta$ and as output a set of induced subgraphs $\Lambda$ of $\Delta$. If such a graph $\Lambda$ is neither spherical nor contained in a wide subgraph, then $W_\Lambda$ is the stabilizer of an infinite Morse boundary component $C$, and $C = (\mb W_\Lambda, W_\Delta)$. In the other direction, every stabilizer of an infinite connected component of $\mb W_\Delta$ is conjugate to a special subgroup of $W_\Delta$ whose defining graph is among the induced subgraphs with which the algorithm terminates.

Moreover,  \cref{thm:point-convergence_intro} implies that an infinite Morse boundary component of $W_\Delta$  is compact if and only if its stabilizer is stable. 
 In the right-angled case, stability of special subgroups was characterized combinatorially by Tran~\cite[Thm 1.11]{Tran}, and \cite[Prop.4.9]{Genevois} and later reproved in \cite[Thm 7.5]{hierarchically}:
A special subgroup $W_\Lambda$ of a RACG $W_\Delta$ is stable if and only if $\Lambda$ does not contain any induced $4$-cycle and is {square-complete}. A subgraph $\Lambda$ of a graph $\Delta$ is  \emph{square-complete} if each induced $4$-cycle of $\Delta$ with two non-adjacent vertices lying in $\Lambda$ is fully contained in $\Lambda$. For general Coxeter groups, a characterization was recently given in Theorem 4.17 of \cite{CashenChuTao}. 

Thus, we obtain the following in the right-angled case. 
One could formulate the analogous statement for general Coxeter groups without a spherical nexus using the conditions in \cite{CashenChuTao}.

\begin{cor} \label{thmQIRACGs2}
Let  $\mb W_\Delta$  and $\mb W_{\Delta'}$  be two quasi-isometric RACGs with an infinite compact component $C$ in their Morse boundary.
Then the stabilizer of $C$ in  $ W_\Delta$ and $ W_{\Delta'}$ are both stable parabolic and quasi-isometric to each other. In particular, the corresponding induced subgraphs of $\Delta$ and $\Delta'$ are square-complete and do not contain any induced 4-cycle. 
\end{cor}

\begin{figure}[h]
    \centering
    \includegraphics[width=\textwidth]{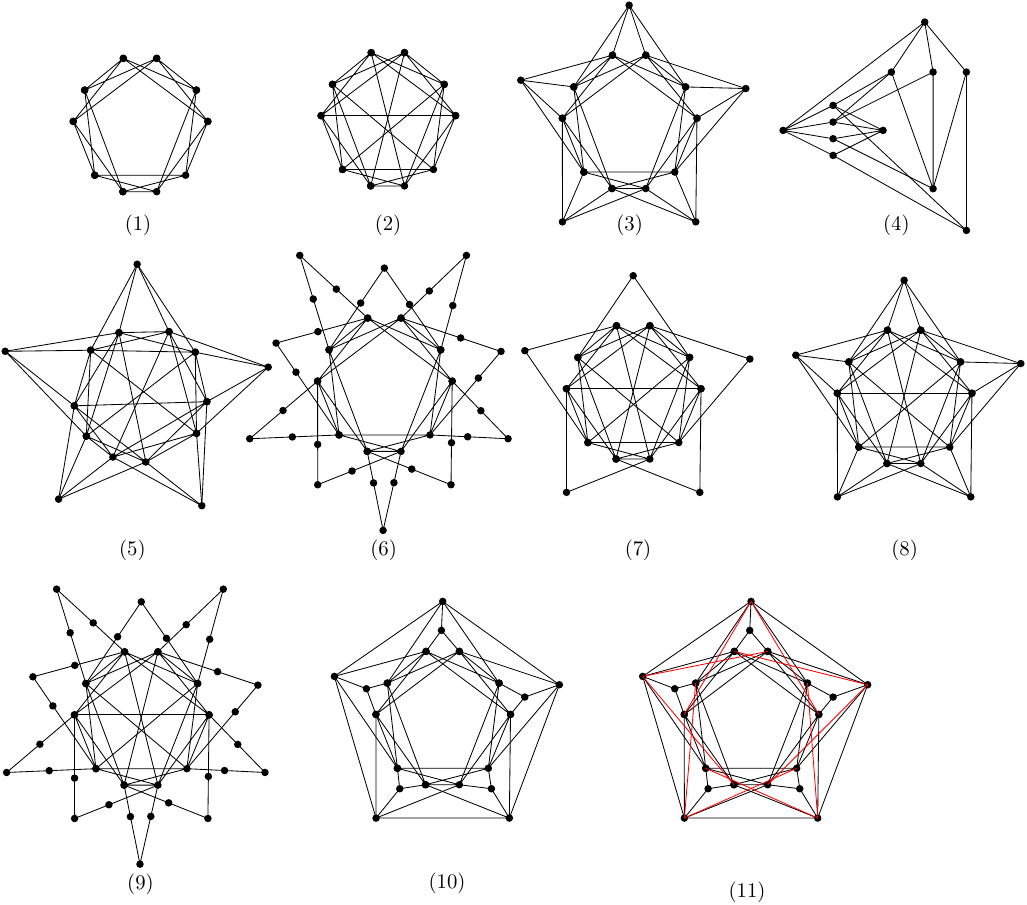}
    \caption{ Graph (1) is quasi-isometric to the 5-cycle RAAG; graph (2) arises from graph (1) through adding 6 edges; graph (3) was studied in \cite[Fig 10]{MihalikCamp}; graph (4) was discussed in \cite[Fig. 13]{RAAGEDYRACGS}; 
    graphs (5)-(11) are variations of graphs (1)-(3).}
    \label{fig:Main_examples}
\end{figure}

The following examples in \cref{fig:Main_examples} illustrate how the new quasi-isometry invariants discussed in this section are helpful for distinguishing RACGs from other RACGs or from RAAGs up to quasi-isometry. 

We briefly describe the information contained in the table. 

Column 2: The Morse boundary column indicates whether there are none, one, or infinitely many infinite non-singleton connected components along with whether those non-trivial components are compact or not. 

Column 3: The stable cycle indicates whether the defining graph contains a stable cycle or not. 

Column 4: In \cite{CorLev25}, the authors give combinatorial conditions on the defining that exactly characterize when the Morse boundary of the RACG is connected and locally connected (LC). These same conditions appear in the earlier paper \cite{MRT07}, where the authors investigate when the visual boundaries for a given RACG are locally connected (LC) versus not locally connected (NLC). We record this in the visual boundary column. They show that if one of these conditions is satisfied, then any visual boundary for the corresponding RACG is not locally connected (NLC). In \cite{MihalikCamp}, Mihalik and Camp construct graph (3) to show that the combinatorial conditions in \cref{thmCorLev} are sufficient, but not necessary for local connectedness of the visual boundary. Specifically, the RACG group with graph (3) as defining graph has LC visual boundary but totally disconnected  Morse boundary.

We now focus on the last three columns. For RACGs, quadratic divergence is equivalent to a combinatorial condition on the defining graph called CFS \cite{DanTho15},\cite{Lev22}. Since quadratic divergence is a sufficient but not necessary condition for a RACG to be quasi-isometric to a RAAG, there has been other work developing finer invariants. In the triangle free case, the definition of CFS was strengthened to \emph{strongly} CFS as a necessary condition for a RACG to be RAAGedy \cite{RAAGEDYRACGS}. 

\begin{center}
\begin{tabular}{ |c|c|c|c|c|c|c| } 
 \hline
graph & Morse   & stable & visual  & RAAG & $\triangle$-free & CFS\\ 
\ & boundary &  cycle & bdry & -edy & & \\
 \hline
1 & tot. disc& no &NLC &yes&yes  & strong\\ 
2 & conn.  &no &LC&no&no & strong\\ 
3 & tot. disc. &no &LC&no&no& no\\ 
4 & $\infty$ non-cpt  &no & NLC&no&yes& strong\\ 
5 &  $\infty$ non-cpt   &no &NLC &no&no&no\\ 
6 & $\infty$ cpt  &yes &NLC &no&yes & no\\ 
7 & $\infty$ non-cpt   &no & NLC&no&no & strong\\ 
8 & $\infty$ non-cpt  &no & LC &no&no & no\\ 
9 & $\infty$ cpct \& non-cpct   &yes & NLC &no&yes & no\\ 
10 & $\infty$ non-cpt  &yes & LC&no &yes & strong\\ 
11 & tot. disc  &no & NLC&?&yes & strong\\ 
 \hline
\end{tabular}
\end{center}

We list the \textbf{new} insights through the study of Morse hierarchies that, according to our knowledge, \emph{cannot} be obtained by other methods here: 
\begin{itemize}
    \item (2) is not quasi-isometric to (7), (4), (10)
    \item (3) is not quasi-isometric to (5), (6), (8)
    \item (4) is not quasi-isometric to (6), (11)
    \item (5) is not quasi-isometric to to (6) 
    \item (6) is not quasi-isometric to (8), (9)
    \item (7) is not quasi-isometric to a RAAG and not quasi-isometric to (1) or (11)
\end{itemize}

Note that we can also apply \cref{thm:point-convergence_intro} and \cref{corSurfaces_intro} to the graphs above, meaning that in all groups quasi-isometric to the examples above, the infinite Morse boundary components arise from their stabilizers. Moreover, all compact components in such groups have stable stabilizers.  In the case of circle components in the boundary, as for instance in graph (6), every group quasi-isometric to that group contains a stable subgroup that is virtually a hyperbolic surface group (see \cref{subsec:exampes_RACG}).

For many more examples that illustrate Morse accessibility, including an illustration of the mentioned algorithm in \cref{fig:proof_morse_accessibility} for RACGs, see \cref{sec:examples}.

\section{Examples} \label{sec:examples} 
\subsection{Morse accessible groups} \label{subsec:list_Morse_accessible}
In this section, we list examples of Morse-accessible groups and outline how to construct Morse hierarchies for them.  We begin with the most basic example of finitely generated free group so that we can illustrate how we will describe Morse hierarchies with notation and figures. 

\begin{example}[Morse hierarchies of free groups] \label{ex:freegroup}
\begin{figure}[h]
    \centering
    \includegraphics[]{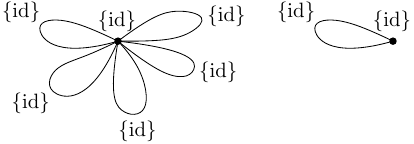}
    \caption{ Two Morse splittings of $\mathbb F_5$ and $\mathbb Z$ respectively. }
    \label{fig:freeHNN}
\end{figure}
 For $\mathbb F_n$ with $n\ge 1$, we obtain a Morse splitting by viewing $\mathbb F_n$ as the fundamental group of a graph of groups with one vertex labeled with the trivial group and $n$ loops where each represents an HNN-extension of the trivial group as in \cref{fig:freeHNN}. In particular, it will be useful to view $\mathbb Z = \{id\}_{\{id\}}$ for extending weak Morse hierarchies with $\mathbb Z$ leaf groups to Morse hierarchies.

When illustrating Morse hierarchies, we will often label the vertices of the hierarchy tree with spaces whose fundamental groups yield the corresponding vertex groups, rather than the groups themselves. The associated graph of groups decomposition can then be indicated by coloring the edge spaces in red. Decomposing the space along these red subspaces yields pieces whose fundamental groups are precisely the vertex groups of the corresponding graph of groups decomposition.   For instance, we can illustrate the Morse hierarchy on the left of \cref{fig:freeHNN} alternatively as in \cref{fig:freeHNN_topol}.
\begin{figure}[h]
    \centering
    \includegraphics[]{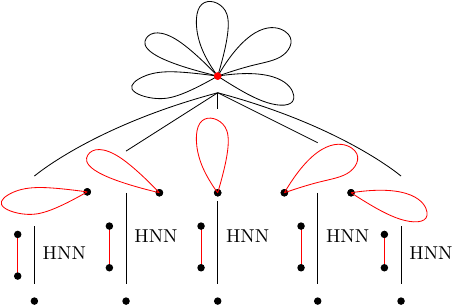}
    \caption{A Morse hierarchy of height 2 of $\mathbb F_5$ illustrated via decompositions of associated spaces.}
    \label{fig:freeHNN_topol}
\end{figure}

In \cref{fig:freegroup}, two  Morse hierarchies for  $\mathbb F_5$ are illustrated where each vertex of the pictured tree is labeled by the corresponding graph of groups. In the hierarchy on the left, we first split $\mathbb F_5$ into $\mathbb Z$ factors, and then we write each such factor as $\mathbb Z = \{id\}_{\{id\}}$. On the right, we split  $\mathbb F_5$ successively into $\mathbb Z$ factors, yielding a hierarchy of height $3$.

\begin{figure}[h]
    \centering
    \includegraphics[]{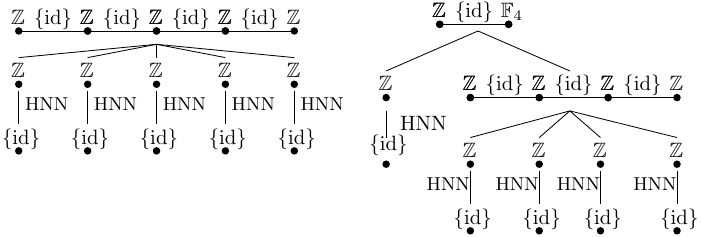}
    \caption{Two Morse hierarchies of $\mathbb F_5$.}
    \label{fig:freegroup}
\end{figure}
\end{example} 

\begin{example}[Hyperbolic groups]
    If $\Gamma$ is a one-ended hyperbolic group, then the Morse boundary of $\Gamma$ coincides with the Gromov boundary of $\Gamma$ and thus we have a trivial Morse hierarchy of height 0.  In the multi-ended case, we apply the Dunwoody accessibility result as in the introduction to obtain a Morse hierarchy. 
\end{example}
\begin{example}[Non-hyperbolic group with height 0 Morse hierarchy]
    If $G$ is a finitely generated group with connected Morse boundary, then the group itself yields a Morse hierarchy of height $0$. A non-hyperbolic example of such a group is a cusped hyperbolic 3-manifold group whose Morse boundary is an $\omega$-Sierpinski curve \cite{ChaCorSis23}.
\end{example}

\begin{example}[Free products]
    Suppose a finitely generated group $G=A_1\star \cdots \star A_n$ splits as a free product. In \cite{Zbinden-freeproducts}, Zbinden proves that $\mb G$ is determined by the Morse boundaries of the factors. Since each $A_i$ is undistorted in $G$, and because every edge group has empty relative Morse boundary, this yields a \emph{weak} Morse splitting of $G$. If, in addition, each $A_i$ has connected relative Morse boundary in $G$, the splitting is Morse. 
\end{example}

\begin{example}[Morse Splittings and  hierarchies of CAT(0) groups]\label{ex_BH}
Many Morse splittings and hierarchies of CAT(0) groups can be obtained by constructing the group iteratively via HNN-extensions and amalgamated free products, checking that the amalgams used satisfy the assumptions of \cite[II.Thm 11.18, II.Thm 11.21  ]{BridsonHaefliger}. Hierarchies constructed this way are: 
\begin{enumerate}
    \item \label{item:cat(0)-1} Weak Morse, if one ensures that every edge group has empty relative Morse boundary.  The easiest way to ensure this is to construct splittings along groups that have empty Morse boundary such as finite groups or direct products of infinite finitely generated groups. 
    \item \label{item:cat(0)-2}Morse, if in addition to (\ref{item:cat(0)-1}), every leaf group has connected relative Morse boundary. One way to ensure this is to build Morse hierarchies in which the leaf groups have connected Morse boundaries that topologically embeds into the ambient group.  The latter is satisfied for leaf groups that are \textit{stable} (see \cref{defstable}) or  \textit{Morse} (see \cref{def:MorseSubgroup}).
\end{enumerate}

For a specific example, we can glue two hyperbolic 3-manifolds $M_1$ and $M_2$ with cusps along one of their torus boundaries $T$. The corresponding fundamental group has a Morse splitting that is obtained by cutting the space along  $T$.  In fact, one can show that $\pi_1(M_1)$ and $\pi_1(M_2)$ are Morse in $G$, and  by  \cite{ChaCorSis23}, the Morse boundaries of  $\pi_1(M_1)$ and $\pi_1(M_2)$ are $\omega$-Sierpinski-curves, and as such non-compact and connected.

However, (\ref{item:cat(0)-1}) and (\ref{item:cat(0)-2}) can also be satisfied without the leaf groups being stable or Morse and the edge groups having empty Morse boundary. We illustrate this in \cref{subsec:Surfaces} and \cref{subsec:comp}. 
\end{example}
 \begin{example}[3-manifold groups]
 \label{3mf}  
In \cite{Zbinden}, Zbinden calculates the Morse boundaries of closed, orientable 3-manifold groups. In the proof of Theorem 7 of that paper, Zbinden implicitly outlines how to obtain a Morse hierarchy for these groups.  Briefly, every closed, orientable 3-manifold $M$ can first be cut along spheres (prime decomposition) and then along tori (JSJ decomposition) such that the remaining pieces are well-understood (Geometrization Theorem). 
Since spheres are finite, and the fundamental groups of tori have empty Morse boundary, this yields a hierarchy of height at most 2.  The Morse boundary and relative Morse boundary of each leaf group is connected, so the considered hierarchy is Morse. More precisely, the Morse boundary and relative Morse boundary of a leaf group is either empty, or a 2-sphere (and the corresponding piece is a closed hyperbolic 3-manifold) or an $\omega$-Sierpinski curve (and the corresponding piece is a hyperbolic 3-manifold with cusp, \cite{ChaCorSis23}).  
\end{example}

\begin{example}\label{ex:RACGsshort}(Right-angled Artin Groups)
A right-angled Artin group (RAAG) $A_\Delta$ is obtained from the RACG $W_\Delta$ by removing the condition that every generator is an involution.  Charney--Cordes--Sisto~\cite{ChaCorSis23} showed that the Morse boundary of every RAAG is totally disconnected. An alternative proof was given in \cite[Thm 7.6]{karrer} via implicitly constructing a weak Morse hierarchy. This weak Morse hierarchy is obtained by successively deleting vertices from the defining graph until one obtains a graph consisting of disconnected cliques. The weak Morse hierarchy can be extended to a Morse hierarchy by writing $\mb \mathbb Z$ as trivial HNN extension $\mb \mathbb Z= id*_{id}$ whenever a clique we end up is a single vertex.  In particular, every RAAG is Morse-accessible. Two such Morse hierarchies are illustrated in the middle and on the right in \cref{fig:MaxAndMin_Morse_hierarchies_RAAGexample}.
\end{example}

\begin{example}[Coxeter groups]  
In \cref{sec:MorseAccessAFCoxeter}, we show that every Coxeter group has a weak Morse hierarchy, which is a genuine Morse hierarchy for Coxeter groups that do not contain a spherical nexus (as for instance in the right-angled case).  The proof is based on an algorithm on the defining graph illustrated in \cref{fig:proof_morse_accessibility}. In  \cref{subsec:introRACGRAAG} we explain how to use this  algorithm to identify all infinite Morse boundary components of a Coxeter group that does not contain a spherical nexus, and how to use it in order to distinguish Coxeter groups up to quasi-isometry.  
\end{example} 

\subsection{Non-uniqueness of Morse hierarchies}\label{subsec:nonunique}
  \begin{figure}[h]
    \centering
    \includegraphics[]{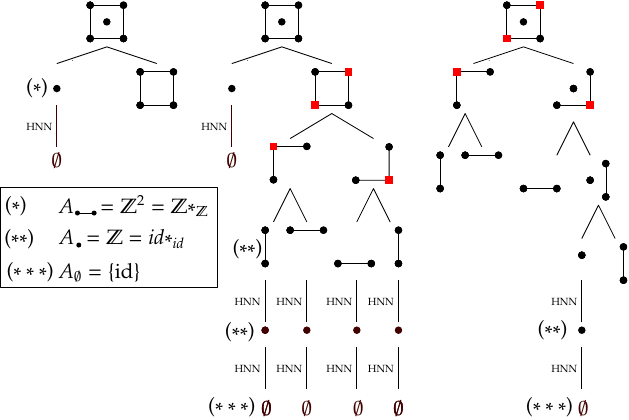}
    \caption{Three Morse hierarchies of a RAAG $A_\Delta = \mathbb Z * (\mathbb F_2 \times \mathbb F_2)$ with $\Delta$ being the union of a 4-cycle and a single vertex.}
    \label{fig:MaxAndMin_Morse_hierarchies_RAAGexample}
\end{figure}
 We have already illustrated the non-uniqueness of Morse hierarchies using $\mathbb F_n$. In this subsection, we give a concrete example of a non-hyperbolic RAAG with  different Morse hierarchies. 

We study the disjoint union of a $4$-cycle and a vertex. This RAAG is isomorphic to $\mathbb Z * \mathbb (\mathbb F_2 \times \mathbb F_2)$. 
The three trees in \cref{fig:MaxAndMin_Morse_hierarchies_RAAGexample} yield three different Morse hierarchies for the corresponding RAAG: Whenever a vertex labeled by a graph $\Delta$ has two descendants $\Delta_1$ and $\Delta_2$, the RAAG $A_\Delta$ can be written as $A_\Delta=A_{\Delta_1}*_{A_{\Delta_1 \cap \Delta_2}} A_{\Delta_2}$. Whenever a vertex has only one descendant, the corresponding RAAG is isomorphic to $\mathbb Z$ which we write as HNN extension $\mathbb Z = id_{id}$. In the hierarchy on the left, we stopped splitting the graph earlier via noticing that the RAAG corresponding to a 4-cycle is the direct product of two free groups of rank 2. As such, its Morse boundary and relative Morse boundary is empty, and we do not need to split this group further.  
Observe that the edge groups in the three Morse hierarchies are isomorphic to $\mathbb Z$, $\mathbb F_2$  or $id$. Note that the Morse boundaries of $\mathbb Z$ and $\mathbb F_2$ are not empty. However, their \textbf{relative Morse boundaries} are empty as they lie in direct products of infinite groups.

\subsection{Comparing  Morse hierarchies of two specific surface amalgams}\label{subsec:Surfaces}
In this subsection compare the Morse hierarchies of two specific surface amalgams in order to illustrate some features and technicalities of Morse hierarchies. 
   \begin{figure}[h]
    \centering
    \includegraphics[width=0.4\textwidth]{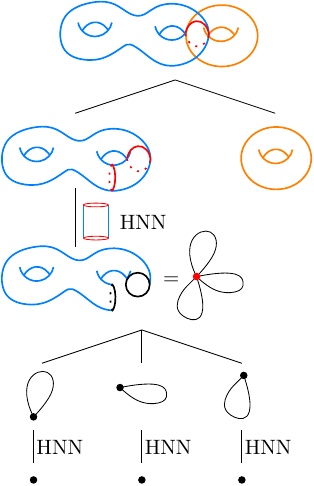}
    \caption{The pictured tree yields a Morse hierarchy for the fundamental group of the space at the top of the picture. }
    \label{fig:surfaceamalgam1}
\end{figure}

 First, we study the space $X$ obtained by gluing a surface $S$ of genus two to a torus $T$ along a simple closed curve $\alpha$ as in \cref{fig:surfaceamalgam1}. 
 In order to create a Morse hierarchy for $\pi_1(X)$, we first decompose $S$ and $T$ along $\alpha$. At the group level, this corresponds to a splitting of $\pi_1(X)$ along an undistorted $\mathbb Z$ subgroup. This $\mathbb Z$ has empty relative Morse boundary in $\mb \pi_1(X)$ because $\alpha$ is embedded in $T$, i.e. the lift of $\alpha$ in $\tilde X$ lives in a flat subspace and is therefore not Morse in $\tilde X$.  

 Observe that although $\mb \pi_1(S)$ is a circle, it's \emph{relative} Morse boundary in $\pi_1(X)$ is totally disconnected (see section 6 in \cite{Charney-Sultan}), thus we must continue to decompose $S$ in order to obtain a Morse hierarchy. We cut $S$ along $\alpha$ and a closed curve $\beta$ parallel to $\alpha$ which allows us to write $\pi_1(S)$ as HNN extension along the curve $\alpha$. This is a valid splitting as the lifts of $\alpha$ and $\beta$ are not Morse in $\tilde X$.  After decomposing $S$ this way, we obtain a space that is homotopy-equivalent to the wedge of three circles and thus has fundamental group $F_3$ which we can decompose as in \cref{ex:freegroup}. 
Accordingly, the described decomposition of $X$ yields a Morse hierarchy of $\pi_1(X)$. In this Morse hierarchy, each leaf group has empty relative Morse boundary, implying that $\partial_*\pi_1(X)$ is totally disconnected by \cref{corMorseaccessiblegroups_Intro}.
   \begin{figure}[h]
    \centering
    \includegraphics[]{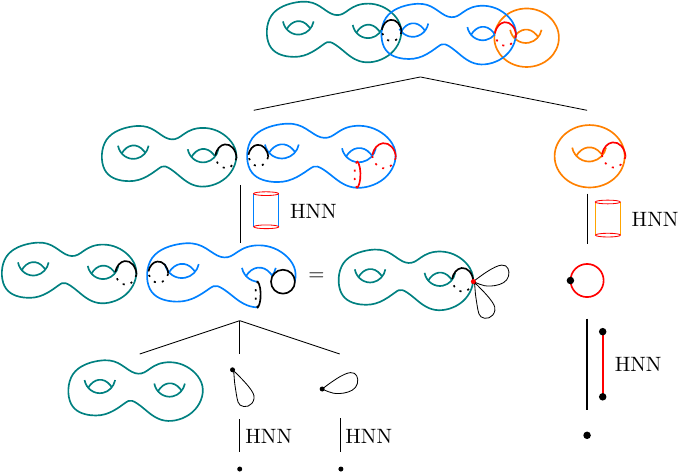}
    \caption{The pictured tree yields a Morse hierarchy for the fundamental group of the space at the top of the picture.}
    \label{fig:surfaceamalgam2}
\end{figure}

Next, we consider a space $Y$ obtained by gluing two surfaces $S_1$ and $S_2$ of genus two to a torus $T$ as seen in \cref{fig:surfaceamalgam2}.  We can decompose $Y$ and $\pi_1(Y)$ similarly to how we did in the first example. However, this time, we decompose the spaces occurring in the hierarchy further whenever possible.  In particular, we do not keep the torus $T$ as a leaf, but we decompose the torus further. On a group theoretic level this means to write $\mathbb Z^2=\mathbb Z *_\mathbb Z$. This way, we  obtain a Morse hierarchy with two different kind of leaf groups: the trivial group, and $\pi_1(S_1)$. The Morse boundary of $\pi_1(S_1)$ is a circle. Since $\pi_1(S_1)$ turns out to be a stable subgroup of $\pi_1(Y)$, the Morse boundary of $\pi_1(S_1)$ is topologically embedded in the Morse boundary of $Y.$ 
Accordingly, $\pi_1(S_1)$ is a leaf group with non-empty connected relative Morse boundary, and stabilizes an infinite component of $\partial_*\pi_1(Y)$ homeomorphic to a circle by \cref{thmMorseaccessiblegroups}. 
Hence, $\pi_1(Y)$ is not quasi-isometric to $\pi_1(X)$, as the latter has totally disconnected Morse boundary.  Moreover, \cref{corSurfaces_intro} yields that every group quasi-isometric to $\pi_1(Y)$ has a stable subgroup that contains as a subgroup of finite index the fundamental group of a closed hyperbolic surface. 

\begin{remark} \label{rem_relMor}
    The discussed hierarchies of $\pi_1(X)$ and $\pi_1(Y)$ show that it is crucial to allow splittings along groups that have \textit{non-empty Morse boundary} but empty \textbf{relative Morse boundary} in the ambient group. Moreover, neither $\pi_1(X)$ nor $\pi_1(Y)$ split as a \textbf{graph of groups} that is \emph{Morse}, i.e. in which every vertex group has connected relative Morse boundary and in which every edge group is undistorted and has empty relative Morse boundary. Accordingly, the study of hierarchies rather than graph of groups is necessary in order to understand these examples.
\end{remark}

\subsection{Comparing  Morse hierarchies of two specific RACGs}\label{subsec:exampes_RACG}

In this subsection we show explicitly how the algorithm outlined in \cref{fig:proof_morse_accessibility} in the proof of \cref{thmRACG_intro} 
 produces Morse hierarchies from the defining graph of a RACG and how it can be used to distinguish RACGs up to quasi-isometry. In order to do so, we apply the algorithm in different ways to two RACGs and compare the Morse hierarchies we obtain this way. We highlight the different features and technicalities that might occur along the way. 

   \begin{figure}[h]
    \centering
    \includegraphics[width=0.75\textwidth]{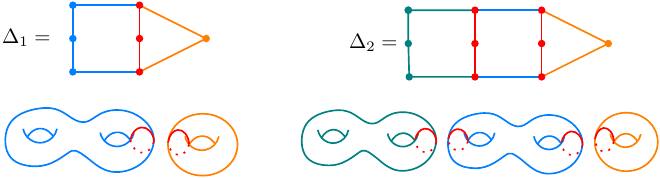}
    \caption{Two examples inspired from Section 6 of \cite{Charney-Sultan}. The groups $W_{\Delta_1}$ and $W_{\Delta_2}$ are quasi-isometric to the fundamental groups of the spaces pictured below the graphs $\Delta_1$ and $\Delta_2$ respectively. }
    \label{fig:Main_examplesCharSul}
\end{figure} 

We apply the algorithm to two toy examples, $\Delta_1$ and $\Delta_2$ in
\cref{fig:Main_examplesCharSul}. Note that these groups are analogues of the groups in \cref{subsec:Surfaces}. The algorithm will show that $\mb W_{\Delta_1}$ is totally disconnected and $\mb W_{\Delta_2}$ has an infinite connected component just as in the corresponding surface amalgams. In particular, the Morse boundaries of $W_{\Delta_1}$ and $W_{\Delta_2}$ are distinct and hence the groups $W_{\Delta_1}$ and $W_{\Delta_2}$ are not quasi-isometric. 

\begin{figure}[h]
    \centering
    \includegraphics[width=\textwidth]{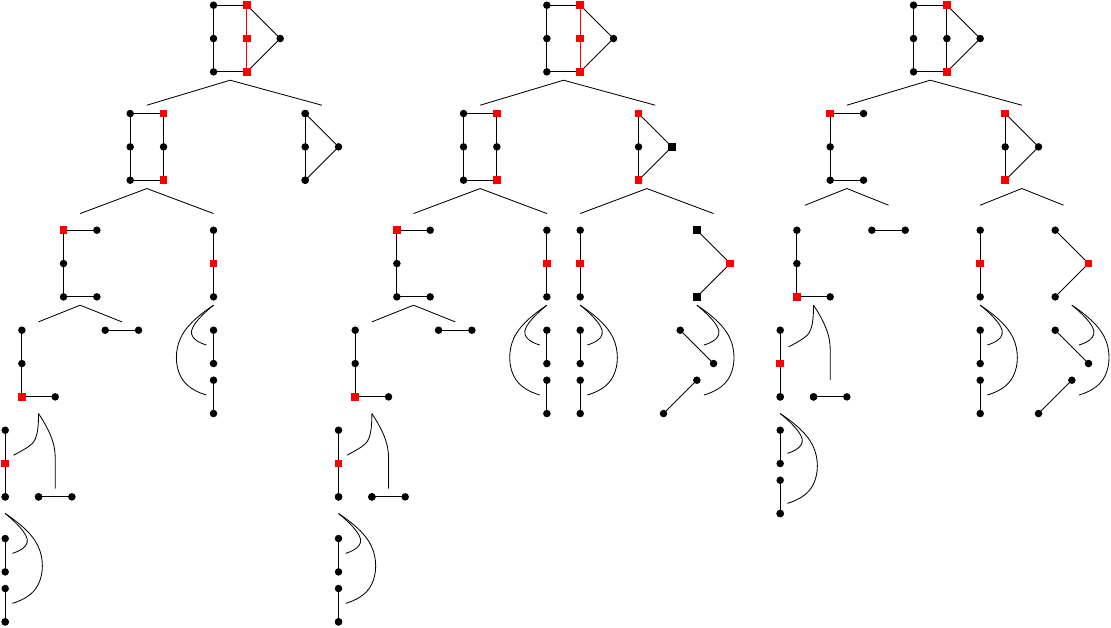}
    \caption{Three trees corresponding to two Morse hierarchies of $W_{\Delta_1}$ with finite leaf groups.}
    \label{fig:two_nice_Max_Morsetrees_one_storey_houses1}
\end{figure} 

The three trees in \cref{fig:two_nice_Max_Morsetrees_one_storey_houses1}, give three different Morse hierarchies of $W_{\Delta_1}$. Analogous to the RAAG example in \cref{subsec:nonunique}, each of the three trees illustrates how the graph $\Delta_1$ can be successively decomposed into a pair of proper induced subgraphs whose intersection graph (each time drawn in red) corresponds to a RACG-subgroup with empty relative Morse boundary. 

We now compare the three Morse hierarchies of $W_{\Delta_1}$, first focusing on the left and middle trees. In both of these trees, we decompose the graph along a two-path in the first step resulting in two subgraphs: a 6-cycle $C_6$ and a 4-cycle $C_4$.  These are valid splittings because the two path in each graph is contained in an induced 4-cycle and thus the corresponding edge group has empty relative Morse boundary. In the graph on the left, we keep the 4-cycle as is so that it corresponds to a leaf group of the Morse hierarchy.  

However, if we follow the algorithm until the end, we are forced to decompose the 4-cycle $C_4$ further until we end up with cliques. All further splittings are valid because every induced subgraph of $C_4$ corresponds to a special subgroup with empty relative Morse boundary. This is illustrated in the middle tree. 

We now consider the 6-cycle $C_6$. Note that it does not contain any separating clique or induced 4-cycle. However, the two end vertices of the two-path along which we decomposed $\Delta_1$ in the first step are contained in a 4-cycle of $\Delta_1$. Hence, the corresponding special subgroup has empty relative Morse boundary in $W_{\Delta_1}$, and so we can decompose the graph along these two vertices. This results in a path of length 4 and the two-path along which we decomposed the graph in the first step. Now, we successively decompose these paths along single vertices until we end up with cliques of size 2. 

{\bf Note:} Each edge of $\Delta_1$ occurs once or twice as a leaf of both trees. The edges occurring twice are the edges in the two-path along which we decomposed $\Delta_1$ in the first step. This two-path is a subgraph of both graphs obtained after the first step of the decomposition and hence is decomposed twice. Such inconvenience can happen when decomposing a graph along cliques or subgraphs of thick joins. In the tree on the right, we apply the algorithm to $\Delta_1$ differently to avoid this inconvenience described above. 

All three hierarchies share the following crucial properties: Every edge of $\Delta_1$ occurs at least once as a leaf of the tree. Only the number of occurrences per edge differs. All leaf groups correspond to cliques, implying that the Morse boundary of $W_{\Delta_1}$ is totally disconnected by \cref{thm:AFCoxeterTotDisconn_intro}. 

   \begin{figure}[h]
    \centering
    \includegraphics[width=\textwidth]{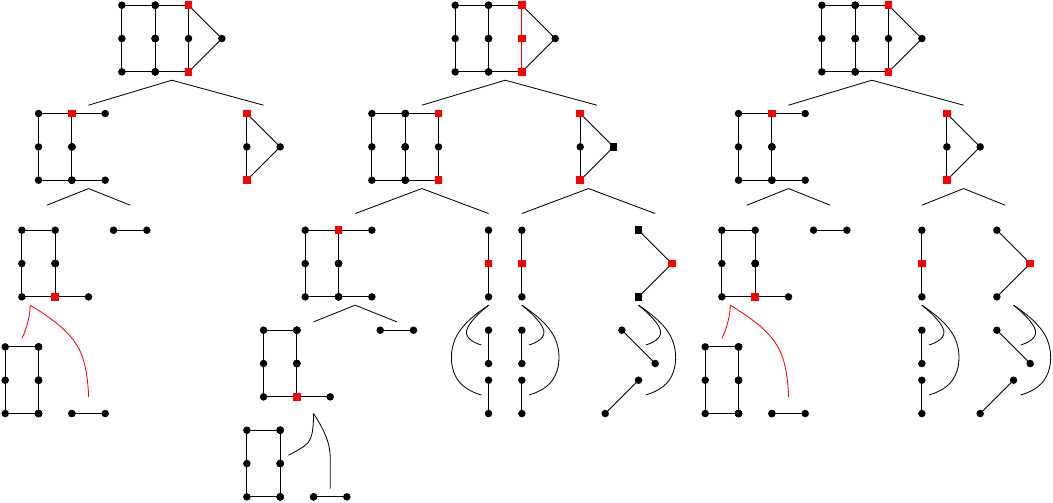}
    \caption{Three trees that yield Morse hierarchies of $W_{\Delta_2}$ in which a leaf group is virtually a surface group.}
    \label{fig:two_nice_Max_Morsetrees_double_houses}
\end{figure}

Next, we examine $W_{\Delta_2}$. \cref{fig:two_nice_Max_Morsetrees_double_houses} shows three trees that correspond to three different Morse hierarchies.  In contrast to $\Delta_1$, these hierarchies have a leaf that corresponds to a $6-$cycle $C_6$. 
Note that $W_{C_6}$ is virtually a hyperbolic surface group whose Morse boundary is a circle. Thus the corresponding 
$W_{C_6}$ subgroup is the stabilizer of an infinite Morse boundary component in $\mb W_{\Delta_2}$. As in the torsion-free analogue  of $W_{\Delta_2}$  in \cref{fig:surfaceamalgam2}, the cycle topologically embeds in the Morse boundary of the ambient group because $W_{C_6}$ is stable. This has strong consequences for the defining graph of any RACG quasi-isometric to $W_{\Delta_2}$. 

Indeed, by \cref{thmQIRACGs2}, if $W_{\Delta'}$ is a RACG quasi-isometric to $W_{\Delta_2}$, then $\mb W_{\Delta'}$ will also contain a circle component whose stabilizer is a  conjugate to a special subgroup of $W_{\Delta'}$. Let $\Delta$ be the corresponding induced subgraph of $\Delta'$. \cref{thmQIRACGs2} also implies that $\Delta$ is square-complete, does not contain an induced $4$-cycle, and that $W_{\Delta}$ is quasi-isometric to $W_{C_6}$. In particular, $W_{\Delta}$ is hyperbolic with Gromov boundary a circle. 
Now by work of \cite{Tukia, CassonJungreis, Gabai} $W_{\Delta}$ is also virtually a closed surface group. Finally, by \cite{DaniThomasJSJ},  $\Delta$ contains an induced cycle of length at least $5$ in which no pair of non-adjacent vertices is contained in an induced 4-cycle in  its defining graph.\\

We observe that the Morse boundaries of  $W_{\Delta_1}$ and $W_{\Delta_2}$ are not compact, but all their connected components are. In the next subsection, we'll consider examples in which non-compact connected components arise. 

\subsection{Non-compact connected components}\label{subsec:comp}
Recall that by \cref{thm:point-convergence_intro}, an infinite Morse boundary component is compact if and only if its stabilizer is stable.  It is naturally to ask if a similar characterization can be given for infinite Morse boundary components that are not compact.
In order to give an intuition for this question, we first consider an easy example of a RACG with infinite non-compact Morse boundary components. 
\begin{example}[A RACG with non-compact Morse boundary components] \label{ex:RACGMorseComp}
The RACG whose defining graph is the one-skeleton of a 3-dimensional cube is a virtual hyperbolic 3-manifold group with torus cusp. 
Thus we can obtain a  RACG example analogous to the 3-manifold example in \cref{ex_BH} by identifying two copies of this along of one of its 4-cycle faces. Thus the Morse boundary of this amalgam consists of single-point components and components homeomorphic to $\omega$-Sierpinski curves.

In general, whenever a hierarchy for a RACG produced by the algorithm  in \cref{fig:proof_morse_accessibility} has a leaf group whose corresponding induced subgraph contains an induced 4-cycle and is square-complete, the corresponding infinite Morse boundary component is not compact and arises from a special subgroup that is Morse, but not stable. This can be concluded from Tran's characterization of stability and Morseness of special subgroups \cite{Tran}.
\end{example}
In view of this example one might hope that every non-compact infinite connected component arises from a group that is Morse but not stable. However, this turns out not be the case. The next example shows that stabilizers of infinite Morse boundary components can behave less nicely than one might expect. In this example, an infinite non-compact Morse boundary component $C$ has a hyperbolic stabilizer $H$ that is not Morse. Moreover, $\mb H$ is not homeomorphic to $C$. 
 \begin{figure}[h]
    \centering
    \includegraphics[width=0.3\textwidth]{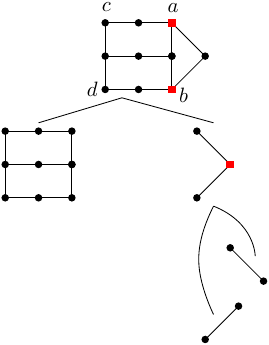}
    \caption{The pictured tree yields a Morse hierarchy.}
    \label{fig:non-Morse-components}
\end{figure}

\begin{example}[Hyperbolic RACG leaf group that is not Morse]
\label{ex:topologies different} 
    Let $\Delta$ be the root of the tree in \cref{fig:non-Morse-components} and $\Lambda$ be the only leaf of the pictured tree that is not a clique.   
    $\Lambda$ is relatively wide-avoidant in $\Delta$, so $(\mb W_\Lambda, W_\Delta)$ is connected, thus the pictured tree yields a Morse hierarchy. Accordingly, $(\mb W_\Lambda, W_\Delta)$ is an infinite connected component of $\mb W_\Delta$, and $W_\Lambda$ is its stabilizer by \cref{thmMorseaccessiblegroups_Intro}. Since $\Lambda$ does not contain any induced 4-cycle, and it is not square-complete,  $W_\Lambda$ is hyperbolic, but not Morse. 

  This also gives an example where $\iota$ induces a bijection $\iota_* \colon \mathcal M(Y,X) \to  (\mb \iota(Y), X)$ that is closed and open, but not continuous. 

     To see this, let $\Upsilon$ be the induced subgraph of $\Lambda$ on the vertices $\{a,b,c,d\}$. We note that this subgraph has no edges, so $W_\Upsilon$ is virtually free and its corresponding Davis complex is isometrically embedded in the Davis complex of $W_\Delta$. In particular, all infinite reduced words in $W_\Upsilon$ are distinct geodesic rays in the Davis complex of $W_\Delta$. 

    Now consider the geodesic rays $\gamma_i=(cd)^i(ab)^i(cd)^\infty$. Every $\gamma_i$ only intersects the flat determined by the unique $4$-cycle containing $a,b$ in $W_\Delta$ in $(ab)^i$. Thus each $\gamma_i$ is $M_i$-Morse, where $M_i$ depends on $i$ and $M_i < M_{i+1}$. First, this shows that $W_\Lambda$ is \emph{not} Morse in $W_\Delta$. Secondly, this means that $\gamma_i$ does \emph{not} converge in $\mb W_\Delta$ and thus does not converge in $(\mb W_\Lambda, W_\Delta)$. Now since $W_\Lambda$ is a hyperbolic group, $\gamma_i$ does converge to the geodesic ray $(cd)^i$, which is in the Morse (= Gromov) boundary $\mb W_\Delta$ of $W_{\Delta}$, and thus converges in the topology on $\mathcal M(\mb W_\Lambda,W_\Delta)$.
    We remark that the same line of argument used here also works for the example in \cite{GraKarLaz21} in which the defining graph studied there does not contain any stable cycle. 
\end{example}

\section{Alternative definitions of Morse-accessibility}\label{sec:AlternativeMorse} \label{rem_varMorse-acc}

The concept of \textit{Morse accessibility} was informally introduced in \cite{karrer2024subgroups}. 
Our definition here is just one way to formalize the notion that worked for our purposes.  We offer two alternative ways one could have defined Morse accessibility that might provide other interesting results.

Recall that in our definition of Morse hierarchy all edge groups and vertex groups are undistorted, and all leaf groups have connected relative Morse boundary. Alternatively keeping the undistorted condition, one could assume

\begin{enumerate}
    \item all leaf groups have connected Morse boundary, or 
    \item all leaf groups have both, connected Morse boundary and connected relative Morse boundary.
\end{enumerate} 

 Let $L$ be a leaf group of a weak Morse hierarchy. If the Morse boundary of $L$ is empty (and hence connected), then the relative Morse boundary is also empty (and hence connected). However, the other implication is not always true. There are numerous examples of weak Morse hierarchies in which a leaf group has an empty relative Morse boundary, yet its Morse boundary is non-empty and disconnected. 
 
 As a simple example for right-angled Coxeter groups, let $\Delta$ be $K_{3,2}$, i.e. the suspension over the graph of three distinct vertices. Suppose $a,b$ are the suspension points and $x,y,z$ are the three other vertices. Then we can view the right-angled Coxeter group on $W_\Delta=A \ast_C B$ as the free product with amalgamation of $A=\langle a,b,x,y\rangle$ (a 4-cycle) with $B=\langle a,b,z\rangle$ (a 2-path) over $C=\langle a,b\rangle$. In this case, $C$ is two-ended and thus has non-empty Morse boundary, however $C$ has empty relative Morse boundary in the ambient group because $a,b$ live in the 4-cycle graph $a,b,x,y$ which gives a virtually $\mathbb{Z}^2$ subgroup. In this case, $W_\Delta$ has empty Morse boundary, however you can carefully attach any graph $\Delta'$ to $\Delta$ along the suspension points to obtain a right-angled Coxeter group with non-empty Morse boundary and $\langle a,b\rangle$ still has the desired properties.

 However, in all examples we know, such hierarchies can be extended to hierarchies in which each leaf group with empty relative Morse boundary has empty Morse boundary as well. Thus $L$ has both a connected relative Morse boundary and a connected Morse boundary. \\ 

 If the Morse boundary  of the leaf group $L$ is not empty yet is connected, then its relative Morse boundary might not be connected. Again, there are numerous examples of weak Morse hierarchies in which a leaf group has a connected Morse boundary, yet its relative Morse boundary is non-empty and disconnected. A good example of this behavior is found in \cref{fig:Main_examplesCharSul}. 
 
 However, in all examples we know, the other implication is satisfied. In all weak Morse hierarchies we know, every leaf group with connected non-empty relative Morse boundary has Morse boundary that is connected and non-empty as well. 

These observations suggest that the answer to the following question might be positive:

\begin{question}
    Suppose every leaf group of a weak Morse hierarchy has connected relative Morse boundary. Can this hierarchy be extended to one in which every leaf group has both, a connected relative Morse boundary and a connected Morse boundary? 
\end{question}

A positive answer of this question would suggest that the last condition is the most natural choice for the notion of Morse accessibility. Nevertheless, for the purposes of the proofs and the underlying theory, it is necessary and sufficient only to assume that the leaf groups have connected relative Morse boundary. Therefore, we have chosen to require only this assumption.

One could also loosen the restrictions on the undistortedness of the leaf and/or edge groups in a Morse hierarchy, but there be dragons. Indeed, there would need to be many new results on splittings over undistorted distorted groups and its relationship to the Morse boundary. Indeed, there is an open question that would need to be answered first:

\begin{question}
    Does there exist infinite groups $H<G$ where $\mb H = \emptyset$ but $(\mb H, G)$ is non-empty?
\end{question}

Such an $H$ cannot be contained in an undistorted subgroup with empty relative Morse boundary. Note there is a related open question from \cite[Question 1.5]{CCGSZ}.

\begin{question}
    Recall that Cashen and Mackay introduced another topology on the Morse boundary that is metrizable and thus coarser than the direct limit topology \cite{cashen2019metrizable}. It is natural to ask: Do Morse hierarchies exist for the Cashen--Mackay topology?
\end{question}

In \cref{appendix} of this paper, the second author and He address this question and give a partial positive answer. \newpage

\appendix 

\section{Cashen--Mackay Morse hierarchies} \label{appendix} 
\begin{center}Annette Karrer, Vivian He \end{center} 
\vspace{0.5 cm}
As explained in \cref{sec:AlternativeMorse}, there are several possible variations of definitions of Morse hierarchies and their consequences. One variant arises by equipping the Morse boundary not with the direct limit topology introduced by Charney and Cordes \cite{Charney-Sultan, Cordes}, but with the metrizable topology introduced by Chashen--Mackay instead \cite{cashen2019metrizable}. In this appendix, we show that, with an additional quasi-convexity condition on the edge groups (see \cref{convexitycondition} in  Section \ref{appsubsec:key-lemma} ), we can transfer the machinery of Morse hierarchies to the Cashen--Mackay topology (see Section \ref{mainsec}); 
we will call these \emph{Cashen--Mackay Morse hierarchies}. 

The quasi-convexity condition we are imposing might feel very restrictive. However, many natural Morse hierarchies of groups satisfy this quasi-convexity condition. Indeed, it holds for almost all hierarchies considered in \cref{subsec:list_Morse_accessible}): hierarchies of RAAGs and Coxeter groups along special subgroups (in fact, they are convex \cite[Proposition~4.8.1]{Davis}),  prime- and JSJ-decompositions of closed oriented  non-positively curved 3-manifolds, and surface amalgams (and more generally hierarchies of certain CAT(0) groups). See \cref{Mainthm}.
 
An important fact we are using is that the direct-limit topology is finer than the Cashen--Mackay topology; so if a set in connected in the direct-limit topology, it is connected in the Cashen--Mackay topology. As a result, if the leaf groups of a Morse hierarchy have connected relative Morse boundaries in the direct-limit topology, they also have connected relative Morse boundaries in the Cashen--Mackay topology. Hence, if the edge groups of a Morse hierarchy satisfy the quasi-convexity condition, then it is also a Cashen--Mackay Morse hierarchy. 

 Most of the arguments concerning Morse hierarchies in \cref{sec:prelims},  \cref{sec:Def_Morse_hierarchy} and \cref{sec:proofThmA}   transfer naturally to the Cashen--Mackay topology. In fact, there is only one key lemma \cite[Lemma~3.5]{fioravanti2022connected} that needs to be adapted and for which we need to impose the additional quasi-convexity condition on the edge groups. This lemma says that a splitting along undistorted groups with empty relative Morse boundary yields a corresponding decomposition of the Morse boundary into clopen sets. This is the Key Lemma of this appendix (\cref{lem:key-lemma}). 

\subsection{Acknowledgments}
We thank Matthew Cordes and Kim Ruane for their helpful comments. 
\subsection{Setting and notation}

Let $X$ be a proper geodesic metric space and fix a base point $p \in X$. The metrizable topology on the Morse boundary we are interested in is defined in \cite{cashen2019metrizable}. Cashen--Mackay work with equivalence classes of quasi-geodesic rays instead of equivalence classes of geodesic rays.  As every Morse quasi-geodesic ray is at finite Hausdorff distance to a Morse geodesic ray, every such equivalence class contains a Morse geodesic ray as representative. We will rewrite the definition of the open neighborhood basis.

\begin{defn}\label{def:metrizable-topology}
Let $[\gamma] \in \mb X$ be an equivalence class of Morse quasi-geodesics, and $\gamma$ be a geodesic representative. For each $r \geq 1$, define $U([\gamma],r)$ to be the set of points $[\eta] \in \mb X$ such that for all $q \geq 1$ and $Q \geq 0$ and every continuous $(q,Q)$--quasi-geodesic ray
$\eta' \in [\eta]$ we have
\[
    d(\eta', \gamma \cap (\mc{N}_r (p))^c ) \leq \kappa(\rho_{[\gamma]},q,Q).
\]
\end{defn}
The function $\rho_{[\gamma]}$ is a sublinear function depending on the equivalence class of Morse rays, and $\kappa$ is a function defined in \cite[Definition~4.1]{cashen2019metrizable}). For more information on the metrizable topology, see \cite[Section 4]{cashen2019metrizable}. \\

Form now on,  we equip the Morse boundary with the Cashen--Mackay topology. In particular, relative Morse boundaries are now subspaces in the Morse boundary equipped with the Cashen--Mackay topology.
Recall from \cref{definition:relative_morse_boundary} that, if $X$ is a subspace of a proper metric space $X$, then $(\mb Y,X)$ denotes the relative Morse boundary of $Y$ in $X$, equipped with the subspace topology inherited from $\mb X$. Recall further from \cref{def:surviving_morse_boundary} that $\mc M(Y,X)$ denotes the surviving Morse boundary of $Y$ with respect to the embedding in $X$ which is the subspace of $\mb Y$ containing those equivalence classes of Morse rays from $Y$ that remain Morse when included in $X$. In the direct-limit topology these two topological spaces might not be homeomorphic, even if $Y$ is quasi-isometrically embedded in $X$  (compare \cref{rem_sub_top}). When working with the Cashen--Mackay topology, this discrepancy vanishes:  

\begin{thm}[Corollary 6.2  in \cite{cashen2019metrizable}]\label{lem_jotaCM}
     Let $(X, d_X)$ and $(Y, d_Y)$ be proper geodesic metric spaces and let $\iota \colon Y\to X$ be a quasi-isometric embedding. Then $\iota$ induces a homeomorphism $\iota_* \colon \mathcal M(Y,X) \to  (\mb \iota(Y), X)$.
\end{thm}

Note that this is a stronger variant of Lemma 2.3 in \cite{fioravanti2022connected}  (see \cref{lem_jota}) for the direct-limit topology that is used in  \cref{sec:prelims}, \cref{sec:Def_Morse_hierarchy} and \cref{sec:proofThmA}.\\

We now recall the setting and the notation in \cite{fioravanti2022connected}:
Let $G$ be a finitely generated group that splits as a graph of group, i.e. $G$ acts non-elliptically and minimally without inversions on a simplicial tree $T$ with edge set $\E(T)$. The edge stabilizers of this action are called \textit{edge groups} and the vertex stabilizers \textit{vertex groups}.  
A \emph{half-space} of $T$ is one of the two subtrees of $T$ obtained by removing the interior of an edge $\mf e$. If $\mf{h}$ denotes one of the half-spaces, we denote the other half-space by $\mf{h}^*$. We denote by $\H(T)$ the set of all possible half-spaces. 
 
Let $X$ be the Cayley graph $G$ with a fixed finite generating set. Let $o$ be a basepoint in the vertex set of $T$ and $f_o\colon G\ra T$ denote the orbit map $f_o(g)=go$. Every $\mf{h}\in\H(T)$ corresponds to a partition $G=f_o^{-1}(\mf{h})\sqcup f_o^{-1}(\mf{h}^*)$. We denote the stabilizer of an edge $\mf e$ in $T$ by $G_{\mf{e}}$.

The following lemma shows that, if every edge group is finitely generated and undistorted, we can give the Cayley graph a \say{tree-like} structure. More specifically, we can associate to every edge of $T$ a subspace $X(\mf e)$ in the Cayley graph $X$ that behaves like an edge in $T$: Every edge space $X(\mf e)$ partitions the Cayley graph $X$ into two disjoint sets.

\begin{lem}[{\cite[Lemma~3.2]{fioravanti2022connected}}]\label{lem:Xe-construction}
Let $G$ be a finitely generated group that splits as graph of groups so that every edge group is finitely generated and undistorted in $G$. 
For each $\mf{e}\in\E(T)$, we can choose a subgraph $X(\mf{e})\cu X$ so that the following hold:
\begin{enumerate}
\item $X(\mf{e})$ is connected, $G_{\mf{e}}$--invariant and $G_{\mf{e}}$--cocompact;
\item $X(g\mf{e})=gX(\mf{e})$ for all $g\in G$ and $\mf{e}\in\E(T)$;
\item if $\mf{h},\mf{h}^*$ are the two half-spaces determined by $\mf{e}$, then $X(\mf{e})$ contains every edge of $X$ with an endpoint in $f_o^{-1}(\mf{h})$ and the other in $f_o^{-1}(\mf{h}^*)$;
\item for every $r\geq 0$, we have $\mc{N}_{r}(f_p^{-1}(\mf{h}))\cap\mc{N}_{r}(f_p^{-1}(\mf{h}^*))\cu\mc{N}_{r}(X(\mf{e}))$, where the neighborhoods are taken in the Cayley graph $X$;
\item each edge $e\cu X$ lies in the subgraph $X(\mf{e})$ for only finitely many edges $\mf{e}\in\E(T)$.
\end{enumerate}
\end{lem}

If $\mf{h}$ is one of the two half-spaces determined by $\mf{e}$,  $X(\mf{h})\cu X$ denotes  the set $f_o^{-1}(\mf{h})\setminus X(\mf{e})\cu G$ along with all (half-open) edges of $X\setminus X(\mf{e})$ that it intersects. 
This yields a $G_{\mf{e}}$--invariant partition of the Cayley graph $X$:
\[ X = X(\mf{h}) \sqcup X(\mf{e}) \sqcup X(\mf{h}^*).\]

\subsection{Proof of Key Lemma} \label{appsubsec:key-lemma}
We now proceed to prove the key lemma. The version of the lemma for the direct limit topology is \cite[Lemma~3.5]{fioravanti2022connected}. In order to transfer it to the Cashen--Mackay topology, we impose the following additional quasi-convexity condition on the edge groups (see the third condition in \cref{lem:key-lemma} below). 

\begin{align}   \text{  $\forall~ \mf{e} \in \mc{E}(T)$, $X(\mf{e})$ in Lemma~\ref{lem:Xe-construction} is quasi-convex in $X$.} \label{convexitycondition}\end{align}
Here, we refer to the usual definition of quasi-convexity of metric spaces:
\begin{defn} 
    A subset $Y$ of a metric space $X$ is \textit{quasi-convex} if there exists a constant $D$ such that every geodesic segment with endpoints in $Y$ is contained in the the closed $D$-neighborhood $\mc N_D(Y)$ of $Y$.
\end{defn}
 It is important to note that the quasi-convexity condition we impose does depend on the choice of the generating set of $G$ defining the Cayley graph $X$, and not only on the edge groups of $G$. \\

The lemma considers the following (well-defined) subset of $\mb G$ for each $\mf{h}\in\H(T)$:
\[
M(\mf{h})=\{[\g]\in\mb G \mid \text{$\g$ is eventually contained in $X(\mf{h})$}\}.
\]

\begin{lem}[key lemma for the Cashen--Mackay topology]
\label{lem:key-lemma}
Let $G$ be a finitely generated group with Cayley graph $X$, and $\mb G$ its Morse boundary equipped with the Cashen--Mackay topology.
Suppose the action of $G$  on a Bass--Serre tree $T$ (as in \cref{conv}) results in a graph of groups splitting such that 
 for each edge $\mf{e} \in \mc{E}(T)$,
\begin{itemize}
    \item the edge group $G_\mf{e}$ is finitely generated and undistorted in $G$,
    \item the relative Morse boundary $(\mb G_\mf{e},G)$ is empty, and,
    \item the graph $X(\mf{e})$ defined in Lemma~\ref{lem:Xe-construction} is quasi-convex in $X$.
\end{itemize}
Then for every $\mf{h}\in\H(T)$, we have a partition $\partial_MG=M(\mf{h})\sqcup M(\mf{h}^*)$ into closed subsets.
\end{lem}
\begin{proof}
Let $\mf{h} \in \mc{H}(T)$ and let $\mf{e} \in \mc{E}(T)$ be the corresponding edge. Fix a basepoint $p \in X(\mf h)$.
It follows from Lemma~3.4 in \cite{fioravanti2022connected} that $M(\mf{h})$ and $M(\mf{h}^*)$ are disjoint and cover $\mb G$. It remains to show that $M(\mf{h})$ is closed in the metrizable topology.
Let $([\gamma_i])_i$ be a sequence of equivalence classes of Morse geodesic rays $\gamma_i$, with $\gamma_i(0)=p$ for all $i$. Suppose that 
    \[
    \lim_{i \to \infty}([\gamma_i])_i = [\gamma]
    \]
    for some Morse geodesic ray $\gamma$ with $\gamma(0)=p$ (in the metrizable topology). 
    Suppose that $[\gamma_i] \in M(\mf h)$ for all $i \in \mathbb N$, i.e. every $\gamma_i$ is eventually contained in $X(\mf h)$. We have to show that $ [\gamma]\in M(\mf h)$.

For each geodesic $\gamma_i$ in the sequence, consider the intersection $\gamma_i \cap X(\mf{h}^*)$. Since $\gamma_i(0)=p \in X(\mf{h})$, $\gamma_i \in M(\mf{h})$,  and every path connecting a vertex in $X(\mf{h})$ and $X(\mf{h^*})$ passes through $X(\mf e)$, each connected component of $\gamma_i \cap X(\mf{h}^*)$ is a geodesic segment whose endpoints are in $X(\mf{e})$. Hence, by quasi-convexity of $X(\mf{e})$, there exists $D>0$ such that 
\[
    \gamma_i \subseteq X(\mf h) \cup \mc{N}_D( X(\mf e)) \qquad \text{for all $i \in \mathbb N$}.
\]
 As $[\gamma_i]$ converge to $[\gamma]$, by definition of the topology, for any $R>0$, there exists $N_R>0$ such that for all $i \ge N_R$, 
\[
    d(\gamma_i,\gamma(R,\infty)) \leq \kappa(\rho_{[\gamma]},1,0).
\]
Note that $\kappa(\rho_{[\gamma]},1,0)$ is a constant. Hence, for all $i \ge N_R$, there exists $t^i_R>R$ with
\begin{equation}\label{eq:gamma-keeps-returning}
    \gamma(t^i_R)\in X(\mf h) \cup \mc{N}_{D'}( X(\mf{e})) \qquad \text{where $D' = D + \kappa(\rho_{[\gamma]},1,0)$.}
\end{equation}
We will show that there exists $R_0>0$ so that $\gamma(t^i_R) \in X(\mf h)$ for all $R \ge R_0$ and $i \ge N_R$.
Then a subsequence $(s_j)_j$ of $\left(t^{N_{R_0+i}}_{R_0+i}\right)_{i \in \mathbb N}$ is an increasing sequence with $\gamma(s_j) \in X(\mf h)$. Thus, $\gamma$ is not eventually contained in $X(\mf h^*)$. Since $\mb G = M(\mf{h}) \sqcup M(\mf{h}^*)$, this implies that $[\gamma] \in M(\mf{h})$. \\

It remains to show that there exists $R_0>0$ so that 
\[
\gamma \left(t^i_R\right) \in X(\mf h) \qquad \text{for all $R \ge R_0$ and $i \ge N_R$.}
\]
Suppose not, then there exists an increasing sequence $(t_i)_i$ so that $\gamma(t_i) \notin  X (\mf h)$. Then \eqref{eq:gamma-keeps-returning} implies that $\gamma(t_i)\in  \mc{N}_D( X(\mf e))$. Let $x$ and $y$ two points in $X(\mf e)$ with 
\[
d(x,\gamma(t_i))\le D' \qquad \text{and} \qquad(y,\gamma(t_i))\le D'.
\]
As in Lemma 3.3 in \cite{fioravanti2022connected}, one argues that there is a quasi-geodesic $\alpha$ connecting $x$ and $y$ inside of $X(\mf e)$ where the quasi-isometry constants are given by the quasi-isometry embedding from the edge group into $X$. Then by \cite[Lemma~2.1.(1)]{fioravanti2022connected}, the Hausdorff distance of $\alpha$ and $\gamma([t_i, t_i+1])$ is bounded by a constant $M$ depending only on the Morse gauge of $\gamma$, the quasi-isometry constants of $\alpha$, and $D'$.
However, this means that $\gamma$ is at bounded distance to $X(\mf e)$ -- a contradiction to the assumption that $X(\mf e)$ has empty relative Morse boundary. 
\end{proof}

\begin{remark}[Direct limit versus Cashen--Mackay topology]
If one works with the direct limit topology, it is not necessary to assume quasi-convexity of $X(\mf e)$: the key distinction between the two topologies is that, if a sequence $[\gamma_i]$ converges in the direct limit topology, there must be a uniform Morse gauge of all geodesic rays $\gamma_i$, whereas in the metrizable topology a sequence of geodesics with increasing Morse gauge can converge. Because of this, we impose the quasi-convexity assumption on $X(\mf e)$ to control $\gamma_i \cap X(\mf{h}^*)$ instead. 
\end{remark}

\begin{remark}[Quasi-convexity condition for graph of spaces] \label{graph_of_spacesvariant}
    In the proof above, we work with the setup and notation in \cite{fioravanti2022connected}, to make the difference of \cref{lem:key-lemma} and \cite[Lemma~3.5]{fioravanti2022connected} clear. However, it is possible to transfer Section 3 in \cite{fioravanti2022connected} to a setting where $G$ acts geometrically on the total space $X$ of a graph of spaces $\mc G$, assuming that the edge stabilizers act geometrically on the edge spaces and have empty relative Morse boundary. In this situation, Lemma 3.5 in \cite{fioravanti2022connected} holds for the Cashen--Mackay topology under the additional assumption that the edge spaces are quasi-convex within the total space $X$. As in the case of Cayley graphs it is important to note that the quasi-convexity of the edge spaces depend on the  total space $X$, and not only on the edge groups themselves. 
\end{remark}

\subsection{Cashen--Mackay Morse hierarchies} \label{mainsec}
In this subsection, we explain how to transfer \cref{sec:prelims}, \cref{sec:Def_Morse_hierarchy} and \cref{sec:proofThmA} to the Cashen-Mackay topology, and show that almost all Morse hierarchies in \cref{sec:examples} are Cashen--Mackay Morse hierarchies. \\

Using \cref{lem:key-lemma}, we conclude the following variant of  \cite[Theorem~A, Corollary~B]{fioravanti2022connected}  for the Cashen--Mackay topology (see  \cref{thmA}) by adding the quasi-convexity condition \ref{convexitycondition} as condition (c).
\begin{thm}
\label{thm:components}
Let $G$ be a finitely generated group with Cayley graph $X$, and $\mb G$ its Morse boundary equipped with the Cashen--Mackay topology.
Suppose the action of $G$  on a Bass--Serre tree $T$ (as in \cref{conv}) results in a graph of groups splitting such that 
 for each edge $\mf{e} \in \mc{E}(T)$,
\begin{enumerate}
    \item[(a)]  the edge group $G_\mf{e}$ is finitely generated and undistorted in $G$,
    \item[(b)] the relative Morse boundary $(\mb G_\mf{e},G)$ is empty, and,
    \item[(c)]  the graph $X(\mf{e})$ defined in Lemma~\ref{lem:Xe-construction} is quasi-convex in $X$. 
\end{enumerate}
Let $x$ and $y$ be distinct vertices of $T$ and $V_x$ and $V_y$ be their stabilizers. 
     Then $(\mb V_x, G)\cap (\mb V_y, G)=\emptyset$. 
     Furthermore, if a connected component $C \subseteq \mb G$ contains at least two points, then there exists a vertex group $V$ so that  $C$ is contained in $(\mb V, G)$.
\end{thm}

\begin{proof}[{Proof of Theorem~\ref{thm:components}}]
Follow section 3 of \cite{fioravanti2022connected}, replacing Lemma 3.5  with 
 \cref{lem:key-lemma}. 
\end{proof}

\begin{cor}
Under the assumptions of Theorem~\ref{thm:components}, suppose additionally that $(\mb V, G)$ is totally disconnected for every vertex group $V\leq G$. Then $\mb G$ is totally disconnected.
\end{cor} 

\begin{remark}[Topology on the relative Morse boundary] 
The notation of the relative Morse boundary $(\mb H, G)$ of an undistorted subgroup $H$ in a finitely generated group $G$ differs from the convention adopted in \cite{fioravanti2022connected}, where $(\mb H,G)$ denotes what we call the \textit{surviving Morse boundary of $H$ with respect to its embedding in $G$} (compare \cref{rem:differentConventions}). However, when working with the Cashen--Mackay topology, these spaces are homeomorphic, and hence, the difference of notation is immaterial (see \cref{lem_jotaCM}).
\end{remark}

The remaining statements in  \cref{sec:prelims} can be adjusted in the same manner as \cref{thm:components}, and lead to the following definition of Cashen--Mackay Morse hierarchies via adding the quasi-convexity condition  \ref{convexitycondition} to \cref{def_Morsehierarchy}:

\begin{defn}[Cashen--Mackay Morse hierarchies] \label{def_CMMorsehierarchy}
Let $G$ be a finitely generated group, and $\mb G$ its Morse boundary equipped with the Cashen--Mackay topology.    A hierarchy $\mc H$ of a finitely generated group $G$ is called  \emph{weak Cashen--Mackay Morse} if it is finite and for each non-leaf vertex $v\in \mathcal V(\mathcal H)$, the action of the associated group $G_v$  on one of its Cayley graphs $X_v$ and on a Bass--Serre tree $T$ (as in \cref{conv})  results in a graph of groups splitting such that 
 for each edge $\mf{e} \in \mc{E}(T)$,
\begin{enumerate}
    \item[(a)]  the edge group $G_\mf{e}$ is finitely generated and undistorted in $G_v$,
    \item[(b)] the relative Morse boundary $(\mb G_\mf{e},G)$ is empty, and
    \item[(c)]  the graph $X(\mf{e})$ defined in Lemma~\ref{lem:Xe-construction} is quasi-convex in $X_v$. 
\end{enumerate}

 A weak Cashen--Mackay Morse hierarchy is called a \emph{Cashen--Mackay Morse hierarchy} if for each leaf $v$ of $\mathcal H$, $(\mb G_v,G)$ is a connected subspace of $\mb G$.

\end{defn}
Using the transferred results of \cref{sec:prelims} for the Cashen--Mackay topology, the statements and proofs in \cref{sec:Def_Morse_hierarchy} and \cref{sec:proofThmA} transfer to the Cashen--Mackay topology by replacing  \cite[Lem 2.3]{fioravanti2022connected}  (see \cref{lem_jota})   with \cite[Cor 6.2]{cashen2019metrizable} (see \cref{lem_jotaCM}). We obtain:

\begin{thm}[\cref{thmMorseaccessiblegroups_Intro} for Cashen--Mackay Morse hierarchies] \label{CMMorsethm}
    Let $G$ be a finitely generated group, and $\mb G$ its Morse boundary equipped with the Cashen--Mackay topoogy. If $G$ has a Cashen--Mackay Morse hierarchy $\mc H$, then the following two sets are identical: 
    \begin{itemize}
        \item $\mathfrak{L}$, the set of leaf groups $L$ of $\mc H$ with $(\mb L, G)\neq \emptyset$
        \item $\mathfrak{S}$, the set of stabilizers of infinite connected components of $\mb G$. 
        \end{itemize}
        Moreover, every infinite connected component $C$ in $\mb G$ satisfies \[C = (\mb \stab (C), G).\]
  \end{thm} 

  \begin{cor}[\cref{corMorseaccessiblegroups_Intro} for Cashen--Mackay Morse hierarchies] 
         Let $G$ be a group with a Cashen--Mackay Morse hierarchy $\mc H$. The Morse boundary of $G$ equipped with the Cashen--Macay topology is totally disconnected if and only if for every leaf group $L$ of $\mc H$ we have $(\mb L, G)=\emptyset$. 
    \end{cor}

\begin{remark}
There are two main advantages when working with the Cashen--Mackay topology rather than the direct-limit topology. Firstly, the metric is metrizable. Secondly, each infinite connected component of the Morse boundary topologically embed into the Morse boundary of its stabilizer \cite[Thm]{cashen2019metrizable} (see \cref{lem_jotaCM}). An inconvenience as in \cref{ex:topologies different} does not occur in the Cashen-Mackay topology. On the other hand, imposing the quasi-convexity condition makes all statements more technical, and adds the need of working with a specific model space rather than only with the group itself. 
\end{remark}

 We finish this section by proving that almost all Morse hierarchies constructed in \cref{sec:examples} are Cashen--Mackay Morse hierarchies.
 
\begin{cor} \label{Mainthm}
If $G$ satisfies one of the following, then $G$ has a Cashen--Mackay Morse hierarchy. 
\begin{enumerate}
    \item $G$ is hyperbolic or has connected Morse boundary,
    \item $G$ is the free product of finitely many finitely generated groups, each of which has empty or connected relative Morse boundary,
     \item $G$ is a CAT(0) group as in \cref{ex_BH},
    \item $G$ is the fundamental group of a closed orientable 3-manifold  that carries a Riemannian metric of nonpositive sectional curvature,
    \item $G$ is a RAAG,
    \item $G$ is a Coxeter group whose Coxeter graph does not contain a spherical nexus.
   
\end{enumerate}
\end{cor}
\begin{proof}
Each of the listed groups has a Morse hierarchy descried in  \cref{subsec:list_Morse_accessible}. 
In order to show that these are also Cashen--Mackay Morse hierarchies, it suffices to show that all the listed examples satisfy the additional quasi-convexity condition  \ref{convexitycondition}. This is clearly true for (1) and (2). For (6), it follows from \cite[Proposition~4.8.1]{Davis}. For the remaining cases (3), (4), and (5), the groups that occur in the Morse hierarchies we construct in \cref{subsec:list_Morse_accessible} are fundamental groups of total spaces of graph of spaces.  By \cref{graph_of_spacesvariant}, it is enough to show that the edge spaces of these total spaces are quasi-convex. For (3), the edge spaces are convex by construction. The hierarchies in (6) can be view as a combination of (3) and (2), and hence satisfy the quasi-convexity conditon we need. Finally, for (4), the quasi-convexity condition is given by\cite[Lem 2.1]{Leeb1995}.
\end{proof}

\begin{remark}
    The listed groups in \cref{Mainthm} contain all examples of groups in \cref{sec:examples} except for 3-manifolds that cannot be equipped with a metric of non-positive curvature. The key property we need a 3-manifold $M$ to satisfy is is that $\tilde M$ can be equipped with a metric in which the lift of the tori in the JSJ decomposition are quasi-convex. A proof that this can be done for non-positively curved manifolds was given by Leeb~\cite[Lem 2.1]{Leeb1995}.
\end{remark}

\bibliographystyle{alpha}
\bibliography{papers}

@article{Leeb1995,
  author    = {Leeb, Bernhard},
  title     = {3-manifolds with(out) metrics of nonpositive curvature},
  journal   = {Inventiones mathematicae},
  year      = {1995},
  volume    = {122},
  number    = {1},
  pages     = {277--289},
  doi       = {10.1007/BF01231445},
  url       = {https://doi.org/10.1007/BF01231445}
}

@book {Serre,
    AUTHOR = {Serre, Jean-Pierre},
     TITLE = {Trees},
      NOTE = {Translated from the French by John Stillwell},
 PUBLISHER = {Springer-Verlag, Berlin-New York},
      YEAR = {1980},
     PAGES = {ix+142},
      ISBN = {3-540-10103-9},
   MRCLASS = {20H10 (05C05 22E50)},
  MRNUMBER = {607504},
}

@article {karrer2024subgroups,
    AUTHOR = {Karrer, Annette and Miraftab, Babak and Zbinden, Stefanie},
     TITLE = {Subgroups arising from connected components in the {M}orse
              boundary},
   JOURNAL = {Math. Z.},
  FJOURNAL = {Mathematische Zeitschrift},
    VOLUME = {313},
      YEAR = {2026},
    NUMBER = {2},
     PAGES = {Paper No. 25},
      ISSN = {0025-5874,1432-1823},
   MRCLASS = {99-06},
  MRNUMBER = {5071803},
       DOI = {10.1007/s00209-026-04019-2},
       URL = {https://doi.org/10.1007/s00209-026-04019-2},
}

@article {Bowditch_peripheral_splitting,
    AUTHOR = {Bowditch, B. H.},
     TITLE = {Peripheral splittings of groups},
   JOURNAL = {Trans. Amer. Math. Soc.},
  FJOURNAL = {Transactions of the American Mathematical Society},
    VOLUME = {353},
      YEAR = {2001},
    NUMBER = {10},
     PAGES = {4057--4082},
      ISSN = {0002-9947,1088-6850},
   MRCLASS = {20F67 (05C25 20F65)},
  MRNUMBER = {1837220},
       DOI = {10.1090/S0002-9947-01-02835-5},
       URL = {https://doi.org/10.1090/S0002-9947-01-02835-5},
}

@article {Bowditch_canonical_splittings,
    AUTHOR = {Bowditch, Brian H.},
     TITLE = {Cut points and canonical splittings of hyperbolic groups},
   JOURNAL = {Acta Math.},
  FJOURNAL = {Acta Mathematica},
    VOLUME = {180},
      YEAR = {1998},
    NUMBER = {2},
     PAGES = {145--186},
      ISSN = {0001-5962,1871-2509},
   MRCLASS = {20F32 (20E06 20E08 57M07)},
  MRNUMBER = {1638764},
MRREVIEWER = {Eric\ M.\ Freden},
       DOI = {10.1007/BF02392898},
       URL = {https://doi.org/10.1007/BF02392898},
}

@article {Haulmark_Hruska_canonical_splittings,
    AUTHOR = {Haulmark, Matthew and Hruska, G. Christopher},
     TITLE = {On canonical splittings of relatively hyperbolic groups},
   JOURNAL = {Israel J. Math.},
  FJOURNAL = {Israel Journal of Mathematics},
    VOLUME = {258},
      YEAR = {2023},
    NUMBER = {1},
     PAGES = {249--286},
      ISSN = {0021-2172,1565-8511},
   MRCLASS = {20F67 (20F65)},
  MRNUMBER = {4682936},
MRREVIEWER = {Michael\ L.\ Mihalik},
       DOI = {10.1007/s11856-023-2472-1},
       URL = {https://doi.org/10.1007/s11856-023-2472-1},
}

@article {Haulmark_cutpoints,
    AUTHOR = {Haulmark, Matthew},
     TITLE = {Local cut points and splittings of relatively hyperbolic
              groups},
   JOURNAL = {Algebr. Geom. Topol.},
  FJOURNAL = {Algebraic \& Geometric Topology},
    VOLUME = {19},
      YEAR = {2019},
    NUMBER = {6},
     PAGES = {2795--2836},
      ISSN = {1472-2747,1472-2739},
   MRCLASS = {20F65 (20F67)},
  MRNUMBER = {4023329},
MRREVIEWER = {Loreno\ Heer},
       DOI = {10.2140/agt.2019.19.2795},
       URL = {https://doi.org/10.2140/agt.2019.19.2795},
}

@article {Papasoglu_Swenson,
    AUTHOR = {Papasoglu, Panos and Swenson, Eric},
     TITLE = {Boundaries and {JSJ} decompositions of {CAT}(0)-groups},
   JOURNAL = {Geom. Funct. Anal.},
  FJOURNAL = {Geometric and Functional Analysis},
    VOLUME = {19},
      YEAR = {2009},
    NUMBER = {2},
     PAGES = {559--590},
      ISSN = {1016-443X,1420-8970},
   MRCLASS = {20F67 (20E06 57M07)},
  MRNUMBER = {2545250},
MRREVIEWER = {Vassilis\ Metaftsis},
       DOI = {10.1007/s00039-009-0012-8},
       URL = {https://doi.org/10.1007/s00039-009-0012-8},
}

@article {MihalikCamp,
    AUTHOR = {Camp, Wes and Mihalik, Michael},
     TITLE = {A classification of right-angled {C}oxeter groups with no
              3-flats and locally connected boundary},
   JOURNAL = {J. Group Theory},
  FJOURNAL = {Journal of Group Theory},
    VOLUME = {17},
      YEAR = {2014},
    NUMBER = {5},
     PAGES = {717--755},
      ISSN = {1433-5883,1435-4446},
   MRCLASS = {20F55 (20F65)},
  MRNUMBER = {3259092},
MRREVIEWER = {Stefan\ Witzel},
       DOI = {10.1515/jgt-2014-0019},
       URL = {https://doi.org/10.1515/jgt-2014-0019},
}

@article {hierarchically,
    AUTHOR = {Russell, Jacob and Spriano, Davide and Tran, Hung Cong},
     TITLE = {Convexity in hierarchically hyperbolic spaces},
   JOURNAL = {Algebr. Geom. Topol.},
  FJOURNAL = {Algebraic \& Geometric Topology},
    VOLUME = {23},
      YEAR = {2023},
    NUMBER = {3},
     PAGES = {1167--1248},
      ISSN = {1472-2747,1472-2739},
   MRCLASS = {20F65 (20F67)},
  MRNUMBER = {4598806},
MRREVIEWER = {Bruno\ P.\ Zimmermann},
       DOI = {10.2140/agt.2023.23.1167},
       URL = {https://doi.org/10.2140/agt.2023.23.1167},
}

@article {Tukia,
    AUTHOR = {Tukia, Pekka},
     TITLE = {Homeomorphic conjugates of {F}uchsian groups},
   JOURNAL = {J. Reine Angew. Math.},
  FJOURNAL = {Journal f\"ur die Reine und Angewandte Mathematik. [Crelle's
              Journal]},
    VOLUME = {391},
      YEAR = {1988},
     PAGES = {1--54},
      ISSN = {0075-4102,1435-5345},
   MRCLASS = {30C60 (11F06 30F35 57S05)},
  MRNUMBER = {961162},
MRREVIEWER = {Gaven\ J.\ Martin},
       DOI = {10.1515/crll.1988.391.1},
       URL = {https://doi.org/10.1515/crll.1988.391.1},
}

@article {Gabai,
    AUTHOR = {Gabai, David},
     TITLE = {Convergence groups are {F}uchsian groups},
   JOURNAL = {Ann. of Math. (2)},
  FJOURNAL = {Annals of Mathematics. Second Series},
    VOLUME = {136},
      YEAR = {1992},
    NUMBER = {3},
     PAGES = {447--510},
      ISSN = {0003-486X,1939-8980},
   MRCLASS = {20H10 (57M07)},
  MRNUMBER = {1189862},
MRREVIEWER = {Alan\ W.\ Reid},
       DOI = {10.2307/2946597},
       URL = {https://doi.org/10.2307/2946597},
}

@article {CassonJungreis,
    AUTHOR = {Casson, Andrew and Jungreis, Douglas},
     TITLE = {Convergence groups and {S}eifert fibered {$3$}-manifolds},
   JOURNAL = {Invent. Math.},
  FJOURNAL = {Inventiones Mathematicae},
    VOLUME = {118},
      YEAR = {1994},
    NUMBER = {3},
     PAGES = {441--456},
      ISSN = {0020-9910,1432-1297},
   MRCLASS = {57M50 (57M60 57N10)},
  MRNUMBER = {1296353},
MRREVIEWER = {G.\ Peter\ Scott},
       DOI = {10.1007/BF01231540},
       URL = {https://doi.org/10.1007/BF01231540},
}

@article {Genevois,
    AUTHOR = {Genevois, Anthony},
     TITLE = {Hyperbolicities in {${\rm CAT}(0)$} cube complexes},
   JOURNAL = {Enseign. Math.},
  FJOURNAL = {L'Enseignement Math\'{e}matique},
    VOLUME = {65},
      YEAR = {2019},
    NUMBER = {1-2},
     PAGES = {33--100},
      ISSN = {0013-8584,2309-4672},
   MRCLASS = {20F65 (20F67 53C23)},
  MRNUMBER = {4057355},
MRREVIEWER = {Thomas\ Haettel},
       DOI = {10.4171/lem/65-1/2-2},
       URL = {https://doi.org/10.4171/lem/65-1/2-2},
}

@article {Tran,
    AUTHOR = {Tran, Hung Cong},
     TITLE = {On strongly quasiconvex subgroups},
   JOURNAL = {Geom. Topol.},
  FJOURNAL = {Geometry \& Topology},
    VOLUME = {23},
      YEAR = {2019},
    NUMBER = {3},
     PAGES = {1173--1235},
      ISSN = {1465-3060,1364-0380},
   MRCLASS = {20F67 (20F65)},
  MRNUMBER = {3956891},
MRREVIEWER = {Jack\ O.\ Button},
       DOI = {10.2140/gt.2019.23.1173},
       URL = {https://doi.org/10.2140/gt.2019.23.1173},
}

@article {Dunwoody,
    AUTHOR = {Dunwoody, M. J.},
     TITLE = {The accessibility of finitely presented groups},
   JOURNAL = {Invent. Math.},
  FJOURNAL = {Inventiones Mathematicae},
    VOLUME = {81},
      YEAR = {1985},
    NUMBER = {3},
     PAGES = {449--457},
      ISSN = {0020-9910,1432-1297},
   MRCLASS = {20E06 (54H99 57M99)},
  MRNUMBER = {807066},
MRREVIEWER = {D.\ E.\ Cohen},
       DOI = {10.1007/BF01388581},
       URL = {https://doi.org/10.1007/BF01388581},
}

@misc{CorLev25,
      title={Connectivity of Coxeter group Morse boundaries}, 
      author={Matthew Cordes and Ivan Levcovitz},
      year={2025},
      eprint={2503.14085},
      archivePrefix={arXiv},
      primaryClass={math.GR},
      url={https://arxiv.org/abs/2503.14085}, 
}

@misc{CashenChuTao,
      title={Hyperbolic models and stability recognition for Coxeter groups}, 
      author={Christopher H. Cashen and Michelle Chu and Jing Tao},
      year={2026},
      eprint={2608.17036},
      archivePrefix={arXiv},
      primaryClass={math.GR},
      url={https://arxiv.org/abs/2608.17036}, 
}

@article{GraKarLaz21,
	author = {Graeber, Marius and Karrer, Annette and Lazarovich, Nir and Stark, Emily},
	doi = {10.1016/j.topol.2021.107645},
	fjournal = {Topology and its Applications},
	issn = {0166-8641,1879-3207},
	journal = {Topology Appl.},
	mrclass = {20F65 (20F55)},
	mrnumber = {4234246},
	mrreviewer = {Hoang\ Thanh\ Nguyen},
	pages = {Paper No. 107645, 3},
	title = {Surprising circles in {M}orse boundaries of right-angled {C}oxeter groups},
	url = {https://doi.org/10.1016/j.topol.2021.107645},
	volume = {294},
	year = {2021}}

@incollection{Behrstockstablecycle,
	author = {Behrstock, Jason},
	booktitle = {Beyond hyperbolicity},
	mrclass = {20F55 (20F65 57M07)},
	mrnumber = {3966609},
	pages = {151--159},
	publisher = {Cambridge Univ. Press, Cambridge},
	series = {London Math. Soc. Lecture Note Ser.},
	title = {A counterexample to questions about boundaries, stability, and commensurability},
	volume = {454},
	year = {2019}}

@article{DanTho15,
	author = {Dani, Pallavi and Thomas, Anne},
	fjournal = {Transactions of the American Mathematical Society},
	issn = {0002-9947},
	journal = {Trans. Amer. Math. Soc.},
	mrclass = {20F65 (20F55)},
	mrnumber = {3314816},
	mrreviewer = {Kevin D. Schreve},
	number = {5},
	pages = {3549--3577},
	title = {Divergence in right-angled {C}oxeter groups},
	url = {https://doi.org/10.1090/S0002-9947-2014-06218-1},
	volume = {367},
	year = {2015}}

@article{Lev22,
	author = {Levcovitz, Ivan},
	doi = {10.1112/topo.12267},
	fjournal = {Journal of Topology},
	issn = {1753-8416,1753-8424},
	journal = {J. Topol.},
	mrclass = {20F65 (20F69)},
	mrnumber = {4584587},
	mrreviewer = {Kevin\ D.\ Schreve},
	number = {4},
	pages = {2143--2173},
	title = {Characterizing divergence and thickness in right-angled {C}oxeter groups},
	url = {https://doi.org/10.1112/topo.12267},
	volume = {15},
	year = {2022}}

@misc{RAAGEDYRACGS,
      title={RAAGedy right-angled Coxeter groups}, 
      author={Christopher H. Cashen and Pallavi Dani and Alexandra Edletzberger and Annette Karrer},
      year={2025},
      eprint={2506.16789},
      archivePrefix={arXiv},
      primaryClass={math.GR},
      url={https://arxiv.org/abs/2506.16789}, 
}

@article {Qing,
    AUTHOR = {Liu, Qing},
     TITLE = {Dynamics on the {M}orse boundary},
   JOURNAL = {Geom. Dedicata},
  FJOURNAL = {Geometriae Dedicata},
    VOLUME = {214},
      YEAR = {2021},
     PAGES = {1--20},
      ISSN = {0046-5755,1572-9168},
   MRCLASS = {53C23 (37D40 37F32)},
  MRNUMBER = {4308266},
MRREVIEWER = {Stefan\ Suhr},
       DOI = {10.1007/s10711-021-00600-7},
       URL = {https://doi.org/10.1007/s10711-021-00600-7},
}

@incollection {Davis,
    AUTHOR = {Davis, Michael W.},
     TITLE = {The geometry and topology of {C}oxeter groups},
 BOOKTITLE = {Introduction to modern mathematics},
    SERIES = {Adv. Lect. Math. (ALM)},
    VOLUME = {33},
     PAGES = {129--142},
 PUBLISHER = {Int. Press, Somerville, MA},
      YEAR = {2015},
      ISBN = {978-1-57146-305-0},
   MRCLASS = {20F55 (57M07)},
  MRNUMBER = {3445448},
MRREVIEWER = {Mikhail\ V.\ Belolipetsky},
}

@article{Zbinden,
      title={Morse boundaries of 3-manifold groups}, 
      author={Stefanie Zbinden},
      year={2024},
      eprint={2403.03939},
      JOURNAL={accepted in Transactions of the AMS},
      primaryClass={math.GR},
}

@article{karrer,
    AUTHOR = {Karrer, Annette},
     TITLE = {Right-angled {C}oxeter groups with totally disconnected
              {M}orse boundaries},
   JOURNAL = {Geom. Dedicata},
  FJOURNAL = {Geometriae Dedicata},
    VOLUME = {217},
      YEAR = {2023},
    NUMBER = {4},
     PAGES = {Paper No. 71, 40},
      ISSN = {0046-5755,1572-9168},
   MRCLASS = {20F65 (20F55 20F67 57M07 57M15)},
  MRNUMBER = {4601993},
MRREVIEWER = {Anthony\ Genevois},
       DOI = {10.1007/s10711-023-00798-8},
       URL = {https://doi.org/10.1007/s10711-023-00798-8},
}

@article {MRT07,
    AUTHOR = {Mihalik, Michael and Ruane, Kim and Tschantz, Steve},
     TITLE = {Local connectivity of right-angled {C}oxeter group boundaries},
   JOURNAL = {J. Group Theory},
  FJOURNAL = {Journal of Group Theory},
    VOLUME = {10},
      YEAR = {2007},
    NUMBER = {4},
     PAGES = {531--560},
      ISSN = {1433-5883},
   MRCLASS = {20F65 (20F55 57M05 57M15)},
  MRNUMBER = {2334758},
MRREVIEWER = {Nadia Benakli},
       DOI = {10.1515/JGT.2007.042},
       URL = {https://doi.org/10.1515/JGT.2007.042},
}

@article {Cordes,
    AUTHOR = {Cordes, Matthew},
     TITLE = {Morse boundaries of proper geodesic metric spaces},
   JOURNAL = {Groups Geom. Dyn.},
  FJOURNAL = {Groups, Geometry, and Dynamics},
    VOLUME = {11},
      YEAR = {2017},
    NUMBER = {4},
     PAGES = {1281--1306},
      ISSN = {1661-7207},
   MRCLASS = {20F65},
  MRNUMBER = {3737283},
MRREVIEWER = {Caglar Uyanik},
       DOI = {10.4171/GGD/429},
       URL = {https://doi.org/10.4171/GGD/429},
}

@article {corrigendum,
    AUTHOR = {Cordes, Matthew and Alessandro Sisto and Zbinden, Stefanie},
     TITLE = {Corrigendum to ``{M}orse boundaries of proper geodesic metric spaces''},
   JOURNAL = {Groups Geom. Dyn.},
  FJOURNAL = {Groups, Geometry, and Dynamics},
    VOLUME = {18},
      YEAR = {2024},
    NUMBER = {4},
     PAGES = {1559–1563},
}

@article {fioravanti2022connected,
title={Connected components of Morse boundaries of graphs of groups; accepted to Pac J Math}, 
      author={Elia Fioravanti and Annette Karrer},
      year={2022},
      eprint={2109.12064},
      archivePrefix={arXiv},
      primaryClass={math.GR}
}

@article{ChaCorSis23,
	author = {Charney, Ruth and Cordes, Matthew and Sisto, Alessandro},
	doi = {10.4171/ggd/669},
	fjournal = {Groups, Geometry, and Dynamics},
	issn = {1661-7207,1661-7215},
	journal = {Groups Geom. Dyn.},
	mrclass = {20F67 (20F65 57K32)},
	mrnumber = {4563334},
	mrreviewer = {Jacob\ Russell},
	number = {1},
	pages = {157--184},
	title = {Complete topological descriptions of certain {M}orse boundaries},
	url = {https://doi.org/10.4171/ggd/669},
	volume = {17},
	year = {2023}}

@article {Caprace-Fujiwara,
    AUTHOR = {Caprace, Pierre-Emmanuel and Fujiwara, Koji},
     TITLE = {Rank-one isometries of buildings and quasi-morphisms of
              {K}ac-{M}oody groups},
   JOURNAL = {Geom. Funct. Anal.},
  FJOURNAL = {Geometric and Functional Analysis},
    VOLUME = {19},
      YEAR = {2010},
    NUMBER = {5},
     PAGES = {1296--1319},
      ISSN = {1016-443X,1420-8970},
   MRCLASS = {20F65 (20E32 20E42 20G44 22E65)},
  MRNUMBER = {2585575},
MRREVIEWER = {Guy\ Rousseau},
       DOI = {10.1007/s00039-009-0042-2},
       URL = {https://doi.org/10.1007/s00039-009-0042-2},
}

@article {Charney-Sultan,
    AUTHOR = {Charney, Ruth and Sultan, Harold},
     TITLE = {Contracting boundaries of {$\rm CAT(0)$} spaces},
   JOURNAL = {J. Topol.},
  FJOURNAL = {Journal of Topology},
    VOLUME = {8},
      YEAR = {2015},
    NUMBER = {1},
     PAGES = {93--117},
      ISSN = {1753-8416},
   MRCLASS = {20F65 (20F55 20F67)},
  MRNUMBER = {3339446},
MRREVIEWER = {Nadia Benakli},
       URL = {http://doi.org/10.1112/jtopol/jtu017},
}

@article{THOMASSEN1993248,
title = {Vertex-Transitive Graphs and Accessibility},
journal = {Journal of Combinatorial Theory, Series B},
volume = {58},
number = {2},
pages = {248-268},
year = {1993},
issn = {0095-8956},
doi = {https://doi.org/10.1006/jctb.1993.1042},
url = {https://www.sciencedirect.com/science/article/pii/S0095895683710427},
author = {C. Thomassen and W. Woess}
}

@article {CapraceWalls,
    AUTHOR = {Caprace, Pierre-Emmanuel},
     TITLE = {Conjugacy of 2-spherical subgroups of {C}oxeter groups and
              parallel walls},
   JOURNAL = {Algebr. Geom. Topol.},
  FJOURNAL = {Algebraic \& Geometric Topology},
    VOLUME = {6},
      YEAR = {2006},
     PAGES = {1987--2029},
      ISSN = {1472-2747,1472-2739},
   MRCLASS = {20F65 (20F55 20F67 57M07)},
  MRNUMBER = {2263057},
MRREVIEWER = {Robert\ W.\ Bell},
       DOI = {10.2140/agt.2006.6.1987},
       URL = {https://doi.org/10.2140/agt.2006.6.1987},
}

@article{cashen2019metrizable,
  title={A metrizable topology on the contracting boundary of a group},
  author={Cashen, Christopher and Mackay, John},
  journal={Transactions of the American Mathematical Society},
  volume={372},
  number={3},
  pages={1555--1600},
  year={2019}
}

@book {BridsonHaefliger,
    AUTHOR = {Bridson, Martin R. and Haefliger, Andr\'e},
     TITLE = {Metric spaces of non-positive curvature},
    SERIES = {Grundlehren der mathematischen Wissenschaften [Fundamental
              Principles of Mathematical Sciences]},
    VOLUME = {319},
 PUBLISHER = {Springer-Verlag, Berlin},
      YEAR = {1999},
     PAGES = {xxii+643},
      ISBN = {3-540-64324-9},
   MRCLASS = {53C23 (20F65 53C70 57M07)},
  MRNUMBER = {1744486},
MRREVIEWER = {Athanase\ Papadopoulos},
       DOI = {10.1007/978-3-662-12494-9},
       URL = {https://doi.org/10.1007/978-3-662-12494-9},
}

@book {bourbaki,
    AUTHOR = {Bourbaki, N.},
     TITLE = {\'El\'ements de math\'ematique. {F}asc. {XXXIV}. {G}roupes et
              alg\`ebres de {L}ie. {C}hapitre {IV}: {G}roupes de {C}oxeter
              et syst\`emes de {T}its. {C}hapitre {V}: {G}roupes engendr\'es
              par des r\'eflexions. {C}hapitre {VI}: syst\`emes de racines},
    SERIES = {Actualit\'es Scientifiques et Industrielles [Current
              Scientific and Industrial Topics]},
    VOLUME = {No. 1337},
 PUBLISHER = {Hermann, Paris},
      YEAR = {1968},
     PAGES = {288 pp. (loose errata)},
   MRCLASS = {22.50 (17.00)},
  MRNUMBER = {240238},
MRREVIEWER = {G.\ B.\ Seligman},
}

@article {CCGSZ,
    AUTHOR = {Charney, Ruth and Cordes, Matthew and Goldsborough, Antoine
              and Sisto, Alessandro and Zbinden, Stefanie},
     TITLE = {({N}on-)existence of {C}annon-{T}hurston maps for {M}orse
              boundaries},
   JOURNAL = {Bull. Lond. Math. Soc.},
  FJOURNAL = {Bulletin of the London Mathematical Society},
    VOLUME = {57},
      YEAR = {2025},
    NUMBER = {2},
     PAGES = {463--471},
      ISSN = {0024-6093,1469-2120},
   MRCLASS = {20F65 (20E07 20F67 57M07)},
  MRNUMBER = {4861892},
MRREVIEWER = {Daniele\ Ettore\ Otera},
       DOI = {10.1112/blms.13207},
       URL = {https://doi.org/10.1112/blms.13207},
}

@article {DaniThomasJSJ,
    AUTHOR = {Dani, Pallavi and Thomas, Anne},
     TITLE = {Bowditch's {JSJ} tree and the quasi-isometry classification of
              certain {C}oxeter groups},
   JOURNAL = {J. Topol.},
  FJOURNAL = {Journal of Topology},
    VOLUME = {10},
      YEAR = {2017},
    NUMBER = {4},
     PAGES = {1066--1106},
      ISSN = {1753-8416,1753-8424},
   MRCLASS = {20F65 (05C25 20F55 57M07)},
  MRNUMBER = {3743069},
MRREVIEWER = {Matthew\ C. B. Zaremsky},
       DOI = {10.1112/topo.12033},
       URL = {https://doi.org/10.1112/topo.12033},
}

@article {CharneyCordesMurray,
    AUTHOR = {Charney, Ruth and Cordes, Matthew and Murray, Devin},
     TITLE = {Quasi-{M}obius homeomorphisms of {M}orse boundaries},
   JOURNAL = {Bull. Lond. Math. Soc.},
  FJOURNAL = {Bulletin of the London Mathematical Society},
    VOLUME = {51},
      YEAR = {2019},
    NUMBER = {3},
     PAGES = {501--515},
      ISSN = {0024-6093,1469-2120},
   MRCLASS = {20F65 (57M07)},
  MRNUMBER = {3964503},
MRREVIEWER = {Jens\ Harlander},
       DOI = {10.1112/blms.12246},
       URL = {https://doi.org/10.1112/blms.12246},
}

@article {Gitik–Mitra–Rips–Sageev,
    AUTHOR = {Gitik, Rita and Mitra, Mahan and Rips, Eliyahu and Sageev,
              Michah},
     TITLE = {Widths of subgroups},
   JOURNAL = {Trans. Amer. Math. Soc.},
  FJOURNAL = {Transactions of the American Mathematical Society},
    VOLUME = {350},
      YEAR = {1998},
    NUMBER = {1},
     PAGES = {321--329},
      ISSN = {0002-9947,1088-6850},
   MRCLASS = {20F32 (57M07)},
  MRNUMBER = {1389776},
MRREVIEWER = {Athanase\ Papadopoulos},
       DOI = {10.1090/S0002-9947-98-01792-9},
       URL = {https://doi.org/10.1090/S0002-9947-98-01792-9},
}

@article {Mj08,
    AUTHOR = {Mj, Mahan},
     TITLE = {Relative rigidity, quasiconvexity and {$C$}-complexes},
   JOURNAL = {Algebr. Geom. Topol.},
  FJOURNAL = {Algebraic \& Geometric Topology},
    VOLUME = {8},
      YEAR = {2008},
    NUMBER = {3},
     PAGES = {1691--1716},
      ISSN = {1472-2747,1472-2739},
   MRCLASS = {20F67 (52C35 57M50)},
  MRNUMBER = {2448868},
MRREVIEWER = {Vassilis\ Metaftsis},
       DOI = {10.2140/agt.2008.8.1691},
       URL = {https://doi.org/10.2140/agt.2008.8.1691},
}

@article {Zbinden-freeproducts,
    AUTHOR = {Zbinden, Stefanie},
     TITLE = {Morse boundaries of graphs of groups with finite edge groups},
   JOURNAL = {J. Group Theory},
  FJOURNAL = {Journal of Group Theory},
    VOLUME = {26},
      YEAR = {2023},
    NUMBER = {5},
     PAGES = {969--1002},
      ISSN = {1433-5883,1435-4446},
   MRCLASS = {20F65 (57M07)},
  MRNUMBER = {4634721},
MRREVIEWER = {Christopher\ H.\ Cashen},
       DOI = {10.1515/jgth-2021-0196},
       URL = {https://doi.org/10.1515/jgth-2021-0196},
}
\end{document}